\documentclass[11pt,reqno]{amsart}
 
 \makeindex

\usepackage{amssymb}

\catcode`\@=11

\long\def\@savemarbox#1#2{\global\setbox#1\vtop{\hsize\marginparwidth 
  \@parboxrestore\tiny\raggedright #2}}
\newcommand\lref[1]{\ref{#1}%
\@ifundefined{r@DisplaY #1}{}{ (#1)}}

\newcommand\fakelabel[2]{\@bsphack\if@filesw {\let\thepage\relax
   \newcommand\protect{\noexpand\noexpand\noexpand}%
\xdef\@gtempa{\write\@auxout{\string
      \newlabel{#1}{{#2}{\thepage}}}}}\@gtempa
   \if@nobreak \ifvmode\nobreak\fi\fi\fi\@esphack}

\def\SL@margintext#1{{\showlabelsetlabel{\tiny\{\SL@prlabelname{#1}\}}}}
\catcode`\@=12

\def\Empty{}
\newcommand\oplabel[1]{
  \def\OpArg{#1} \ifx \OpArg\Empty {} \else
        \label{#1}
  \fi}
\newtheorem{theoremSt}{Theorem}[section]

\newtheorem{exampleSt}[theoremSt]{Example}
\newtheorem{exerciseSt}[theoremSt]{Exercise}

\newcommand\MakeStEnv[1]{
  \newenvironment{#1}[1]{
  \begin{#1St} \oplabel{##1}%
  \global\def\CrntSt{\thetheoremSt}%
}{ 
  \end{#1St} }
  \newenvironment{#1+}[1]{
  \begin{#1St} \label{##1}%
  \label{DisplaY ##1}%
  \global\def\CrntSt{\thetheoremSt}%
  \def\Labl{##1}\ifx\Labl\Empty{} \else {\em (\Labl)\,}\fi%
}{ 
  \end{#1St} }
}
\MakeStEnv{theorem}
\MakeStEnv{corollary}
\MakeStEnv{proposition}
\MakeStEnv{lemma}
\MakeStEnv{definition}
\MakeStEnv{conjecture}
\MakeStEnv{problem}
\MakeStEnv{question}

\long\def\State#1#2#3{
\medskip\par\noindent
{\bf #1} {\rm(#2)}
{\it #3}
\par\medskip
}

\long\def\state#1#2{
\medskip\par\noindent
{\bf #1} 
{\it #2}
\par\medskip
}

\newcommand\restate[3]{
\medskip\par\noindent
{\bf #1 \ref{#2}} 
{\it #3}
\par\medskip
}

\newlength{\saveu}

\newenvironment{pf*}[1]{%
 \begin{proof}[#1]%
}{ 
 \end{proof}
}

\newcommand{\finishproof}[1]{ 
  \def\FPArg{#1}
  \ifx\FPArg\Empty
        \newcommand\FPArg{\CrntSt}  \fi
  \smallbreak\noindent\makebox[\textwidth]{\hfill\fbox{\FPArg}}
  \medbreak\noindent
}

\newcommand\CC{{\mathcal C}}

\newcommand\FF{{\mathcal F}}
\newcommand\GG{{\mathcal G}}

\newcommand\KK{{\mathcal K}}
\newcommand\LL{{\mathcal L}}
\newcommand\MM{{\mathcal M}}
\newcommand\NN{{\mathcal N}}
\newcommand\OO{{\mathcal O}}
\newcommand\PP{{\mathcal P}}
\newcommand\QQ{{\mathcal Q}}

\newcommand\SSS{{\mathcal S}}

\newcommand\UU{{\mathcal U}}
\newcommand\VV{{\mathcal V}}

\newcommand\PMF{{\PP\kern-2pt\MM\FF}}

\newcommand\PML{{\PP\kern-2pt\MM\LL}}

\newcommand\half{{\textstyle{\frac12}}}

\newcommand\Mod{\operatorname{Mod}}

\newcommand\ep{\epsilon}

\newcommand\hhat{\widehat}

\newcommand\U{{\mathbf U}}
\newcommand\union{\cup}
\newcommand\intersect{\cap}
\newcommand\bbR{{\mathord{\text{I\kern-2pt R}}}}        
\newcommand\bbH{{\mathord{\text{I\kern-2pt H}}}}        

\newcommand\Z{{\mathbb Z}}
\newcommand\R{{\mathbb R}}
\newcommand\N{{\mathbb N}}

\newcommand\Hyp{{\mathbb H}}

\newcommand\SL[1]{\text{SL}_{#1}}

\newcommand\bigrightarrow[1]{\hbox to #1{\rightarrowfill}}
\newcommand\bigleftarrow[1]{\hbox to #1{\leftarrowfill}}

\newcommand\homeo{\cong}
\newcommand\boundary{\partial}
\newcommand\semidir{\mathrel{\hbox{\vrule depth-.03ex height1.1ex\kern-0.15em$\times$}}}

\newcommand{\diam}{\operatorname{diam}}

\numberwithin{equation}{section}

\catcode`\@=11

\def\subsection{\@startsection{subsection}{2}%
  \z@{.5\linespacing\@plus.7\linespacing}{.5em}%
  {\normalfont\bfseries\centering}}

\def\section{\@startsection{section}{1}%
  \z@{.7\linespacing\@plus\linespacing}{.5\linespacing}%
  {\normalfont\large\bfseries\centering}}

\def\subsubsection{\@startsection{subsubsection}{3}%
  \z@{.5\linespacing\@plus.7\linespacing}{-.5em}%
  {\normalfont\bfseries}}

\catcode`\@=12

\newcommand{\dist}{\operatorname{dist}}

\newcommand{\fsubd}{\mathrel{{\scriptstyle\searrow}\kern-1ex^d\kern0.5ex}}
\newcommand{\bsubd}{\mathrel{{\scriptstyle\swarrow}\kern-1.6ex^d\kern0.8ex}}
\newcommand{\fsubeq}{\mathrel{\raise-.7ex\hbox{$\overset{\searrow}{=}$}}}
\newcommand{\bsubeq}{\mathrel{\raise-.7ex\hbox{$\overset{\swarrow}{=}$}}}

\newcommand{\base}{\operatorname{base}}

\newcommand{\bbar}{\overline}

\newcommand{\tsh}[1]{\left\{\kern-.9ex\left\{#1\right\}\kern-.9ex\right\}}
\newcommand{\Tsh}[2]{\tsh{#2}_{#1}}

\newcommand\Teich{{\mathcal T}}

\newcommand\interior{{\rm int}}

\newcounter{enumitemp}
\newenvironment{enumeratecontinue}{
  \setcounter{enumitemp}{\value{enumi}}
  \begin{enumerate}
  \setcounter{enumi}{\value{enumitemp}}
}
{
  \end{enumerate}
}

\newcommand{\ra}{\to}
\newcommand{\rank}{\operatorname{rank}}

\newcommand\AM{{\mathcal{M}_\omega}}

\def\MCG{\mathcal {MCG}}

\def\ulim{\lim_\omega}

\def\cone{{\rm{Cone}}_{\omega}}
\def\dist{{\rm{dist}}}
\def\co{\colon}

\newcommand{\oaen}{$\omega$-a.e.\ ~$n$}
\def\diam{{\rm{diam}}}
\newcommand\uprod[1]{\langle #1 \rangle}
\newcommand{\sd}{\operatorname{Sd}}

\def\hull{\operatorname{hull}}

\newcommand\seq[1]{\mbox{\boldmath$#1$}}
\newcommand\subseq[1]{\mbox{\boldmath$\scriptstyle#1$}}

\newcommand\suchthat{\bigm|}
\newcommand\Flat{\FF}

\def\QQomega{\QQ_\omega}

\def\aomega{a_\omega}
\def\bomega{b_\omega}

\def\iotaomega{\iota_\omega}

\renewcommand\Mod{\MCG}

\newtheorem{remarkSt}[theoremSt]{Remark}

\renewcommand\paragraph[1]{\smallskip\textbf{#1\/}}
\newcommand\inject\hookrightarrow
\newcommand\wt\widetilde

\newcommand\bdy\partial
\newcommand\F{\mathcal{F}}
\newcommand\abs[1]{\left| #1 \right|}
\newcommand\Ann{\mathcal{A}}
\DeclareMathOperator\Homeo{Homeo}
\newcommand\Teichmuller{Teichm\"uller}
\newcommand\infinity{\infty}

\newcommand\from\colon
\newcommand\<\langle
\renewcommand\>\rangle

\newcommand\MMOne{\cite{masur-minsky:complex1}}
\newcommand\MMTwo{\cite{masur-minsky:complex2}}
\newcommand\BMRank{\cite{behrstock-minsky:rank}}
\newcommand\truncate[2]{\Tsh{#2}{#1}}
\newcommand\cross\times
\DeclareMathOperator\supp{supp}
\DeclareMathOperator\open{open}
\newcommand\restrict{\bigm|}

\newcommand\essunion{\mathrel{\mathchoice
{\cup\mkern-9.45mu\raise1.4pt\hbox{$\scriptstyle\circ$}\mkern2mu}
{\cup\mkern-9.45mu\raise1.4pt\hbox{$\scriptstyle\circ$}\mkern2mu}
{\cup\mkern-10mu\raise1.4pt\hbox{$\scriptscriptstyle\circ$}\mkern2mu}
{\cup\mkern-10mu\raise.8pt\hbox{$\scriptscriptstyle\circ$}\mkern2mu}}}
\newcommand\essint{\mathrel{\mathchoice
{\cap\mkern-9.3mu\raise0.8pt\hbox{$\scriptstyle\circ$}\mkern2mu}
{\cap\mkern-9.3mu\raise0.8pt\hbox{$\scriptstyle\circ$}\mkern2mu}
{\cap\mkern-10mu\raise0.8pt\hbox{$\scriptscriptstyle\circ$}\mkern2mu}
{\cap\mkern-10.3mu\raise0.2pt\hbox{$\scriptscriptstyle\circ$}\mkern2mu}}}
\newcommand\esssubset{\mathrel{\mathchoice
{\subset\mkern-14mu\raise0.7pt\hbox{$\scriptstyle\circ$}\mkern2mu}
{\subset\mkern-14mu\raise0.7pt\hbox{$\scriptstyle\circ$}\mkern2mu}
{\subset\mkern-10mu\raise0.35pt\hbox{$\scriptscriptstyle\circ$}\mkern2mu}
{\subset\mkern-11mu\hbox{$\scriptscriptstyle\circ$}\mkern2mu}}}

\newcommand\essUnion{\operatorname{{\bigcup\mkern-17mu\circ\mkern6mu}}}

\newcommand\cC{{C}}

\newcommand\parenref[1]{(\ref{#1})}
\newcommand\pref[1]{\parenref{#1}}

\newcommand\isotopic{\stackrel{i}{\equiv}}
\newcommand\TM{\mathcal{TM}}
\newcommand\reals{\mathbf{R}}

\DeclareMathOperator\Length{Length}
\DeclareMathOperator\QI{QI}
\DeclareMathOperator\Fill{Fill}

\newcommand\ceq[1]{\stackrel{\scriptscriptstyle{#1}}{\approx}}

\title{Geometry and rigidity of mapping class groups}
\date{\today}

\author{Jason Behrstock}
\address{Lehman College, CUNY}
\email{jason@math.columbia.edu}
\author{Bruce Kleiner}
\address{Yale University}
\email{bruce.kleiner@yale.edu}
\author{Yair Minsky}
\address{Yale University}
\email{yair.minsky@yale.edu}
\author{Lee Mosher}
\address{Rutgers University at Newark}
\email{mosher@andromeda.rutgers.edu}

\thanks{Partially supported by NSF grants DMS-0812513, 
DMS-0701515,
  DMS-0504019, DMS-0706799}

\begin{document}

\begin{abstract}
    We study the large scale geometry of mapping class groups
    $\MCG(S)$, using hyperbolicity properties of curve complexes. 
    We show that any self quasi-isometry of 
    $\MCG(S)$ (outside a few sporadic cases) is a bounded distance away from a
    left-multiplication, and as a consequence obtain quasi-isometric 
    rigidity for $\MCG(S)$, namely that groups quasi-isometric to $\MCG(S)$
    are equivalent to it up to extraction of finite-index subgroups and
    quotients with finite kernel. (The latter theorem was proved by Hamenst\"adt
    using different methods). 
    
    As part of our approach we obtain several
    other structural results: a description of
    the tree-graded structure on the asymptotic cone of $\MCG(S)$;
    a characterization of the image of the curve complex projections
    map from $\MCG(S)$ to $ \prod_{Y\subseteq S} \CC(Y)$; and a construction of
    \emph{$\Sigma$-hulls} in $\MCG(S)$, an analogue of 
    convex hulls.
\end{abstract}

\maketitle
\setcounter{tocdepth}{1}
\tableofcontents
\setlength{\parskip}{1.1ex}


\section{Introduction}
\label{intro}

In this paper we investigate the coarse geometry of the mapping class
group $\MCG(S)$ of an oriented finite-type surface $S$, with our main
goal being a proof of quasi-isometric rigidity for the group. Along
the way we develop a number of tools which we hope will be of
independent interest. 

The following classification theorem for quasi-isometries applies for
all but a standard collection of exceptional cases 
(definitions are in
\S\ref{prelims}, and complete statements
of these theorems handling the exceptions are in \S\ref{endgame}.)
\begin{theorem}{QI classification}
{\rm (Classification of Quasi-Isometries)} 
If $S$ has complexity at least 2 and is not a two-holed torus,
 quasi-isometries of $\MCG(S)$ are uniformly
close to isometries induced by left-multiplication.
\end{theorem}
Note that an immediate consequence of this theorem is that, barring a few
exceptional cases, $\MCG(S)$ is isomorphic to its own quasi-isometry group,
i.e. the group of quasi isometries modulo those that are finite distance
from the identity (see Corollary \ref{QIGroup}).

Theorem \ref{QI classification} also implies the following rigidity 
property (see Section~\ref{endgame} for an alternate statement of 
this result which also covers the complexity $2$ cases):
\begin{theorem}{QI rigidity}
{\rm (Quasi-Isometric Rigidity)} If $S$ has complexity greater than 2 and 
$\Gamma$ is a finitely generated group 
quasi-isometric to $\MCG(S)$ then there is a homomorphism from $\Gamma$ 
to  $\MCG(S)$, with finite kernel and finite-index image.
\end{theorem}

Theorem \ref{QI rigidity} was proved by Hamenst\"adt in
\cite{hamenstadt:III}. The two proofs have a similar flavor
in broad outline, although the underlying machinery supporting the outline
is different. Mosher-Whyte \cite{mosher:qir-puncture} previously
established the theorem for once-punctured surfaces of genus at least 2, using quite different
methods.

The study of coarse properties of groups can be said to have started
with Milnor and \v{S}varc
\cite{milnor:growth,milnor:solvable,svarc:growth}.  This advanced
considerably with Stallings' theorem on ends of groups
\cite{stallings:ends}, the Mostow
rigidity theorem \cite{mostow:monograph}, and Gromov's theorem on
groups of polynomial growth  \cite{gromov:polynomial}. Gromov
\cite{gromov:asymptotic} proposed an ambitious program for the study
and classification of groups by their coarse geometric properties,
which continues to guide current research.   One branch of this
program
is the classification of groups up to quasi-isometry, i.e. determining
up to isomorphism (or commensurability) the groups in a quasi-isometry
class. In the last twenty years much progress has been made on
different cases of this program:

\begin{itemize}
\item Groups quasi-isometric to lattices in semi-simple Lie groups
(Gromov,
Sullivan,
Tukia,
Casson-Jungreis,
Gabai,
Hinkkanen,
Pansu,
Cannon-Cooper,
Schwartz,
Farb-Schwartz,
Kleiner-Leeb,
Eskin
\cite{CannonCooper, 
CassonJungreis,
eskin:qirigid,
farb-schwartz:hilbert,
Gabai:convergence,
Gromov:ICMAddress,
gromov:asymptotic, 
GromovPansu:lattices,
hinkkanen:quasisymmetric2, 
KleinerLeeb:buildings,
Pansu:CC,
schwartz:rankone,
sullivan:rigidity, 
Tukia:QCextension}).

\item Groups quasi-isometric to Euclidean buildings 
(Kleiner-Leeb \cite{KleinerLeeb:buildings}).
\item Groups quasi-isometric to Fuchsian buildings 
(Bourdon-Pajot, Xie \cite{BourdonPajot:HyperbolicBuildings,Xie:QIrigidFuchsian}).
\item Graphs of   $n$-dimensional Poincar\'e duality groups, for some fixed $n$
(Farb-Mosher,
Whyte,
Mosher-Sageev-Whyte,
Papasoglu
\cite{FarbMosher:BSOne,FarbMosher:BSTwo, 
MosherSageevWhyte:QTTwo, MosherSageevWhyte:QTOne,
Papasoglu:Zsplittings,whyte:higher-bs}).

\item Certain solvable semi-direct
products $\R^k\semidir \R^l$  
(Eskin-Fisher-Whyte, 
Dymarz, 
Peng 
\cite{Dymarz,EskinFisherWhyte,peng:coarsediff1}).

\end{itemize}

Research in mapping class groups has long been motivated by a drive to find
analogies with lattices. The list of such results is intriguing but
incomplete: see for example
McCarthy \cite{mccarthy:tits} and Ivanov \cite{ivanov:tits} (a Tits-alternative), Ivanov
\cite{ivanov:arithmetic} and
Farb-Masur \cite{farb-masur:superrigidity} (non-arithmeticity and
superrigidity), Harvey \cite{harvey:boundary} (analogy to Bruhat-Tits
buildings),
Harer \cite{harer:stability} (homology stability),
Farb-Lubotzky-Minsky \cite{farb-lubotzky-minsky} (rank 1 behavior)
and Andersen \cite{andersen:noT} (failure of property
(T)). The expectation has been that quasi-isometric rigidity should
hold
for these groups, but the standard tools have not previously found
purchase.

\subsection*{Structural results}
Our analysis builds on the geometry of curve complexes of surfaces
developed in
\cite{masur-minsky:complex1,masur-minsky:complex2,behrstock:thesis,behrstock-minsky:rank},
but introduces a number of new ideas and techniques, which have some
further applications (for example the
rapid-decay property for $\MCG(S)$, Behrstock-Minsky \cite{behrstock-minsky:RD}).
Among these are a theorem characterizing the
curve-complex projection image of $\MCG(S)$; a coarse convex hull
construction in $\MCG(S)$; and a theory of ``jets,'' which analyzes
infinitesimal directions in the asymptotic cone of $\MCG(S)$
in order to control separation properties.
We summarize
these below (precise statements appear later in the paper):

\subsubsection*{Consistency Theorem}
To the isotopy class of each essential subsurface $W\subseteq S$ is
associated a $\delta$-hyperbolic complex $\CC(W)$, its curve complex
(or twist complex if $W$ is an annulus). There is a
coarse-Lipschitz projection $\pi_W:\MCG(S)\to\CC(W)$, and 
combining these over all $W$ we obtain the \emph{curve complex
  projections map}\index{curve complex projections
  map}\index{1aacurvecomplexprojection@$\Pi$, curve complex
  projections map} 
$$ \Pi\from \MCG(S) \to \prod_W\CC(W).
$$
The Quasidistance Formula of
\cite{masur-minsky:complex2} (see Theorem~\ref{distance formula} below) shows 
that this map is, in a limited sense,
like a quasi-isometric embedding to the $\ell^1 $ metric on the
product.  

Behrstock established {\em consistency inequalities} (see
\S\ref{consistency}) that are satisfied by the image of $\Pi$, and in
this paper we prove that these give a coarse characterization of the image:

\restate{Theorem}{consistency suffices}{The consistency inequalities give a
  necessary and sufficient condition for a point in $\prod_W\CC(W)$ to
  be near the image of $\Pi$.}

This theorem
makes it possible to state and analyze many constructions in $\MCG(S)$
simply in terms of what they look like in the projection image, where
the hyperbolicity of the $\CC(W)$'s  can
be exploited.

\subsubsection*{$\Sigma$-hulls}
An important construction that follows from the Consistency Theorem is
a coarse substitute for a convex hull in $\MCG(S)$.  In a
$\delta$-hyperbolic space, 
the union of all geodesics joining
a set of points is quasi-convex; let us call this simply the ``hull''
of the set. 
For a finite subset $A\subset\MCG(S)$ (and suitable $\ep>0$) we describe a set
$\Sigma_\ep(A)\subset\MCG(S)$, which is exactly the set of points that 
project in each factor of 
$\prod_W\CC(W)$ to within $\ep$ of the hull of the image of $A$.
This notion of hull satisfies a number of useful properties. 

\state{Proposition \ref{hull retraction} and Lemma \ref{coarse Sigma properties}.}{%
Fixing the cardinality of $A$, $\Sigma_\ep(A)$
depends in a coarse-Lipschitz way on $A$ with respect to  the
Hausdorff metric. Its diameter is controlled by $\diam(A)$, and it
admits a coarse retraction $\MCG(S)\to\Sigma_\ep(A)$, which itself
has Lipschitz dependence on $A$. 
}

These hulls, and their rescaling limits in the asymptotic
cone, give us a way to build and control singular chains in the cone,
which we use a la Kleiner-Leeb \cite{KleinerLeeb:buildings} to describe top-dimensional
flats using local homology arguments. 

A different application of $\Sigma$-hulls will appear in
Behrstock-Minsky \cite{behrstock-minsky:RD} where $\MCG(S)$ is
shown to have the rapid decay property of Haagerup/Jolissaint.

\subsubsection*{Jets and separation}
In Behrstock-Minsky \cite{behrstock-minsky:rank}, it was shown that
any pair of distinct points in the asymptotic cone of $\MCG(S)$ could
be separated by an ultralimit of subsets associated to mapping class
groups of subsurfaces. This enabled an inductive argument to
compute the compact topological dimension of the cone.  Here we refine
our understanding of the components of the complement of such a set,
introducing the notion of a {\em jet}. A jet is a sequence of 
geodesic segments in curve complexes of subsurfaces of $S$ (modulo an ultrafilter), such that
asymptotic behavior of projections to these segments determines the
division of the complement into connected components. In the outline below we
motivate this notion via an analogy to $\delta$-hyperbolic spaces. 

As an application of these ideas we also pause in Section~\ref{pieces}
to give a characterization of the tree-graded structure (a la Dru\c
tu-Sapir \cite{drutu-sapir:treegraded}) of the cone of $\MCG(S)$.
This characterization has been used in the recent proof that any
finitely generated group with Property~(T) has at most finitely many
pairwise non-conjugate homomorphisms into a mapping class group
(Behrstock-Dru\c tu-Sapir \cite{BehrstockDrutuSapir:MCGsubgroups}).

\subsection*{Outline of the proof}

In this section we discuss the broad structure of the paper and the
proof of Theorem \ref{QI classification} (the proof of 
Theorem \ref{QI rigidity} from \ref{QI classification} is fairly
standard and will be given in Section \ref{endgame}). 

Definitions and other preliminaries will be given in Section
\ref{prelims}.

\subsubsection*{Coarse preservation of Dehn twist flats}

In order to control a quasi-isometry $f\from G\to G$ of any group, we
wish to identify structures in $G$ which are robust enough to be
preserved by $f$, and intricate enough that they can only be preserved
in the obvious ways.  In the case of $\MCG(S)$, these structures are
{\em (maximal) Dehn twist flats}, which are cosets of maximal-rank 
free-abelian subgroups generated by Dehn twists.  Later we will
actually work with equivalent subsets of the {\em marking complex}
$\MM(S)$, which is our geometric model of choice for $\MCG(S)$. 

Theorem \ref{coarse flat theorem} states that a
quasi-isometry $f\from \MCG(S)\to\MCG(S)$ {\em coarsely preserves} the
set of Dehn twist flats.  That is, the image of any such flat is
within finite Hausdorff distance of another such flat, with the 
bound depending only on the quality of the quasi-isometry.

Once Theorem \ref{coarse flat theorem} is established, we can apply
known results to prove Theorem~\ref{QI classification}.  The coarse
permutation of flats induces an automorphism of the curve complex of
$S$ which by a theorem of Ivanov--Korkmaz--Luo 
\cite{ivanov:complexes2,korkmaz:complex,luo:complex}
is induced by
some mapping class $\phi \in \MCG(S)$ (when $S=S_{1,2}$ this is not
quite right but we ignore this for now).
This gives us the desired element of $\MCG(S)$ and it is then not hard
to show that left-multiplication by $\phi$ is uniformly close to the
quasi-isometry $f$.

\subsubsection*{Preservation of asymptotic Dehn twist flats}

Theorem \ref{coarse flat theorem} will be proven, following
Kleiner-Leeb \cite{KleinerLeeb:buildings}, by reduction to the 
{\em  asymptotic cone} of $\MCG(S)$.  The asymptotic cone, which we denote
$\MM_\omega $, is a limit obtained by unbounded rescaling of the word
metric on $\MCG(S)$.  Extracting this limit requires the choice of an
ultrafilter, although our results hold for any choice --- see \S
\ref{cone defs} for details.

A quasi-isometry of $\MCG(S)$ converges after rescaling to a
bilipschitz homeomorphism of $\MM_\omega $, and Dehn twist flats limit
to bilipschitz-embedded copies of Euclidean space.  Thus our goal in
this context is to show that these asymptotic Dehn twist flats are
permuted by the limiting map:

\state{Theorem \ref{asymptotic flat theorem}.}{%
Except when $S$ is a one-holed torus or four-holed sphere, 
any homeomorphism $f\co\AM(S) \to \AM(S)$ 
 permutes the 
Dehn twist flats in $\AM(S)$. 
}

The proof of  Theorem
\ref{asymptotic flat theorem} will take up most of the
paper.

The reduction of Theorem \ref{coarse flat theorem} to Theorem
\ref{asymptotic flat theorem} will be discussed in detail in
Section~\ref{asymptotic to coarse}.  Let us now sketch the proof of
Theorem~\ref{asymptotic flat theorem}.

\subsubsection*{Structure of $\MCG(S)$ via hyperbolicity}
We begin in Sections \ref{SectionCubes}--\ref{hulls} by refining the
tools developed in Masur--Minsky
\cite{masur-minsky:complex1,masur-minsky:complex2} and Behrstock
\cite{behrstock:asymptotic} to study the coarse structure of $\MCG(S)$
using properties of curve complexes.

In Section \ref{SectionCubes}, we analyze the structure of subsets of
$\MCG(S)$ that come, in a simple way, from restrictions on some of the
coordinates of the subsurface projection maps.  In particular what we
call a {\em product region} is a set $\QQ(\mu)$ corresponding to all
markings of $S$ that contain a fixed partial marking $\mu$.  Coarsely
this is the same as a coset in $\MCG(S)$ of a stabilizer of the
partial marking.  The factors in the product structure are indexed by
the different components of the subsurface of $S$ on which the partial
marking $\mu$ is undefined. In this language a Dehn twist
flat\index{Dehn twist!flat} is defined as a set $\QQ(\mu)$ where $\mu$
is a pants decomposition whose curves are unmarked by any choice of
transversal. Choosing transversals gives one real-valued parameter for
each curve in $\mu$, making a Dehn twist flat quasi-isometric to
$\reals^{\xi(S)}$, where $\xi(S)$ is the number of components of a pants
decomposition of $S$. 

A {\em cube} is a special subset of a product region which is in fact quasi-isometric to a Euclidean cube in a way compatible with the product structure.  We show that product regions and cubes are quasi-isometrically embedded subsets of $\MCG(S)$, generalizing a result of \cite{masur-minsky:complex2}. In Lemmas~\ref{juncture} and~\ref{cube junctures} we analyze the sets, called \emph{junctures}, along which two of these regions are close --- junctures are generalizations of the coarse intersections of quasiconvex sets in a hyperbolic space.

In Section \ref{consistency}, 
we prove the Consistency Theorem \ref{consistency
  suffices}, which characterizes the image of the curve complexes projection map.
This theorem can be applied in many cases to
supplant the use of the
\emph{hierarchy paths} from
Masur--Minsky \cite{masur-minsky:complex2}.
Although useful, these paths are
technical to define and to work with. Hence, for the most part we 
have avoided using them and
Theorem~\ref{consistency suffices} is one of the main tools that 
allows us to do this. 

In Section \ref{hulls}, we use Theorem \ref{consistency suffices} to
define and study $\Sigma$-hulls.

\subsubsection*{Local homology via hulls}
In Section \ref{contractibility} we
use $\Sigma$-hulls in order to study the local homology properties of
the asymptotic cone.  The coarse properties established in
Proposition \ref{hull retraction} and Lemma 
\ref{coarse Sigma properties} imply, in the
cone, that ultralimits of $\Sigma$-hulls are contractible and have
controlled geometry (Lemmas \ref{Sigma properties} and \ref{hull
contractibility}). This allows us to use them to build singular
chains which \emph{$\Sigma$-compatible}, meaning that each simplex is
contained in the $\Sigma$-hull of its 0-skeleton.  With this, and the results of
Behrstock-Minsky \cite{behrstock-minsky:rank}, we prove local
homological dimension bounds, a result originally established by 
Hamenst\"adt \cite{hamenstadt:III}. We also obtain Corollary \ref{chain control},
an analogue of a result of Kleiner-Leeb \cite{KleinerLeeb:buildings},
which controls embedded top-dimensional manifolds in the cone.
$\Sigma$-compatible chains will also be crucial in Sections
\ref{separation} and \ref{finiteness2}.

\subsubsection*{Separation via jets}
In Section \ref{separation}, we refine the results of
\cite{behrstock:asymptotic} and \cite{behrstock-minsky:rank} to analyze
separation properties of the asymptotic cone. 

One can consider the $\delta$-hyperbolic case in order to describe the
basic intuition behind these separation arguments.  If $X$ is
$\delta$-hyperbolic and $\{g_n\}$ is a sequence of geodesic segments
of lengths going to infinity, we obtain an ultralimit $g_\omega $ in
the asymptotic cone $X_\omega $, which may be a geodesic segment or a
point (assume the latter, for simplicity).  Nearest-point projection
$\pi_n\from X\to g_n$ yields a relation on $X_\omega \setminus
g_\omega $: say that two ultralimits $x_\omega $ and $y_\omega $ are
equivalent if for representative sequences $x_n$ and $y_n$, the
sequence of distances $d(\pi_n(x_n),\pi_n(y_n))$ is $\omega$-almost
everywhere bounded.

It is a nice exercise to check directly from $\delta$-hyperbolicity
that this gives a well-defined equivalence relation, whose equivalence
classes are open, and hence these classes are separated from each
other by $g_\omega $.  In particular if $x_n$ and $y_n$ are always
projected to opposite ends of $g_n$ and neither $x_\omega $ nor
$y_\omega $ lies in $g_\omega $ then they are separated by it. When
$g_\omega$ is a point we call the sequence $\{g_n\}$ a 
{\em microscopic jet}. 

Theorem \ref{projection separation} gives an analogous statement for
$\MCG(S)$, where $\delta$-hyperbolicity of individual curve complexes
is exploited in a similar way, and separating sets are not points but
product regions.

We will also have need to think about the setting where, appealing
again to our $\delta$-hyperbolic analogy, one of our points $x_\omega
$ or $y_\omega$ may lie in $g_\omega$.  For this we introduce a finer
analysis of what we call macroscopic jets with either \emph{gradual} or
\emph{sudden growth}, and prove a suitable separation result in Theorem
\ref{Lambda acyclic}.  In this case we show that certain of the
components are acyclic, and this is where $\Sigma$-hulls come into the
proof.

Section~\ref{pieces} is a digression in which we use these ideas to
characterize the {\em pieces} of the \emph{tree-graded structure} of the
asymptotic cone of $\MCG(S)$, in the sense of Dru\c tu-Sapir
\cite{drutu-sapir:treegraded}.  Although this is not needed for the
rest of the proof, it is a structural fact which follows directly from
our techniques and is likely of independent interest.

\subsubsection*{Finiteness for manifolds in the asymptotic cone}
In Section \ref{finiteness2} we apply the foregoing results to prove a
local finiteness theorem for manifolds in the asymptotic cone.

Theorem \ref{Sigma cubes} shows that the $\Sigma$-hull of a finite set
in a top-dimensional manifold in the cone is always contained in a
finite number of cubes.  Most of the work is done in 
Theorem~\ref{simultaneous trim}, which uses the separation theorems to control
which subsurface projections of the finite set can grow without 
bound. This allows us  
to control the structure of paths connecting points in the
set which behave in an efficient way with respect to their curve
complex projections, e.g., hierarchy paths.

Finally, Theorem \ref{local cube finiteness} states that any
top-dimensional manifold is, locally at any point, contained in a finite
number of cubes.  This uses the results of
Section~\ref{contractibility} --- in particular Theorem~\ref{Sigma
cubes} and a triangulation argument allow us to approximate any sphere
in the manifold as the boundary of a chain supported in finitely many
cubes, and Corollary~\ref{chain control} implies the ball in the manifold 
bounded by the sphere is therefore contained (except for a small error
near the boundary) in these cubes as well.

\subsubsection*{Orthant complex}
In Section \ref{orthant defs} we use the finiteness theorem to study
the local structure of manifolds in the asymptotic cone,
reducing it to a combinatorial question about the {\em complex of
orthants,} which is the complex of germs of cubes with a corner at a
given point $\seq x$. The starting point, using
Theorem \ref{local cube finiteness}, is the fact that the germ of any
top-dimensional manifold at $\seq x$ is equal to a finite collection of orthants.
This allows us to characterize the structure of the complex of orthants
using purely topological properties,
and in particular (Corollary~\ref{top char of twist germs}) to
characterize the germs of Dehn-twist flats in the cone. 
This means that any 
homeomorphism of the cone must permute the germs of Dehn-twist flats.
A simple local-to-global argument
gives the proof of Theorem~\ref{asymptotic flat theorem}.

\section{Preliminaries}
\label{prelims}

In Section~\ref{SectionCC} we review the foundations of curve complexes and marking complexes. In Section~\ref{SectionProjections} we review projection maps between curve complexes of subsurfaces, and how they are used to obtain a quasidistance formula for the marking complex. The main references are \MMTwo\ and \BMRank. In Section~\ref{cone defs} we review asymptotic cones.

\subsection{Curves and markings}
\label{SectionCC}

\subsubsection{Basic definitions.} 
\label{SectionBasicDefs}
A \emph{finite type surface}\index{finite type surface} $X$ is an
oriented surface homeomorphic to a closed surface minus a finite set
of points.  The missing points are in one-to-one correspondence with
the ends of $X$, and these are referred to as the \emph{punctures}
of~$X$.  If $X$ is connected we denote $X=S_{g,n}$ where $g$ is the
genus and $n$ the number of punctures, and we quantify the 
\emph{complexity}\index{complexity} of $S_{g,n}$ by $\xi(S_{g,n}) =
\max\{3g-3+n,0\}$.\index{1aapants@$\xi(S)$, number of curves in a
pants decomposition} All surfaces in this paper will be of finite
type.

Throughout the paper we will consider a single connected ``ambient'' surface $S$ such that $\xi(S) \ge 2$ --- so $S$ is not a sphere with $\le 4$ punctures nor a torus with $\le 1$ puncture. We will also consider subsurfaces of $S$ for which $\xi\ge 1$, as well as subannuli of~$S$.

The \emph{(extended) mapping class group}\index{mapping class group} of $S$ is the group $$\MCG(S) =
\Homeo(S) / \Homeo_0(S)
$$
where $\Homeo(S)$ is the group of homeomorphisms of $S$, and $\Homeo_0(S)$ is the normal subgroup of homeomorphisms isotopic to the identity. We will often implicitly consider isotopy classes, i.e.\ $\Homeo_0(S)$-orbits, of various objects such as subsurfaces and simple closed curves. When this relation is explicit we will denote it by $\isotopic$.\index{1aaisotopic@$\isotopic$, isotopic}

\subsubsection{Curves.} \label{SubsubsectionCurves}
An \emph{essential curve}\index{essential!curve} on a finite type surface $X$ is an embedded circle~$\gamma$ such that if $X$ is not an annulus then no component of $X-\gamma$ is a disc or a once-punctured disc, and if $X$ is an annulus then $\gamma$ is a core of~$X$. An \emph{essential curve system}\index{essential!curve system} on $X$ is a nonempty collection $\cC$ of finitely many pairwise disjoint essential curves. 

If $X$ is connected then a curve system $\cC$ on $X$ is called a \emph{pants decomposition}\index{pants decomposition} if each component of $X-\cC$ is a three-punctured sphere, a \emph{pair of pants}. A nonempty pants decomposition exists on $X$ if and only if $\xi(X) \ge 1$, in which case its number of curves is $\xi(X)$.\index{1aapants@$\xi(S)$, number of curves in a pants decomposition} 

Given any finite type surface $X$, all pairwise nonisotopic, maximal, essential curve systems on $X$ have the same cardinality, a number denoted $r(X)$\index{1aarank@$r(X)$, rank} and called the \emph{rank} of $X$, equal to the sum of the $\xi$-values of the components of $X$ plus the number of annulus components. When $X$ is an essential subsurface of $S$ (see below), this number $r(X)$ equals the locally compact dimension of certain subsets of the asymptotic cone of $\MCG(S)$ associated to $X$; see Lemma~\ref{DimensionLemma}.

Given two essential curve systems $\cC,\cC'$, we may always isotope
one of them so that they are in \emph{efficient position},\index{efficient position} which means that $\cC,\cC'$ are transverse and no component of $S-(\cC \union \cC')$ has closure which is a \emph{bigon}, a nonpunctured disc whose boundary consists of an arc of $\cC$ and an arc of~$\cC'$. We say that $\cC$ and $\cC'$ \emph{overlap}\index{overlap!of two curves} if, after putting them in efficient position, the intersection is nonempty. The concept of overlap will be generalized below.

\subsubsection{The lattice of essential subsurfaces.} 
An \emph{essential subsurface}\index{essential!subsurface} of $S$ is a subsurface $Y \subset S$ with the following properties.
\begin{itemize}
\item $Y$ is a union of (not necessarily all) complementary components of an essential curve
  system $C$. Denote $C\intersect \overline Y$ by $\bdy Y$, the
  boundary curves of $Y$. 
\item No two
components of $Y$ are isotopic --- equivalently, no two annulus
components are isotopic.  
\item Each nonannulus component of $Y$ has
$\xi \ge 1$, equivalently, no component is a 3-punctured sphere.
\end{itemize}
Essential subsurfaces of $S$ are identified when they are isotopic in $S$. Note that two isotopic essential subsurfaces need not be ambient isotopic, for instance the complement of a single essential curve $c$ is isotopic to but not ambient isotopic to the complement of an annulus neighborhood of $c$.

Given an essential subsurface $X$ of $S$, let $\Gamma(X)$ denote the set of isotopy classes in $S$ of simple closed curves that are essential in $X$. Note that a boundary curve of $X$ has isotopy class in $\Gamma(X)$ if and only if it is isotopic to the core of an annulus component of $X$. Because we have excluded 3-holed spheres, $\Gamma(X)$ is empty if and only if $X$ is empty. Note that $r(X)$ equals the maximum cardinality of a curve system in $S$ whose elements are in $\Gamma(X)$.

On the set of essential subsurfaces define a relation $X \esssubset Y$,\index{1aaessentialsubset@$\esssubset$, essential subsurface} read ``$X$ is an essential subsurface of $Y$'', to mean $\Gamma(X) \subset \Gamma(Y)$.

\begin{lemma}{lattice of subsurfaces}
The relation $\esssubset$ is a partial order on the set of isotopy classes of essential subsurfaces of $S$ (including $\emptyset$). In particular, $X$ is isotopic to $Y$ if and only if $\Gamma(X)=\Gamma(Y)$. Moreover there exist binary operations $\essint,$ and $\essunion$, called \emph{essential union}\index{essential!union} and \emph{essential intersection},\index{essential!intersection} which have the following properties: 
\begin{enumerate}
\item $X\essunion Y$\index{1aaessentialunion@$\essunion$, essential union} is the unique $\esssubset$-minimal essential subsurface $Z$ such that $X\esssubset Z$ and $Y\esssubset Z$. 
\item $X\essint Y$\index{1aaessentialintersection@$\essint$, essential intersection} is the unique $\esssubset$-maximal essential subsurface $Z$ such that $Z\esssubset X$ and $Z\esssubset Y$. 
\end{enumerate}
In other words we have a lattice whose partial order is $\esssubset$ and whose meet and join operations are $\essint$ and $\essunion$,
respectively. 
\end{lemma}

\begin{proof}
To define these operations it is helpful to fix a complete hyperbolic
metric on $S$. Every essential curve has a unique geodesic
representative. Every connected essential subsurface $X$ which is not an annulus is represented by the appropriate component of the complement of the union of the geodesic representatives of $\bdy X$. Every essential subannulus is represented by the geodesic representative of its core. We call this the \emph{geodesic representative} of a connected essential subsurface. Note that disjoint components of an essential subsurface have disjoint geodesic representatives, even when annuli are involved. 

Now we can see that $\Gamma(X)$ determines $X$ as follows. If
$C\subset \Gamma(S)$, then for any finite subset of $C$ we can
take a regular neighborhood of the union of geodesic representatives,
fill in disk or punctured-disk components of the complement, and
obtain an essential subsurface. Any exhaustion of $C$ by finite sets gives an increasing sequence of such subsurfaces, which must therefore eventually stabilize up to isotopy. This uniquely determines an essential subsurface which we call $\Fill(C)$. One easily shows
$X \isotopic \Fill(\Gamma(X))$ provided $\Gamma(X)\ne\emptyset$, that is, if $X$ is not a pair of pants. It follows immediately that $\esssubset$ is a partial order.

Let us now show that $\essunion$ is defined. Given $X$ and $Y$, let
$Z=\Fill(C)$ where $C=\Gamma(X) \union \Gamma(Y)$. 
Any curve in $C$ either overlaps some other
curve in $C$, in which case it is essential in a nonannulus component of $Z$, or it does not, in which case it is the core of an annulus component of $Z$ and again essential. Therefore $\Gamma(X)\union\Gamma(Y)\subseteq\Gamma(Z)$, so that $X\esssubset Z$
and $Y\esssubset Z$. $Z$ is minimal 
with respect to this property 
because if $Z'$ is a competitor then every finite subset of $C$ is
realized geodesically in the geodesic representative of $Z'$, and
hence $Z\esssubset Z'$. Uniqueness follows from the fact that
$\esssubset$ is a partial order. We therefore set $X\essunion Y = Z$. 

In fact we note that $\essunion$ is defined for arbitrary collections
$\{X_i\}$, merely by letting $C = \union \Gamma(X_i)$. Now
we can obtain $X\essint Y$ satisfying (3) by taking the essential union of $\{Z:Z\esssubset X \ \text{and}\ Z\esssubset Y\}$.  
\end{proof}

Here are a few remarks on the proof.

Notice that $X\esssubset Y$ if and only if the geodesic representative of each component of $X$ is pointwise contained in the geodesic representative of some component of $Y$. This is in turn equivalent to saying that each component of $X$ is isotopic to an essential subsurface of a component of $Y$ (where we allow an annulus to be an essential subsurface of itself).

It is helpful to notice that $\Gamma(X\essint Y) = \Gamma(X) \intersect \Gamma(Y)$. This is because any element $\gamma$ in
$\Gamma(X)\intersect\Gamma(Y)$ is the core of an essential annulus $A$ in both, hence $A\esssubset X\essint Y$ by (2), so
$\gamma\in\Gamma(X\essint Y)$. The other direction follows from the fact that $X \essint Y$ is essentially contained in both $X$ and $Y$.

We also define the {\em essential complement}\index{essential!complement}\index{1aaessentialcomplement@$X^c$, essential complement} $X^c$ to be the maximal essential subsurface $Z$ whose geodesic representative is disjoint from that of $X$.  More concretely $X^c$ is the union of complementary components of $X$ that are not 3-holed spheres, together with an annulus for each component of $\boundary X$ that is not isotopic into an annulus of $X$. (This definition agrees with that in Behrstock-Minsky \cite{behrstock-minsky:rank}). Note that essential complement does not behave like a true lattice theoretic complement operator, in that $(X^c)^c$ need not be isotopic to $X$, and $X\essunion X^c$ is usually not $S$; for example, if $X$ is a regular neighborhood of a pants decomposition on $S$ then $X^c = \emptyset$. 

On the other hand, the essential complement does satisfy the following easily verified formula, which in the asymptotic cone will allow us to make sense of codimension:\index{codimension}

\begin{proposition}{CodimensionProp} For each essential subsurface $X \subset S$ we have
$$r(X) + r(X^c) = r(S) = \xi(S)
$$
\qed\end{proposition}

\medskip

\subsubsection{Curve complex.} \index{curve complex}
We associate a simplicial complex $\CC(Y)$\index{1aacurvecomplex@$\CC$, curve complex}
 with each connected surface $Y$ with $\xi(Y)\ge 1$, as well as for each essential subannulus of our ambient surface $S$. For $\xi(Y) \ge  1$, the vertex set $\CC_0(Y)$ of $\CC(Y)$ is $\Gamma(Y)$, the isotopy
classes of essential curves, and for $\xi(Y)>1$, $k$-simplices correspond to sets of $k+1$ vertices with disjoint representatives. Hence $\dim\CC(Y) = \xi(Y)-1$. When $\xi(Y)=1$, we place an edge between any two vertices whose geometric intersection number is the smallest possible on $Y$ (1 when $Y=S_{1,1}$ and 2 when $Y=S_{0,4}$). See Harvey \cite{harvey:boundary}, Ivanov \cite{ivanov:complexes1} and Masur-Minsky
\cite{masur-minsky:complex1}.

If $Y$ is a connected essential subsurface of 
$S$ and $\xi(Y) \ge 1$ then the inclusion $Y \inject S$ naturally
induces an embedding $\CC_0(Y) \inject \CC_0(S)$, whose image we
identify with $\CC_0(Y)$. 
If furthermore $\xi(Y) \ge 2$ then the
embedding $\CC_0(Y) \inject \CC_0(S)$ extends to a simplicial
embedding $\CC(Y) \inject \CC(S)$, whose image we identify with
$\CC(Y)$.

As in \cite{masur-minsky:complex2}, we define $\CC(Y)$ for an essential annulus $Y \subset S$ by considering the annular cover of $S$ to which $Y$ lifts homeomorphically, and which has a natural compactification to a compact annulus $\Ann_Y$ (inherited from the usual compactification of the universal cover $\Hyp^2$). Define an \emph{essential arc} in $\Ann_Y$ to be an embedded arc with endpoints on different components of $\bdy\Ann_Y$. We define $\CC(Y)$ to be the graph whose vertices $\CC_0(Y)$ are isotopy classes rel endpoints of essential arcs in $\Ann_Y$, with an edge for each pair of distinct vertices represented by essential arcs with disjoint interiors.

Note that if $Y,Y'$ are isotopic essential annuli then $\CC(Y),
\CC(Y')$ are \emph{the same} complex.

Given an essential curve $\gamma$ in $S$ we let $\CC(\gamma)$ denote $\CC(Y)$ for any essential annulus $Y \subset S$ with core curve isotopic to $\gamma$.  

The mapping class of the Dehn twist about $\gamma$ acts naturally on
$\CC(\gamma)$ as follows: choose the twist $\tau_\gamma$ to be
supported on an annulus neighborhood $Y$ of $\gamma$, lift
$\tau_\gamma$ through the covering map $\interior(\Ann_Y) \to S$ to a
homeomorphism $\wt\tau_\gamma \from \Ann_Y \to \Ann_Y$ that is
supported on the preimage of $Y$, and let this lift act on the essential arcs in
$\Ann_Y$.  The following properties of this action are easy; for
details see~\MMTwo.

\begin{lemma}{PropAnnulusComplex}
For any essential curve $\gamma$ in $S$ and any vertex $v \in
\CC_0(\gamma)$ the orbit map $n \mapsto \tau^n_\gamma(v)$ is a
quasi-isometry between $\Z$ and $\CC(\gamma)$.  The action of the
infinite cyclic group $\<\tau_\gamma\>$ on $\CC(\gamma)$ has a
fundamental domain of diameter~2.
\end{lemma}

We will need to use  the main result of \MMOne:

\begin{theorem}{hyperbolicity of C}
For each surface $Y$ with $\xi(Y)\ge 1$,  the curve complex $\CC(Y)$ is an infinite diameter $\delta$-hyperbolic metric space, with respect to the simplicial metric.
\end{theorem}

\subsubsection{Markings and partial markings.} We define markings and the marking complex for any connected surface $Y$ with $\xi \ge 1$, as well as for any essential subsurface of $S$, including those which are disconnected and/or have some annulus component. We also define partial markings. (In \cite{masur-minsky:complex2}, partial markings are called markings, and markings are called complete markings).

Suppose that $Y$ is connected and $\xi(Y) \ge 1$.  A \emph{partial marking}\index{marking!partial} $\mu = (\base(\mu),t)$ on $Y$ consists of a simplex $\base(\mu)$\index{1aabase@$\base(\mu)$, base of a partial marking} in $\CC(Y)$ together with a choice of element $t(b) \in \CC_0(b)$,\index{1aatransversal@$t(b)$, transversal of a partial marking}
 which we call a \emph{transversal}, for some (possibly none) of the vertices $b\in\base(\mu)$; by convention we allow the empty set $\emptyset$ as a partial marking of~$Y$.  If $t(b)$ is defined then we say that $b$ is \emph{marked (by $\mu$)},\index{marked} otherwise $b$ is \emph{unmarked (by $\mu$)}.\index{unmarked}  A \emph{marking}\index{marking} (sometimes \emph{full marking})\index{marking!full} is a maximal partial marking, one for which $\base(\mu)$ is a pants decomposition and every $b \in \base(\mu)$ is marked. Given two partial markings $\mu=(\base(\mu),t)$ and $\mu'=(\base(\mu'),t')$, we write $\mu \subset \mu'$ to mean that $\base(\mu)\subset\base(\mu')$ and, for each $b \in \base(\mu)$, $b$ is marked by $\mu$ only if it is marked by $\mu'$ in which case $t(b)=t'(b)$; we also write $\mu \isotopic \mu'$ to mean $\mu \subset \mu'$ and $\mu' \subset \mu$.

Consider now a connected essential subsurface $Y$ of $S$. If $\xi(Y) \ge 1$ then we have defined above markings and partial markings of $Y$. If $Y$ is an annulus then a marking of $Y$ is is just a vertex of $\CC(Y)$, and a \emph{partial marking of $Y$} is either a marking of $Y$ or $\emptyset$.

Finally, given an arbitrary essential subsurface $Y$ of $S$, a marking of $Y$ simply means a choice of marking on each component of $Y$.

\subsubsection{The marking complex}\index{marking complex} We define the marking complex of any connected surface $Y$ with $\xi \ge 1$, and of any essential subsurface of $S$. 

First, given an essential subannulus $Y$ of $S$, define the marking complex of $Y$ to be $\MM(Y) = \CC(Y)$.\index{1aamarkingcomplex@$\MM$, marking complex}

Suppose now that $Y$ is connected and $\xi(Y) \ge 1$. The vertices of $\MM(Y)$ are just the markings of $Y$. To define the edges we first need this notion: if $b$ and $c$ are overlapping essential curves, we denote $\pi_b(c)\in\CC(b)$ to be the set of lifts of $c$ to essential arcs in the annular cover associated to $b$. The diameter of this set is bounded in $\CC(b)$, with a bound depending only on the topology of $Y$. The map $\pi_b$ is an example of a subsurface projection map, defined below in a more general setting. 

Edges in $\MM(Y)$ correspond to  \emph{elementary moves} among markings on $Y$, which come in two flavors: \emph{twist moves} and \emph{flip moves}.  To define them, consider a marking $\mu$ on $Y$ and a curve $b \in \base(\mu)$. A marking $\mu'$ is said to be obtained from $\mu$ by a \emph{twist
move} about $b$ if $\base(\mu)=\base(\mu')$, $\mu,\mu'$ have the same
transversals to each curve in 
$\base(\mu)\setminus\{b\}=\base(\mu')\setminus\{b\}$, and the
transversals $t(b)$ in $\mu$ and $t'(b)$ in $\mu'$ satisfy
$d_{\CC(b)}(t,t') \le 2$.  A marking $\mu' = (\base(\mu'),t')$ is said to be obtained from $\mu$
by a \emph{flip move} along $b$ if there exists $b' \in \base(\mu')$
such that $\base(\mu)\setminus\{b\} = \base(\mu')\setminus\{b'\}$, 
$\Fill(b,b')$ is a one-holed torus or 4-holed sphere $W$ such that
$d_{\CC(W)}(b,b') = 1$, 
$d_{\CC(b)}(\pi_b(b'),t(b))\le 2$, and 
$d_{\CC(b')}(t'(b'),\pi_{b'}(b))\le 2$.

Finally, for any essential subsurface $Y$ of $S$ with components $Y_1,\ldots,Y_n$, define $\MM(Y)$ to be the 1-skeleton of the cartesian product $\MM(Y_1) \cross \ldots \cross \MM(Y_n)$ with the usual CW-product structure. To put it another way, the vertices of $\MM(Y)$ are the markings of $Y$, and there is an edge between two markings $\mu,\mu'$ of $Y$ if and only if $\mu,\mu'$ are isotopic outside of a certain component $Y_i$, and the restrictions of $\mu,\mu'$ to $Y_i$ are connected by an edge of $\MM(Y_i)$.

The marking complex is locally infinite because of the structure of
transversals, but it is still quasi-isometric to a locally finite complex, as follows. Recall from \cite{masur-minsky:complex2} that a \emph{clean marking} is a marking $\mu=(\base(\mu),t)$ with the following properties: for each $b \in \base(\mu)$, $t(b)$ is $\pi_b(c)$ where $c=c(b)$ is an essential curve in the component $F$ of $Y\setminus(\base(\mu)\setminus\{b\})$ containing $b$; and the curves $b$ and
$c$ have minimal nonzero intersection number in $F$.  The complex of clean markings is a connected, $\MCG(Y)$-invariant subcomplex of $\MM(Y)$ whose vertices are the clean markings.

In fact the clean marking complex is
what is usually referred to as the marking complex, 
see e.g., \cite{behrstock:asymptotic} 
and \cite{behrstock-minsky:rank}. Because the full complex is more
convenient for our purposes, we record this quasi-isometry: 

\begin{proposition}{M qi MCG}
The marking complex $\MM(Y)$ is quasi-isometric to $\MCG(Y)$ and to
the subcomplex of clean markings.  More
precisely, for each $\mu_0 \in \MM(Y)$ the orbit map $\phi \mapsto
\phi(\mu_0)$ is a quasi-isometry from $\MCG(Y)$ to $\MM(Y)$.
\end{proposition}

In particular $\MM(Y)$ is connected, which may not have been obvious from the definition.

\begin{proof}[Sketch of proof]
As noted in \MMTwo, the clean marking complex is locally finite, the
action of $\MCG(Y)$ is properly discontinuous and cocompact, and so by
the Milnor-\v{S}varc lemma \cite{milnor:growth} \cite{svarc:growth} the
orbit map is a quasi-isometry between $\MCG(Y)$ and the complex of
clean markings.  The inclusion of the complex of clean markings into
$\MM(Y)$ is an $\MCG(Y)$-equivariant quasi-isometry, because for each
marking there is a clean marking within a uniformly bounded distance
$C$ by Lemma \ref{PropAnnulusComplex}, and for each edge of $\MM(Y)$
connecting two markings $\mu_0,\mu_1$, if $\mu'_0,\mu'_1$ are two
clean markings within distance $C$ of $\mu_0,\mu_1$ respectively then
the distance between $\mu'_0,\mu'_1$ in the clean marking complex is
uniformly bounded --- this is checkable explicitly from the definition
in \MMTwo\ of the edges allowed between clean markings.
\end{proof}

\paragraph{Remark on notation.} In any context where $\CC(Y)$ is under consideration the essential subsurface $Y$ is assumed to be connected, whereas when $\MM(Y)$ is being considered then $Y$ can be disconnected.

\subsubsection{Overlap.} We define a symmetric binary relation of \emph{overlap}\index{overlap} for objects on~$S$, denoted $\pitchfork$,\index{1aaoverlap@$\pitchfork$, overlap} as follows.

We have already defined overlap of two essential curve systems, in Section~\ref{SubsubsectionCurves}.

Overlap of an essential curve $\gamma \subset S$ and an essential subsurface $Y \subset S$,\index{overlap!of surface and curve} denoted $\gamma \pitchfork Y$ and $Y\pitchfork \gamma$, means that $\gamma$ cannot be isotoped into the complement of~$Y$.  Equivalently, after isotoping $\gamma$ to intersect $\bdy Y$ efficiently, the intersection $\gamma \intersect Y$ is either a non-boundary-parallel essential curve in $Y$ (the core of an annulus component is not allowed) or a nonempty pairwise disjoint union of \emph{essential arcs in $\overline Y$}, each a properly embedded arc $\alpha \subset \overline Y$ that is not homotopic rel endpoints into $\bdy Y$.

Given an essential curve system $\cC$ and an essential
subsurface $Y$, define $\cC \pitchfork Y$ to mean that there exists a
component $\gamma$ of $\cC$ such that $\gamma \pitchfork Y$.

Overlap of two essential subsurfaces $X,Y \subset S$,\index{overlap!of two surfaces} denoted $X \pitchfork Y$, means that $\bdy Y \pitchfork X$ and $\bdy X \pitchfork Y$. Equivalently, after
$X,Y$ are isotoped so that $\bdy X$, $\bdy Y$ intersect efficiently,
some component of $\overline Y \intersect \bdy X$ is an essential
curve or arc in $Y$ and some component of $\bdy Y \intersect \overline X$ is an essential curve or arc in $X$.  It is also equivalent to say that neither of $Y$ or $X$ is isotopic into the other and no
matter how $Y,X$ are chosen in their isotopy classes their
intersection is nonempty. Notice for example that if $\bdy X \pitchfork \bdy Y$ then $X \pitchfork Y$, but the converse can fail.

Overlap of an essential subsurface $Y$ and a partial marking $\mu$, \index{overlap!of surface and partial marking} denoted $\mu\pitchfork Y$ and $Y\pitchfork \mu$, means that either $\base(\mu) \pitchfork Y$ or there exists $b \in \base(\mu)$ such that $b$ is marked by $\mu$ and some component of $Y$ is an annulus neighborhood $b$.  

\subsubsection{Open subsurface and support.} \qquad Given a partial marking $\mu=(\base(\mu),t)$ on $S$, its \emph{open subsurface},\index{open subsurface of a partial marking}\index{1aaopensubsurface@$\open(\mu)$, open subsurface of a partial marking} denoted $\open(\mu)=\open_S(\mu)$, is defined to be the essential union of all subsurfaces $Z$ such that $Z\not\pitchfork\mu$. Equivalently, $\open(\mu)$ is the largest essential subsurface which does not overlap $\mu$. One can also describe it explicitly as the union of the components $F$ of $S-\base(\mu)$ such that $\xi(F) \ge 1$, and the annuli (if any) homotopic to the {\em unmarked} $b \in \base(\mu)$. Note that each boundary curve of $\open(\mu)$ is isotopic to a curve of $\base(\mu)$. We usually drop the subscript $S$ in the notation $\open_S(\mu)$, unless we want to emphasize the surface in which the operation takes place, as in the proof of Lemma \ref{juncture}.

The \emph{support}\index{support of a partial marking}\index{1aasupport@$\supp(\mu)$, support of a partial marking} of a partial marking $\mu$ of $S$, denoted $\supp(\mu)=\supp_S(\mu)$, is defined to be $\open(\mu)^c$, the essential complement of $\open(\mu)$. The support of $\mu$ does not always behave as might at first be expected: for example, if no transversals are defined in $\mu$ then $\supp(\mu) = \emptyset$, even if $\mu \ne \emptyset$. 

We note several properties of $\open(\mu)$ and $\supp(\mu)$. Let $\mu = \mu^u \union \mu^m$ where $\mu^u$ consists of the unmarked curves of $\base(\mu)$ and $\mu^m$ consists of the marked curves and their transversals.
\begin{itemize}
\item Each component of $\bdy\supp(\mu)$ is isotopic to a component of $\base(\mu)$.
\item $\mu^m$ restricts to a (full) marking of each component of $\supp(\mu)$.
\item $\supp(\mu)$ is characterized up to isotopy as the maximal essential subsurface of $S$ with respect to the previous two properties. 
\item The curves $\mu^u$ are precisely the cores of the annulus components of $\open(\mu)$.
\end{itemize}

\subsubsection{Dehn twist flats.} \label{SectionDehnTwistFlats}
Given a pants decomposition $\mu$ on $S$, the set of markings whose base is $\mu$ is denoted $\QQ(\mu) \subset \MM(S)$, and is called the \emph{Dehn twist flat}\index{Dehn twist!flat} corresponding to $\mu$. The terminology ``flat'' comes from the fact, proved in \BMRank\ Lemma~2.1, that $\QQ(\mu)$ is a quasi-isometrically embedded copy of $\reals^{\xi(S)}$; see also Proposition~\ref{Q product structure}~(1) below, which generalizes this  to allow $\mu$ to be any partial marking.

We state here without proof two other ways to view Dehn twist flats, in two other models of the quasi-isometric geometry of $\MCG(S)$: the word metric on $\MCG(S)$; and the thick part of \Teichmuller\ space.

Under the quasi-isometry $\MM(S) \leftrightarrow \MCG(S)$ of Proposition~\ref{M qi MCG}, a Dehn twist flat $\QQ(\mu)$ corresponds uniformly to a left coset of a maximal-rank free-abelian subgroup generated by Dehn twists. To be precise, consider the finite set of $\MCG(S)$ orbits of pants decompositions, and choose one representative $\mu_1,\ldots,\mu_K$ from each orbit.  There is a constant $A \ge 0$ such that for each pants decomposition $\mu = \Phi \cdot \mu_k$, the image of $\QQ(\mu)$ in $\MCG(S)$ has Hausdorff distance $\le A$ from the left coset $\Phi \, T(\mu_k)$ where $T(\mu_k)$ is the rank $\xi(S)$ free abelian group generated by Dehn twists about the components of $\mu_k$.

Another quasi-isometric model of $\MCG(S)$ is the thick part of the
\Teichmuller\ space $\Teich(S)$. By Margulis' Lemma there is a constant
$\epsilon_0>0$ independent of $S$ such that in any hyperbolic structure on
$S$, the set of simple closed geodesics of length $\le \epsilon_0$ is pairwise
disjoint,  and every other simple closed geodesic has length $\ge 2
\epsilon_0$. The thick part $\Teich_{\text{thick}}(S)$ is defined to be the
set of hyperbolic structures whose shortest closed curve has length $\ge
\epsilon_0$. The action of $\MCG(S)$ on $\Teich_{\text{thick}}(S)$ is properly
discontinuous and cocompact and so, by the \v{S}varc-Milnor lemma, there is a
quasi-isometry $\MM(S) \leftrightarrow \Teich_{\text{thick}}(S)$, with respect
to any equivariant proper geodesic metric on
$\Teich_{\text{thick}}(S)$. Applying Fenchel--Nielsen coordinates one sees
that $\Teich_{\text{thick}}(S)$ is a \emph{manifold-with-corners}, locally
modelled on the closed orthant $\{x \in \reals^{2 \xi(S)} \suchthat \,  \,
\text{$x_i \ge 0$ for all $i$}\}$. Each curve family $C$ of cardinality $j$
corresponds to a codimension-$j$ facet $F(C)$ consisting of hyperbolic
structures on which the curves of $C$ are precisely the curves of length
$\epsilon_0$. Note that for a pants decomposition $\mu$, the Fenchel-Nielsen
length coordinates in the facet $F(\mu)$ are all fixed to be $\epsilon_0$, and
so the remaining coordinates are just the twists around the curves of $\mu$,
making $F(\mu)$ homeomorphic to $\reals^{\xi(S)}$. Putting this all together,
there is a constant $A \ge 0$ such that for each pants decomposition $\mu$,
the image of $\QQ(\mu)$ in $\Teich_{\text{thick}}(S)$ has Hausdorff distance
$\le A$ from $F(\mu)$.

\subsection{Projections}
\label{SectionProjections}

\subsubsection{Subsurface projections.} \index{subsurface projections} Following \MMTwo, \cite{behrstock:asymptotic}, and \BMRank, given a surface $S$  and an essential subsurface $Y$ we shall define several projection maps which take curves and markings in $S$ to curves and markings in $Y$. When the target is the marking complex $\MM(Y)$ we use $\pi_{\MM(Y)}$\index{1aaprojectionmarking@$\pi_{\MM(Y)}$, projection to marking complex} for the projection, and when the target is the curve complex $\CC(Y)$ we use $\pi_{\CC(Y)}$\index{1aaprojectioncurve@$\pi_{\CC(Y)}$, projection to curve complex} or more briefly $\pi_Y$.\index{1aaprojectioncurvebrief@$\pi_Y$, projection to curve complex}  Because these definitions require choices, for example choosing a vertex among a set of vertices, formally speaking we define the image of each map to be the set of all choices. However in all cases the map is \emph{coarsely well-defined}\index{coarsely well-defined} (see Lemma~\ref{coarse definition}), which means that the set of choices has uniformly bounded diameter. In practice we may sometimes abuse terminology and treat these maps as if the image of a point is a point. We may also abuse distance notation between the images of two points, confusing minimum distance, Hausdorff distance, and distance between any two representative elements of the images (see below under ``Notation for various distances''), because those quantities all differ by a uniformly bounded amount. Furthermore, these subsurface projection maps are coarsely Lipschitz (Lemma~\ref{coarse lipschitz}) and they behave well with respect to composition (Lemma~\ref{coarse composition}).

\paragraph{Projecting curves to (sets of) curves.} 
Suppose that $Y\esssubset S$ is connected and not an annulus. 
If $\gamma\in\CC(S)$, we define $\pi_{\CC(Y)}(\gamma)$ to be the set
of vertices of $\CC(Y)$ obtained from essential arcs or curves of
intersection of $\gamma$ with $Y$, by the process of surgery along
$\boundary Y$. To be more precise, put $\gamma$ in efficient position
with respect to $\bdy Y$, choose a component $\alpha$ of $\gamma
\intersect Y$, and consider a component of the boundary of a
regular neighborhood of $\alpha \union \bdy Y$; the set of all essential curves in $Y$ obtained in this way is $\pi_{\CC(Y)}(\gamma)$. 

If $Y$ is an annulus, we let $\pi_{\CC(Y)}(\gamma)$ be the set of vertices of $\CC(Y)$ obtained as lifts of $\gamma$ to the annular cover associated to $Y$. This operation was denoted $\pi_b$ in the earlier discussion of marking complexes, where $b$ is the core of~$Y$.

Note in both cases that $\pi_{\CC(Y)}(\gamma)$ is nonempty if and only if $\gamma\pitchfork Y$. See Lemma~\ref{coarse definition} for coarse well-definedness of $\pi_{\CC(Y)}$. 

Notation: we often write $\pi_Y$ for any projection map whose target
is $\CC(Y)$. When the target needs to be emphasized we
write $\pi_{\CC(Y)}$.

The {\em bounded geodesic projection theorem} from \MMTwo\ will be important for us:

\begin{theorem}{bounded geodesic projection}
Let $Y$ be a connected essential subsurface of $S$ and let $g$ be a geodesic segment in $\CC(S)$ such that $v\pitchfork Y$ for every vertex $v$ of $g$. Then 
$$\diam_{\CC(Y)}(g) \le B
$$
where $B$ is a constant depending only on $\xi(S)$.
\end{theorem}

\paragraph{Projecting (partial) markings to curves.} 
We define a projection $\pi_{\CC(Y)}(\mu) \subset \CC(Y)$ for a partial marking $\mu$
in $S$  as follows:
When $\base(\mu) \pitchfork Y$ we let
$\pi_{\CC(Y)}(\mu)$ be the union of $\pi_{\CC(Y)}(b)$ over all $b \in
\base(\mu)$. When $Y$ is an annulus neighborhood of a marked $b
\in \base(\mu)$ then we define $\pi_{\CC(Y)}(\mu) = t(b)$. In all other cases
$\pi_{\CC(Y)}(\mu) = \emptyset$. 

\paragraph{Projecting (partial) markings to (partial) markings.}
If $\mu$ is a partial marking in $S$  and $Y$ an essential subsurface,
we will define a partial marking $\pi_{\MM(Y)}(\mu)$ in $Y$. If $\mu$ is a marking of $S$ then $\pi_{\MM(Y)}(\mu)$ will be a marking of $Y$, so that we will obtain a coarse Lipschitz map $\pi_{\MM(Y)} \from \MM(S) \to \MM(Y)$ (see Lemma~\ref{coarse lipschitz}). 

Write $\mu =(\base(\mu),t)$. If $Y\not\pitchfork\mu$ then
$\pi_{\MM(Y)}(\mu) =\emptyset$. From now on we may assume
$Y\pitchfork\mu$. If $Y$ is disconnected we can view partial markings
as tuples of partial markings in the components and define 
$\pi_{\MM(Y)}$ componentwise. So we may also assume that $Y$ is connected.

When  $Y$ is an annulus let $\pi_{\MM(Y)}(\mu)$ denote any choice of
element of $\pi_{\CC(Y)}(\mu)$, recalling that $\MM(Y) = \CC(Y)$. 

If $\xi(Y)\ge 1$, let $b$ be any choice of element in
$\pi_{\CC(Y)}(\mu)$, let $A$ be an annulus with core $b$, and let
$Y_{b}$ be the union of $A$ with its essential complement in $Y$. 
Now inductively define
$$
\pi_{\MM(Y)}(\mu) = b \union \pi_{\MM(Y_{b})}(\mu)
$$
where the second term of the union is interpreted as a union over the
components of $Y_{b}$. Note that, at the bottom of the induction, the
annulus case provides transversals for all the base curves of $\pi_{\MM(Y)}(\mu)$ that overlap $\mu$.

There are choices at each stage of this construction, but when $\mu$ is a marking the final output is coarsely well-defined, as proved in
\cite{behrstock:asymptotic}. See Lemma \ref{coarse definition} for a statement.

For a partial marking $\mu$, recall that
$\open(\mu)$ is the unique maximal essential subsurface $Y$ such that
$\mu\not\pitchfork Y$, or equivalently such that
$\pi_{\MM(Y)}(\mu)$ is empty.
The following lemma characterizes the (relative) open subsurface of the projection of a
partial marking, as the maximal subsurface that doesn't overlap the marking:
\begin{lemma}{characterize open projection}
If $Y$ is an essential subsurface of $S$ and $\mu$ a partial marking
in $S$, then 
$$
\open_Y(\pi_{\MM(Y)}(\mu)) = \essUnion\left\{Z\esssubset Y : Z\not\pitchfork\mu\right\}.
$$
\end{lemma}

\begin{proof}
Let $\mu' = \pi_{\MM(Y)}(\mu)$. Every base curve in the inductive construction of $\mu'$ is either a base curve of $\mu$ itself, or an element of a subsurface projection of $\mu$ into some subsurface of $Y$. The induction terminates when the complementary subsurfaces 
of the base have no more overlap with $\mu$, and when the base curves are either marked by $\mu'$ or are base curves of $\mu$ that have no transversals. It follows that $\open_Y(\mu')$ does not overlap $\mu$. 

Conversely, let $Z\esssubset Y$ be a subsurface that does not overlap
$\mu$. If $Z$ is an annulus around $b\in\base(\mu)$ then $b$ is
unmarked by $\mu$. Otherwise $Z$ is disjoint from all vertices of
$\pi_{\CC(Y)}(\mu)$. Hence in the first step of the construction of
$\mu'$, the subsurface $Z$ does not overlap the chosen base curve. Continuing by induction, $Z$ does not overlap $\mu'$. Hence $Z\esssubset \open_Y(\mu')$. 

We have shown that $\open_Y(\mu')$ is among the
set of subsurfaces of $Y$ that don't overlap $\mu$, and that
every subsurface of $Y$ that doesn't overlap $\mu$ is essentially contained
in $\open_Y(\mu')$. Hence the two sides are equal. 
\end{proof}

\paragraph{Notation for various distances:}\index{distance notation} In a metric space $M$, when $p,q$ are subsets of $M$ we use $d_M(p,q)$ to denote \emph{minimum distance},\index{minimum distance}\index{1aaminimumdistance@$d_M$, minimum distance in $M$} meaning the infimum of the distance between an element of $p$ and an element of $q$. We also use $d_{H,M}(p,q)$ to denote \emph{Hausdorff distance},\index{Hausdorff distance}\index{1aahausdorffdistance@$d_{H,M}$, Hausdorff distance in $M$}\index{1aahausdorffdistance@$d_H$, Hausdorff distance} the infimum of all $\epsilon \ge 0$ such that $p$ is contained in the $\epsilon$-neighborhood of $q$ and $q$ is contained in the $\epsilon$-neighborhood of $p$. When the context is clear we abbreviate $d_{H,M}(p,q)$ to $d_H(p,q)$.

Given an essential subsurface $Y \subset S$, and any objects $a$ and $b$ in the domain of $\pi_{\CC(Y)}$ or of $\pi_{\MM(Y)}$ we denote 
$$d_{\CC(Y)}(a,b) = d_{\CC(Y)}(\pi_{\CC(Y)}(a),\pi_{\CC(Y)}(b))
$$
and 
$$d_{\MM(Y)}(a,b) = d_{\MM(Y)}(\pi_{\MM(Y)}(a),\pi_{\MM(Y)}(b))
$$
as long as the right hand side makes sense, for example when $a,b$ are full markings in the latter case. When the context is clear we often abbreviate $d_{\CC(Y)}$ to $d_Y$.

\subsubsection{Quasidistance formula.} \index{quasidistance formula} Knowing that $\MCG(S)$ is quasi-isometric to the marking complex $\MM(S)$, we can study the asymptotic geometry of $\MCG(S)$ by having a useful quasidistance formula on $\MM(S)$, which is provided by the following.

Given two numbers $d \ge 0$, $A \ge 0$, denote the truncated distance\index{1aatruncateddistance@$\truncate{d}{A}$, truncated distance} 
$$\truncate{d}{A} = \begin{cases}
d & \text{if}\,\, d \ge A \\
0 & \text{otherwise}
\end{cases}
$$
Given $r,s \ge 0$, $K \ge 1$, $C \ge 0$ we write $r \ceq{K,C} s$ to
mean that $\frac{1}{K} s - C \le r \le Ks+C$.  We also write $r
\approx s$ to mean that $r \ceq{K,C} s$ for some \emph{constants of
approximation} $K,C$ which are usually specified by the context, and
we similarly write $r \lesssim s$ to mean $r \le Ks+C$.

The following result is proved in \MMTwo.

\begin{theorem}{distance formula}
There exists a constant $A_0 \ge 0$ depending only on the topology of
$S$ such that for each $A \ge A_0$, and for any $\mu,\mu' \in \MM(S)$ 
we have the estimate 
$$d_{\MM(S)}(\mu,\mu')
\approx \sum_{Y\subseteq S} \truncate{d_{\CC(Y)}(\mu,\mu')}{A}
$$
and the constants of approximation depend on $A$ and on the topology
of~$S$. 
\end{theorem}

The constant $A$ in this theorem is usually called the \emph{threshold
constant}.

\paragraph{Remark:}
In summations and other expressions with index $Y$ as in the above
theorem, the convention will be that the index set consists of one
representative $Y$ in each isotopy class of {\em connected} essential
subsurfaces, perhaps with some further restriction on the isotopy
class; see for example Proposition~\ref{Q product structure}(2) below.

There are two ways that Theorem \ref{distance formula} is applied.
First, we can raise 
the threshold with impunity, which can make some terms drop out in a
way that the remaining terms are more easily described.  Second, if
each term in the sum is replaced by another term differing by at most
a uniform constant $C \ge 0$ then, after raising the threshold above
$2C$, we may make the replacement at the cost of a multiplicative
factor of at most~2. As an example we have the following:

\begin{corollary}{C bounds M}
For any $r$ there exists $t$ such that for any $\mu,\nu\in\MM(S)$, if
$d_{\CC(W)}(\mu,\nu) \le r$ for all $W\subseteq S$, then 
$d_{\MM(W)}(\mu,\nu) \le t$. 
\end{corollary}

The proof is simply to raise the threshold above $r$, so that the
right hand side of the quasidistance formula becomes 0.  In fact a
more careful look at the machinery of \cite{masur-minsky:complex2}
yields $t = O(r^{\xi(S)})$.

\subsubsection{Basic properties.} We conclude this section with a brief summary of some of the basic properties of projections: 

\begin{lemma}{coarse definition} 
Subsurface projections are {\em coarsely well defined}:\index{coarsely well-defined}
\begin{itemize}
\item The diameter of $\pi_{\CC(Y)}(x)$, where $x$ is a curve or marking in $S$, is uniformly bounded. 
\item
Similarly, the diameter of all possible choices in the construction of
$\pi_{\MM(Y)}(\mu)$, for $\mu\in\MM(S)$, is uniformly bounded. 
\end{itemize}
\end{lemma}

\smallskip

\begin{lemma}{coarse lipschitz}
Subsurface projections are {\em coarsely Lipschitz} in the following
sense: 
\begin{itemize}
\item If $x,y\in\CC(S)$ with $d(x,y)=1$ and both $x\pitchfork Y$ and
$y\pitchfork Y$, then $\diam(\pi_{\CC(Y)}(x)\union\pi_{\CC(Y)}(y))$ is
uniformly bounded. 
\item 
Similarly $d_{\MM(Y)}(\mu,\nu)$ is uniformly bounded for any $\mu,\nu$ in $\MM(S)$ with $d(\mu,\nu) = 1$. 
\end{itemize}
\end{lemma}

In Lemmas~\ref{coarse definition} and~\ref{coarse lipschitz}, ``uniformly bounded'' means bounded by a constant depending only on $\xi(S)$. In fact, for the $\CC(Y)$ bounds of these lemmas, there is a bound of $3$ when $Y$ is an annulus and of $2$ when $\xi(Y) > 1$ --- see \cite{masur-minsky:complex2} for details.

\begin{lemma}{coarse composition}
Subsurface projections for nested subsurfaces are {\em coarsely
composable} in the following sense. Let  $X,Y \subset S$ be essential subsurfaces such that $X\esssubset Y$. 
\begin{itemize}
\item For any $\gamma \in \CC_0(S)$, the curve $\gamma$ overlaps  $X$ if and only if $\gamma$ overlaps $Y$ and $\pi_Y(\gamma)$ contains at least one element $\alpha$ which overlaps
$X$. In this case, $\pi_X$ is coarsely equivalent to $\pi_X\circ\pi_Y$, meaning that
$$
\diam_{\CC(X)}(\pi_X(\gamma)\union\pi_X(\pi_Y(\gamma)))
$$
is uniformly bounded.
\item 
Similarly, if $\mu\in\MM(S)$ then
$$
d_{\MM(X)}(\pi_{\MM(X)}(\mu), \pi_{\MM(X)}(\pi_{\MM(Y)}(\mu)))
$$
is uniformly bounded.
\end{itemize}
\end{lemma}

We remark that for all three of these lemmas, the statements for curve
complex projections are elementary from the definitions, and the
statements for marking projections follow easily from the quasidistance
formula.

\subsection{Asymptotic cones}
\label{cone defs}

The asymptotic cone of a metric space is a way to encode the geometry 
of that space as seen from arbitrarily large distances. 
We will discuss this construction and the notation we will be using (see \cite{dries-wilkie, gromov:polynomial} for
further details). 
 
To start, we recall that a \emph{(non-principal) ultrafilter}\index{ultrafilter} is a finitely additive probability measure $\omega$ defined on the power set of the natural numbers, which takes values only $0$ or $1$, and for which every finite set has zero measure. If two sequences coincide on a set of indices whose $\omega$-measure is equal to~1 then they they are said to be \emph{$\omega$-equivalent}, or to coincide $\omega$-a.e.\index{1ssz@$\omega$-a.e.} or $\omega$-a.s.\index{1ssz@$\omega$-a.s.} We will have a tendency, especially later in the paper, to abbreviate the terminology by speaking of $\omega$-equivalent objects as being ``equal'' or ``the same''.

The \emph{ultraproduct}\index{ultraproduct} of a sequence of sets $X_n$ is the quotient $\Pi_n X_n/\sim$ of the cartesian product identifying two sequences $(x_n)$, $(y_n)$ if they are $\omega$-equivalent. We will often use the notation $\bbar X$ for the ultraproduct, and we use $\bbar x$ or $\uprod{x_n}$\index{1aaultraequivalence@$\bbar x$, ultraproduct equivalence class}\index{1aaultraequivalence@$\uprod{x_n}$, ultraproduct equivalence class} for the $\omega$-equivalence class of a sequence $(x_n)$, also called its \emph{ultraproduct equivalence class}.\index{ultraproduct equivalence class}

In  a topological space $X$, the  \emph{ultralimit}\index{ultralimit!of a sequence of points}\index{1aaultralimit@$\ulim$, ultralimit} of a sequence of points $(x_{n})$ is $x$, denoted $x=\ulim x_{n}$, if for every neighborhood $U$ of $x$ the set  $\{n:x_n \in U\}$ has $\omega$-measure equal to 1. Ultralimits are unique when they exist, moreover two $\omega$-equivalent sequences in $X$ have the same ultralimit. In this language the Bolzano-Weierstrass Theorem says that when $X$ is compact every sequence has an ultralimit. 

The \emph{ultralimit}\index{ultralimit!of a sequence of metric spaces} of a sequence of based metric spaces $(X_{n}, x_{n}, \dist_{n})$ is defined as follows: for $\bbar y,\bbar z\in \bbar X$, we define
$\dist(\bbar y,\bbar z) = \ulim\dist_n(y_n,z_n)$, 
where the ultralimit is taken in the compact set $[0,\infty]$. 
We then let
$$
\ulim (X_{n},x_{n}, \dist_{n})\equiv\{\bbar y: \dist(\bbar y,\bbar x) <
\infty\}/\sim
$$
where we define $\bbar y\sim \bbar y'$ if 
$\dist(\bbar y,\bbar y')=0$. Clearly $\dist$ makes this quotient into a
metric space called the \emph{ultralimit} of the $X_{n}$.

Given a sequence of positive constants  $s_{n}\to\infty$ and a sequence $(x_n)$ of basepoints in a fixed metric space $(X,\dist)$, we may consider the rescaled space $(X,x_n,\dist/s_n)$. The ultralimit of this sequence is called the \emph{asymptotic cone of $(X,\dist)$ relative to the ultrafilter $\omega$, scaling constants $s_{n}$, and basepoint $\bbar x = \uprod{x_{n}}$}:\index{asymptotic cone}\index{1aaasymptoticcone@$\cone$, asymptotic cone}
$$\cone(X,(x_{n}),(s_{n}))=\ulim (X, x_{n}, \frac{\dist}{s_{n}}).
$$
For the image of $\bbar y$ in the asymptotic cone, the \emph{rescaled ultralimit}\index{rescaled ultralimit} we use the notation either $y_{\omega}$ or $\seq y$.\index{1aarescaledultralimit@$y_\omega$, rescaled ultralimit of points}\index{1aarescaledultralimit@$\seq y$, rescaled ultralimit of points} Given a sequence of subsets $A_n \subset X$ we use the notation $A_\omega$ \index{1aarescaledultralimit@$A_\omega$, rescaled ultralimit of sets} to denote the subset of  $\cone(X,(x_n),s_n))$ consisting of all $y_\omega$ for sequences $(y_n)$ such that $y_n \in A_n$ for all $n$.

The rescaling limit works equally well for a sequence $(X_n,x_n,\dist_n)$ in which the $X_n$ are not all the same metric space, and so we call $\lim_\omega(X_n,x_n,\dist_n/s_n)$ the asymptotic cone of the sequence. 

\paragraph{Convention:} For the remainder of the paper we fix a non-principal ultrafilter $\omega$. Usually also the scaling sequence $s_n\to \infty$ and the basepoint $\mu_0$ for $\MM(S)$ are implicitly fixed, and we write $\AM=\AM(S)$\index{1aaasymptoticmarking@$\AM$, asymptotic cone of marking complex} to denote an asymptotic cone of  $\MM(S)$ with respect to these choices. Any choice of scaling sequence and basepoint will do, but in the last section we will need the flexibility of varying the choice of the scaling sequence~$(s_n)$. Further, properties of linear or sub-linear growth of a non-negative function $f(n)$ are always taken with respect to the choice of $\omega$ and $(s_{n})$, that is, we say $f(n)$ has \emph{linear growth}\index{growth!linear} if $0<\ulim \frac{f(n)}{s_{n}}<\infty$ and \emph{sublinear growth}\index{growth!sublinear} if $\ulim \frac{f(n)}{s_{n}}=0$.

Note that since $\MM(S)$ is quasi-isometric to the group $\MCG(S)$ with any finitely generated word metric, and since the isometry group of $\MCG(S)$ acts transitively, the asymptotic cone is independent of the choice of basepoint.

Any essential connected subsurface $W$ inherits a basepoint $\pi_{\MM(W)}(\mu_0)$, canonical up to bounded error by Lemma~\ref{coarse lipschitz}, and we can use this to define the asymptotic cone $\AM(W)$ of its marking complex $\MM(W)$. For a disconnected subsurface $W=\sqcup_{i=1}^{k}W_{i}$ we have $\MM(W)=\Pi_{i=1}^k\MM(W_i)$ and we may similarly construct $\AM(W)$ which can be identified with $\Pi_{i=1}^{k}\AM(W_{i})$ --- this follows from the general fact that the process of taking asymptotic cones commutes with finite products. Note that for an annulus $A$ we have defined $\MM(A) =\CC(A)$ which is quasi-isometric to $\Z$, so $\AM(A)$ is bilipschitz equivalent (isometric) to~$\R$.

For a sequence $(W_n)$ of subsurfaces we can similarly form the ultraproduct of $(\MM(W_n))$, which we denote by $\MM(\bbar W)$, where $\bbar W = \uprod{W_n}$. The asymptotic cone of this sequence with the inherited basepoints is denoted $\AM(\bbar W)$.  We also let $\bbar S$ denote the constant sequence $(S,S,\ldots)$ so that $\MM(\bbar S)$ is the ultraproduct of $(\MM(S),\ldots)$ and $\AM(\bbar S)$ is the same as $\AM(S)$. 

Any sequence in a finite set $A$ is $\omega$-a.e.\ constant: given $(a_n\in A)$ there is a unique $a\in A$ such that $\omega(\{n:a_n = a\}) = 1$. It follows that the ultralimit $\bbar a$ of this sequence is naturally identified with $a$. For example if $(W_n)$ is a sequence of essential subsurfaces of $S$ then the topological type of $W_n$ is $\omega$-a.e.\ constant, so we call this the topological type of $\bbar W$. Similarly the topological type of the pair $(S,W_n)$ is $\omega$-a.e.\ constant. We can  moreover interpret expressions like
$\bbar U \subset \bbar W$ to mean $U_n\subset W_n$ for $\omega$-a.e.\ $n$, and so on. Note that $\AM(\bbar W)$ can be identified with $\AM(W)$, where $W$ is a surface homeomorphic to $W_{n}$ for $\omega$-a.e.~$n$; this identification is not natural, however, because it depends up to $\omega$-equivalence on choosing a homeomorphism between $W$ and each $W_n$.

For two sequences of sets $(A_n)$ and $(B_n)$ and a sequence of functions $f_n \co A_n \to B_n$,
passing to the ultralimit gives rise to a single function $\bbar f\co\bbar A \to \bbar B$, and $\bbar f$ determines $f_n$ up to $\omega$-equivalence in the ultraproduct of the sequence $(B_n^{A_n})$. 
We can therefore think of a sequence of projection maps
$\pi_{\MM(W_n)}\co\MM(S) \to \MM(W_n)$ as a single map
$$\pi_{\MM(\bbar W)}\co\MM(\bbar S) \to \MM(\bbar W)
$$
Upon rescaling, this map descends to a map of the asymptotic cones, the \emph{rescaled ultralimit of the projection maps} \index{1aarescaledprojection@$\pi_{\AM(\bbar W)}$, rescaled ultralimit of projection maps}
$$
\pi_{\AM(\bbar W)}\co\AM(S)\to \AM(\bbar W).
$$
Note by Lemma \ref{coarse lipschitz} that this is a Lipschitz map.
This sort of notation will be used heavily in Section~\ref{finiteness2}.

\section{Product regions and cubes}
\label{SectionCubes}

In this section we define and study subsets of the marking complex
obtained by holding fixed one part of the surface and varying the
rest. These will be called \emph{product regions}, because of the
product structure described in Lemma~\ref{Q product structure}. A
special case of a product region is a \emph{Dehn twist flat}. 
We will also consider particular subsets of product regions called
\emph{cubes}, which are in fact naturally quasi-isometrically
parametrized by cubes in Euclidean space.

The metric relation between a pair of product regions or cubes will be described in terms of \emph{junctures}, the part of each set which comes closest to the other.  Later, when we pass to the asymptotic cone, these
junctures will become intersections, and will be important in
understanding the structure of orthants in the cone.

\subsection{Product regions}
\label{SectionProducts}

For each partial marking $\mu$ of $S$, define its associated \emph{product region}\index{product region}\index{1aaproductregion@$\QQ(\mu)$, product region of a partial marking}
$$
\QQ(\mu) = \{\mu'\in\MM(S): \mu \subset \mu'\}, 
$$
that is, the set of all (full) markings that extend $\mu$. 

For example, if $\mu$ is a curve system with no transversals then $\QQ(\mu)$ is the collection of all markings whose base curves contain $\mu$. In this case $\QQ(\mu)$ is uniformly quasi-isometric to a left coset of a subgroup stabilizing a certain curve system, just as Dehn twist flats are uniformly quasi-isometric to left cosets of certain maximal rank Dehn twist subgroups, as explained in Section~\ref{SectionDehnTwistFlats}. 

In particular, if $\mu$ is a pants decomposition with no transversals then $\QQ(\mu)$ \emph{is} a Dehn twist flat,\index{Dehn twist!flat} and so our notation is consistent with the notation for Dehn twist flats introduced in Section~\ref{SectionBasicDefs}.

\subsubsection*{Product structure on $\QQ(\mu)$.} An element $\mu' \in\QQ(\mu)$ is specified by choosing, for each component $Y$ of
$\open(\mu)$, a marking on $Y$ which we denote $\mu' \restrict Y$, and which may be identified with $\pi_{\MM(Y)}(\mu')$. Hence $\QQ(\mu)$ is naturally identified with $\MM(\open(\mu))$, which is a product
$$ \prod_{Y\in \abs{\open(\mu)}} \MM(Y)
$$
where $\abs{Z}$\index{1aacomponents@$\abs{Z}$, component set of an essential subsurface} denotes the set of components of $Z$. 

With this in mind, given two partial markings $\mu,\nu$ on $S$, define the \emph{extension of $\mu$ by $\nu$}\index{1aaextensionmarking@$\mu \rfloor \nu$, extension of partial marking by partial marking}\index{extension of a partial marking} to be the partial marking given by:
\begin{align*}
\mu\rfloor\nu  = \mu \union \pi_{\MM(\open(\mu))}(\nu).
\end{align*}
If $\pi_{\MM(\open(\mu))}(\nu)$ is a full marking in $\open(\mu)$ then $\mu\rfloor\nu\in\QQ(\mu)$. This applies for example when $\nu$ itself is a full marking on $S$, in which case item~(\ref{ItemClosest}) of the following proposition tells us that $\mu \rfloor \nu$ is (coarsely) the closest point projection of $\nu$ onto $\QQ(\mu)$.

The following result generalizes the case considered in \BMRank\ 
where $\mu$ was a curve system without transversals.

\begin{proposition}{Q product structure} Let $\xi(S)\ge 1$ and
let $\mu$ be a partial marking of~$S$.
\begin{enumerate}
\item \label{ItemProdQI}
The map
\begin{equation*}
\MM(\open(\mu)) = \prod_{Y\in\abs{\open(\mu)}}\MM(Y) \qquad\to\qquad \MM(S)
\end{equation*}
induced by the identification with $\QQ(\mu)$ is a quasi-isometric embedding, with constants depending only on the topology of $S$.  
\item There is a constant $A_0$ depending only on
the topology of $S$ such that for each $A\ge A_0$, and for each $x \in
\MM(S)$, the minimum distance $d_{\MM(S)}(x,\QQ(\mu))$ from $x$ to $\QQ(\mu)$
in $\MM(S)$ is estimated by 
$$
d_{\MM(S)}(x,\QQ(\mu)) \approx \sum_{Y
\pitchfork \mu} \truncate{d_{\CC(Y)}(x,\mu)}{A}
$$
where the constants of approximation depend on $A$ and on the topology
of $S$.
\item \label{ItemClosest}
Moreover, again with uniform constants, for each $x\in\MM(S)$
$$
d_{\MM(S)}(x,\QQ(\mu)) \approx d_{\MM(S)}(x,\mu\rfloor x).
$$
\end{enumerate}
\end{proposition}

As a consequence of (1), combined with the ordinary quasidistance
formula for components of $S\setminus\base(\mu)$, if we let $\Delta$ be the
subset of $\base(\mu)$ consisting of those curves for which no
transversal is defined then we have a quasi-isometry $$\QQ(\mu)
\,\, \approx \,\, \Z^\Delta \cross \prod_{ \substack{Y \in\abs{\open(\mu)} \\ \xi(Y) \ge
1}} \MM(Y)
$$

As another example of (1), given an essential subsurface $W$, we have
$$\QQ(\bdy W) \approx \MM(W) \cross \MM(W^c)
$$

As an example of (2), if $\mu$ is a full marking then $\QQ(\mu) =
\{\mu\}$, $\open(\mu) = \emptyset$, $\supp(\mu)=S$, and (2) is just
the ordinary quasidistance formula to $\mu$.

\begin{proof}[Proof of Proposition \ref{Q product structure}] Suppose
that $F$ is an essential subsurface of $S$ that is \emph{not}
essentially contained in $\open(\mu)$, and recall that $\mu \pitchfork F$. It follows that for any $\mu',\mu'' \in \QQ(\mu)$ each of $\pi_{\CC(F)}(\mu')$, $\pi_{\CC(F)}(\mu'')$ is within a uniformly bounded distance of $\pi_{\CC(F)}(\mu)$ --- this is a consequence of coarse well-definedness of $\pi_{\CC(F)}$, Lemma \ref{coarse definition}. We can therefore make the $\pi_{\CC(F)}$ term drop out of the quasidistance formula for $d(\mu',\mu'')$ by raising the threshold. The remaining terms can be collected to give the sums of the quasidistance formulas for the projections of $\mu',\mu''$ to the components of $\open(\mu)$, and we obtain (1) as an immediate consequence.

Now consider a connected essential subsurface $Y \subset S$
such that $Y \pitchfork \mu$.  For each $\mu' \in \QQ(\mu)$, from
coarse well-definedness (Lemma \ref{coarse definition}) it follows
that $\pi_{\CC(Y)}(\mu)$ and $\pi_{\CC(Y)}(\mu')$ are within uniformly bounded distance of each other.  Therefore, for each $x \in \MM(S)$ we have $d_{\CC(Y)}(x,\mu) \ceq{1,C} d_{\CC(Y)}(x,\mu')$ where $C$ depends only on the topology of $S$. The terms $d_{\CC(Y)}(x,\mu)$ comprise the right side of (2) and the terms $d_{\CC(Y)}(x,\mu')$ are among the terms in the quasidistance formula for $d_{\MM(S)}(x,\mu')$, so after raising the threshold by $2C$ we obtain the $\gtrsim$ direction of (2).

Next, let $\mu' = \mu\rfloor x$. For each $Y$, if
$\mu \not\pitchfork Y$ then $d_{\CC(Y)}(x,\mu')$ is bounded above by a
constant depending only on the topology of $S$, because by Lemma
\ref{coarse composition} 
$\pi_{\CC(Y)}(\mu')$ is within uniformly bounded distance of
$\pi_{\CC(Y)}(\pi_{\CC(\open(\mu))} (\mu'))$, which equals
$\pi_{\CC(Y)}(\pi_{\CC(\open(\mu))}(x))$,
which is itself within uniformly bounded distance of $\pi_{\CC(Y)}(x)$.
We can then raise the threshold above this constant, so that all of
these terms drop out of the quasidistance formula for $d_{\MM(S)}(x,\mu')$,
leaving only the terms where $\mu\pitchfork Y$. This proves the
$\lesssim$ direction of (2), as well as (3). 
\end{proof}

\subsubsection*{Junctures}\index{junctures}
Let $\FF$ be a family of subsets of a metric space $\MM$.  We say that
$\FF$ \emph{has junctures} if the following holds: for any $X,Y$ in
$\FF$ there exist $E(X,Y)\subset X$ and $E(Y,X)\subset Y$, both
members of $\FF$ as well, such that:
\begin{enumerate}
\item The Hausdorff distance $d_H(E(X,Y),E(Y,X))$ is finite.  \item If
$x\in X, y\in Y$ then
$$
d(x,y)\gtrsim d(x,E(X,Y)) + d_H(E(X,Y),E(Y,X)) + d(E(Y,X),y).
$$
with constants of approximation being uniform over the family $\FF$.
\end{enumerate}
The sets $E(X,Y)$ and $E(Y,X)$ are called the \emph{junctures} of
$X$ and~$Y$.  Note that the junctures are ``parallel'' in the
sense not just of the bound on Hausdorff distance, but the inequality
the other way $d(x,y)\gtrsim d_H(E(X,Y),E(Y,X))$ which by (2) holds for
all $x\in E(X,Y)$ and $y\in E(Y,X)$.

The motivating example of a family having junctures is the family
of geodesics (finite or infinite) in a $\delta$-hyperbolic space.
Here the implicit constants depend on $\delta$.  This example has the
feature that for any $X,Y$ either $E(X,Y)$ and $E(Y,X)$ are points, or
$d_H(E(X,Y),E(Y,X)) \approx 0$.  Junctures for the family
$\QQ(\mu)$ will not have this feature.

Junctures are examples of the general concept of ``coarse intersection'', as we now explain. In a metric space $\MM$, given subsets $A,B \subset \MM$, we say that the \emph{coarse intersection of $A,B$ is well-defined} if there exists $R \ge 0$ such that any two elements of the collection of subsets $\{\NN_r(A) \intersect \NN_r(B) \suchthat r \ge R\}$ have finite Hausdorff distance. For any subset $C \subset \MM$ which has finite Hausdorff distance from any one of these sets, we also say that the coarse intersection of $A$ and $B$ is \emph{represented by} the set $C$. We can define coarse intersection of a finite number of sets in the same way. 

If the family $\F$ of subsets of $\MM$ has junctures then the junctures are representatives of coarse intersections. To prove this, choose $K \ge 1$, $C \ge 0$ to be constants of approximation as in the definition of junctures and let $L = 1 + \frac{2}{K}$. For any $X,Y \in \F$, let $D = d_H(E(X,Y),E(Y,X))$, and suppose that $r \ge L D$. We have $E(X,Y), E(Y,X) \subset \NN_r(X) \intersect \NN_r(Y)$ because $r \ge D$. If $q \in \NN_r(X) \intersect \NN_r(Y)$ then we may choose $x \in X$, $y \in Y$ such that $d(x,q), d(q,y) \le r$, and we have
\begin{align*}
2r &\ge d(x,y) \\
    &\ge K \left( d(x,E(X,Y) + D + d(E(Y,X),y) \right) + C \\
\frac{2r}{K}-D &\ge d(x,E(X,Y)) \\
d(q,E(X,Y)) &\le d(q,x) + d(x,E(X,Y)) \\
   &\le Lr  - D
\end{align*}
and similarly for $d(q,E(Y,X))$ (notice that $Lr-D \ge 0$ because $r \ge LD \ge \frac{D}{L}$). We record this as:

\begin{lemma}{LemmaCoarseIntJunct} If $\MM$ is a metric space and $\F$ is a family of subsets of $\MM$ having junctures then the coarse intersection of any two elements of $\F$ is well-defined and is represented by their junctures. More precisely, there exists $L \ge 1$, depending only on the constants of approximation, such that for all $X,Y \in \F$, letting $D = d_H(E(X,Y),E(Y,X))$, for all $r \ge L D$ we have 
$$E(X,Y) \union E(Y,X) \subset \NN_r(X) \intersect \NN_r(Y) \subset \NN_{L r - D}(E(X,Y)) \intersect \NN_{Lr-D}(E(Y,X))
$$
\qed\end{lemma}

\subsubsection*{Junctures for the family $\{\QQ(\mu)\}$.}

\begin{lemma}{juncture}
The family of subsets $\QQ(\mu) \subset \MM(S)$ has junctures:
for any partial markings $\mu_0,\mu_1$ of~$S$, the junctures for
$\QQ(\mu_0)$ and $\QQ(\mu_1)$ are
\begin{align*}
E_{01} &= E(\mu_0,\mu_1) = \QQ(\mu_0 \rfloor \mu_1) \subset \QQ(\mu_0)
\\
E_{10} &= E(\mu_1,\mu_0) = \QQ(\mu_1 \rfloor \mu_0) \subset \QQ(\mu_1)
\end{align*}
More precisely we have:
\begin{enumerate}
\item \label{ItemOpenWD} The subsurfaces $\open(\mu_0 \rfloor \mu_1)$,
$\open(\mu_1 \rfloor \mu_0)$, and $\open(\mu_0) \essint \open(\mu_1)$
are all isotopic.  Let $\open(\mu_0,\mu_1)$ denote a surface in this
isotopy class.
\item \label{J Hausdorff}
The Hausdorff distance $d_{H}(E_{01},E_{10})$ in $\MM(S)$ is estimated
by
\begin{equation*}
d_{H}(E_{01},E_{10})\approx d_{\MM(\supp(\mu_0,\mu_1))}(\mu_0 \rfloor
\mu_1,\mu_1 \rfloor \mu_0)
\end{equation*}
where we denote $\supp(\mu_0,\mu_1) = \open(\mu_0,\mu_1)^c$.
\item \label{J divergence}
For $x_i \in \QQ(\mu_i)$ we have
\begin{align*}
d_{\MM(S)}(x_0,x_1) &\approx d_{\MM(S)}(x_0,E_{01}) + d_{\MM(S)}(x_1,E_{10}) \\
&+ d_{\MM(\supp(\mu_0,\mu_1))}(\mu_0 \rfloor \mu_1,\mu_1 \rfloor
\mu_0) + d_{\MM(\open(\mu_0,\mu_1))}(x_0,x_1)
\end{align*}
\end{enumerate}
\end{lemma}

\begin{proof} Part \pref{ItemOpenWD} follows (by symmetry) from the general identity
\begin{equation}\label{open intersection}
  \open(\mu \rfloor \nu) \isotopic \open(\mu) \essint \open(\nu)
\end{equation}
for any two partial markings. The subsurface $Z=\open(\mu\rfloor\nu)$ is the maximal essential subsurface that does not
overlap $\mu\rfloor\nu$, hence $Z$ does not overlap $\mu$ so
$Z\esssubset \open(\mu)$. Also, $Z$ does not overlap
$\pi_{\MM(\open(\mu))}(\nu)$, so by Lemma \ref{characterize open projection}, $Z$ does not overlap $\nu$ and therefore $Z\esssubset \open(\nu)$. We conclude that $\open(\mu \rfloor \nu) \esssubset \open(\mu) \essint \open(\nu)$.

Conversely letting $X=\open(\mu)\essint\open(\nu)$, from Lemma  
\ref{characterize open projection} it follows that  $X\subset
\open_{\open(\mu)}(\pi_{\MM(\open(\mu)}(\nu))$ and so $X$ does not overlap
$\mu\rfloor\nu$. We conclude that 
$\open(\mu)\essint\open(\nu)\esssubset\open(\mu\rfloor\nu)$, and
(\ref{open intersection}) follows.

The proofs of \pref{J Hausdorff} and \pref{J divergence} will be
applications of the quasidistance formula.
Note that, now that we know that 
$\mu_0 \rfloor \mu_1$ and $\mu_1 \rfloor \mu_0$ have the same support surface
$\supp(\mu_0,\mu_1)$, 
the distance between these markings in the complex
$\MM(\supp(\mu_0,\mu_1))$ is defined, so that
\pref{J Hausdorff} makes sense. 

To obtain the inequality $\gtrsim$ in \pref{J Hausdorff}, consider
any term in the quasidistance formula for
$d_{\MM(\supp(\mu_0,\mu_1))}(\mu_0 \rfloor \mu_1, \mu_1 \rfloor
\mu_0)$, indexed by $Y \subset \supp(\mu_0,\mu_1)$. 
This term is within uniform distance of 
$d_{\CC(Y)}(x,x')$ for any $x \in E_{01}$ and $x' \in E_{10}$, since $x$
contains $\mu_0\rfloor\mu_1$ and $x'$
contains $\mu_1\rfloor\mu_0$. 
Hence this term contributes to a lower bound for the 
quasidistance formula for 
$d_{\MM(S)}(x,x')$. As before, raising the threshold eliminates the
effect of the additive errors. 

To prove the inequality $\lesssim$ of \pref{J Hausdorff}, note that each
$x\in E_{ij}$ contains $\mu_i\rfloor\mu_j$.  If we replace this part
of $x$ by $\mu_j\rfloor\mu_i$, holding the part $x \restrict
\open(\mu,\nu)$ constant, we obtain a point $x' \in E_{ji}$.  For any $Y
\subset S$ which does \emph{not} index the quasidistance formula for
the right hand side of \pref{J Hausdorff}, the term $d_{\CC(Y)}(x,x')$
is uniformly bounded, as we see by enumerating cases.  If $Y$
essentially intersects $\bdy\open(\mu,\nu)$ then $\pi_{\CC(Y)}(x)$ and
$\pi_{\CC(Y)}(x')$ are uniformly close to
$\pi_{\CC(Y)}(\bdy\open(\mu,\nu))$.  If $Y$ does not essentially
intersect $\bdy\open(\mu,\nu)$ then $Y$ is isotopic into
$\open(\mu,\nu)$ or its complement.  If $Y$ is isotopic into the
complement of $\open(\mu,\nu)$ then either $Y$ is an annulus component
of $\open(\mu,\nu)$ or $Y$ is an index for the right hand side of
\pref{J Hausdorff}.  If $Y$ is an annulus component of
$\open(\mu,\nu)$, or if $Y$ is an essential subsurface of a component
of $\open(\mu,\nu)$, then $\pi_{\CC(Y)}(x)$ and $\pi_{\CC(Y)}(x')$ are
within uniformly bounded distance of the projection of $x \restrict
\open(\mu,\nu) = x' \restrict \open(\mu,\nu)$.  This exhausts all
cases.  By raising the threshold, it follows that $d_{\MM(S)}(x,x')$ reduces to
the right hand side of \pref{J Hausdorff}, proving the inequality
$\lesssim$.

To prove \pref{J divergence}, let $x_i\in \QQ(\mu_i)$ for $i=0,1$.  We
just need to check that each term in the quasidistance formula for
$d_{\MM(S)}(x_0,x_1)$ contributes to one of the four summands on the right hand
side.

The first summand $d_{\MM(S)}(x_0,E_{01})$, by Proposition~\ref{Q product
structure}(2), is estimated by
\begin{equation*}
d_{\MM(S)}(x_0,E_{01}) \approx \sum_{Y\pitchfork \mu_0\rfloor\mu_1} \Tsh
A{d_{\CC(Y)}(x_0,\mu_0\rfloor\mu_1)}
\end{equation*}
However, note that if $Y\pitchfork \mu_0$ then
$d_{\CC(Y)}(x_0,\mu_0\rfloor\mu_1)\approx 1$ since both markings
contain $\mu_0$.  On the other hand, if $Y \pitchfork \mu_0 \rfloor
\mu_1$ and $Y \not\pitchfork \mu_0$ then $Y \esssubset \open(\mu_0)$ and
$Y \pitchfork \mu_1$ by Lemma \ref{characterize open projection}; and
the converse holds as well. Therefore by 
raising the threshold $A$ we get
\begin{equation*}
d_{\MM(S)}(x_0,E_{01}) \approx \sum_{\substack{Y \esssubset \open(\mu_0) \\
Y\pitchfork \mu_1}} \Tsh A{d_{\CC(Y)}(x_0,\mu_0\rfloor\mu_1)}
\end{equation*}
Each term $d_{\CC(Y)}(x_0,\mu_0\rfloor\mu_1)$ is within a uniformly
bounded distance of $d_{\CC(Y)}(x_0,\mu_1)$, by
Proposition~\ref{coarse composition}; it follows that by raising the
threshold above twice this bound, at the cost of another
multiplicative factor of~2, we get
\begin{equation}\label{divergence term 1b}
d_{\MM(S)}(x_0,E_{01}) \approx \sum_{\substack{Y \esssubset \open(\mu_0) \\
Y\pitchfork \mu_1}} \Tsh A{d_{\CC(Y)}(x_0,\mu_1)}
\end{equation}

We obtain a similar expression for the second summand:
\begin{equation}\label{divergence term 2b}
d_{\MM(S)}(x_1,E_{10}) \approx \sum_{\substack{Y \esssubset \open(\mu_1) \\
Y\pitchfork \mu_0}} \Tsh A{d_{\CC(Y)}(x_1,\mu_0)}
\end{equation}

The third summand is given by 
$$
d_{\MM(\supp(\mu_0,\mu_1))}(\mu_0
\rfloor \mu_1,\mu_1 \rfloor \mu_0) \approx \sum_{Y \esssubset
\supp(\mu_0,\mu_1)} \Tsh A{d_{\CC(Y)}(\mu_0 \rfloor \mu_1,\mu_1
\rfloor \mu_0)}
$$
If $Y \esssubset \supp(\mu_0,\mu_1)$ and $Y \not\pitchfork\mu_0$ then $Y
\pitchfork \mu_1$ and both $\pi_{\CC(Y)}(\mu_0 \rfloor \mu_1)$ and
$\pi_{\CC(Y)}(\mu_1 \rfloor \mu_0)$ are within uniformly bounded
distance of $\pi_{\CC(Y)}(\mu_1)$, so these terms may be dropped by
raising the threshold.  Similarly, if $Y \esssubset \supp(\mu_0,\mu_1)$
and $Y \not\pitchfork\mu_1$ then $Y \pitchfork \mu_0$ and these terms
may be dropped.  

If $Y \not\esssubset \supp(\mu_0,\mu_1)$, and if $Y \pitchfork \mu_0$
and $Y \pitchfork \mu_1$, then $Y \pitchfork \bdy \supp(\mu_0,\mu_1)$;
for if not, $Y$ would be isotopic to the complement of
$\supp(\mu_0,\mu_1)$ and so to overlap $\mu_0$ and $\mu_1$, $Y$ would
have to be an annulus isotopic to a boundary curve of
$\supp(\mu_0,\mu_1)$, which is marked by both $\mu_1$ and $\mu_0$. But
in this case, by definition the annulus would be a component of
$\supp(\mu_0,\mu_1)$ so $Y\esssubset\supp(\mu_0,\mu_1)$ after all.  In
this situation both $\pi_{\CC(Y)}(\mu_0 \rfloor \mu_1)$ and
$\pi_{\CC(Y)}(\mu_1 \rfloor \mu_0)$ are within uniformly bounded
distance of $\pi_{\CC(Y)}(\bdy \supp(\mu_0,\mu_1))$, so that
$d_{\CC(Y)}(\mu_0\rfloor\mu_1,\mu_1\rfloor\mu_0)$ is uniformly bounded. Thus although these
terms do not appear in the sum, by raising the threshold we may
formally put them into it with only a bounded change to the estimate.

At this stage, the sum is indexed by the set of all $Y \subset S$ such
that $Y \pitchfork \mu_0$ and $Y \pitchfork \mu_1$.  For such $Y$,
$\pi_{\CC(Y)}(\mu_0 \rfloor \mu_1)$ is within uniformly bounded
distance of $\pi_{\CC(Y)}(\mu_0)$, and $\pi_{\CC(Y)}(\mu_1 \rfloor
\mu_0)$ is within uniformly bounded distance of $\pi_{\CC(Y)}(\mu_1)$,
and so the $d_{\CC(Y)}(\mu_0 \rfloor \mu_1,\mu_1 \rfloor \mu_0)$ is
approximated within a uniform additive error by
$d_{\CC(Y)}(\mu_0,\mu_1)$.  By raising the threshold above twice this
error we obtain
\begin{equation}\label{divergence term 3}
d_{\MM(\supp(\mu_0,\mu_1))}(\mu_0 \rfloor \mu_1,\mu_1 \rfloor \mu_0)
\approx \sum_{\substack{Y\pitchfork \mu_0\\ Y\pitchfork\mu_1}} \Tsh
A{d_{\CC(Y)}(\mu_0,\mu_1)}
\end{equation}

The fourth summand is, by the quasidistance formula in
$\open(\mu_0,\mu_1)$, approximated by
\begin{equation}\label{divergence term 4}
d_{\MM(\open(\mu_0,\mu_1))}(x_0,x_1) \approx \sum_{Y \esssubset
\open(\mu_0,\mu_1)} \Tsh A{d_{\CC(Y)}(x_0,x_1)}
\end{equation}

Now putting these four sums \pref{divergence term 1b},
\pref{divergence term 2b}, \pref{divergence term 3}, \pref{divergence
term 4} together, and recalling that $Y \esssubset \open(\mu_i)$ if and
only if $Y \not\pitchfork \mu_i$, it follows that each $Y \subset S$
appears in exactly one of these four sums.  Moreover, whenever $\mu_i$
appears it can be replaced by $x_i$ with a bounded additive change in
the term.  Raising the threshold above twice the value of this change,
we see that the sum is approximated by the quasidistance formula for
$d_{\MM(S)}(x_0,x_1)$.
\end{proof}

\subsection{Dimension and nonseparation in the asymptotic cone}
\label{SectionNonseparation}

We shall apply Proposition~\ref{Q product structure} and the results of \cite{behrstock-minsky:rank} to compute the dimension of a product region in the asymptotic cone, and combined with Alexander duality we will obtain nonseparation results.

Recall that the \emph{topological dimension}\index{topological dimension} of a topological space $X$ is the least $n$ such that every open cover of $X$ can be refined to an open cover $\U$ having the property that for any subset of $\U$ whose members have nonempty intersection in $X$, the cardinality of the subset is $\le n+1$. The \emph{locally compact topological dimension}\index{locally compact topological dimension} is the least $n$ such that each locally compact subset of $X$ has topological dimension $\le n$. It follows immediately that the \v{C}ech cohomology of each locally compact subset of $X$ vanishes in dimensions above $n$.

Throughout this paper we usually use the word \emph{dimension}\index{dimension} to refer to locally compact topological dimension. Also, given a subset $A \subset X$, we use the term \emph{codimension}\index{codimension} to refer to the locally compact dimension of $X$ minus the locally compact dimension of $A$. Proposition~\ref{CodimensionProp} is useful in some contexts for computing codimension.

The following theorem is the main result of \cite{behrstock-minsky:rank}:

\begin{theorem}{DimensionTheorem}
For each connected, finite type surface $F$ with $\xi(F) \ge 1$, the locally compact topological dimension of $\AM(F)$ is $\xi(F)$.
\end{theorem}

Consider a sequence of partial markings $(\mu_n)$. Since an essential subsurface of $S$ can have only finitely many topological types, the type of $\open(\mu_n)$ is $\omega$-a.e.\ constant, and so the number $r(\open(\mu_n))$ is $\omega$-a.e.\ constant, a number we denote $r(\open(\bar\mu))$.

By combining Theorem~\ref{DimensionTheorem} with Proposition~\ref{Q product structure}, in particular the coarse cartesian product formula for $\QQ(\mu_n)$ given after the statement of the proposition, plus the fact that finite cartesian products commute with rescaled ultralimits, we immediately obtain:

\begin{lemma}{DimensionLemma}
If $(\mu_n)$ is a sequence of partial markings such that the rescaled ultralimit $\QQ_\omega(\bbar \mu)$ of $(\QQ(\mu_n))$ is nonempty in $\AM(S)$, then $\QQ_\omega(\bbar\mu)$ has locally compact topological dimension $r(\open(\bbar\mu))$.
\qed\end{lemma}

Given a topological space $X$ and two subsets $A,B$, we say that \emph{$A$ separates $B$ in $X$} if $B$ has nonempty intersection with at least two components of $X-A$; in particular, $B-A$ is disconnected.

\begin{lemma}{NonseparationLemma}
Let $(\mu_n)$ be as in Lemma~\ref{DimensionLemma}. If $E \subset \AM(S)$ is a connected oriented manifold of dimension $D \ge r(\open(\bbar\mu)) + 2$ then $\QQ_\omega(\bbar\mu)$ does not separate $E$ in $X$.
\end{lemma}

\begin{proof}
Since $\QQ_\omega(\bbar\mu)$ is closed and $E$ is locally compact, the set $E \intersect\QQ_\omega(\bbar\mu)$ is locally compact, and therefore by Lemma~\ref{DimensionLemma} is of topological dimension $\le r(\open(\bbar\mu)) < D-1$. By Alexander Duality, 
$$\widetilde H_0(E-\QQ_\omega(\bbar\mu)) \approx \check{H}^{D-1}(E \intersect \QQ_\omega(\bbar\mu))
$$
using reduced homology on the left hand side and \v{C}ech cohomology on the right hand side, but the right hand side is trivial. 
\end{proof}

\subsection{Cubes and their junctures}
\label{SectionCubejunctures}
Cubes are subsets of $\MM(S)$ obtained by choosing an essential subsurface of $S$ whose components all satisfy $\xi \le 1$, marking the complement of that subsurface, and choosing a geodesic in the marking complex of each component of the subsurface; the cube is parameterized by the product of the chosen geodesics. A special case is a ``Dehn twist $k$-flat''\index{Dehn twist!$k$-flat} in which $S$ has $k$ components, each one an annulus, and the chosen geodesic in the marking complex of each annulus is bi-infinite; see Section~\ref{SectionKFlat} for more details on Dehn twist $k$-flats. In this language, a ``Dehn twist flat''\index{Dehn twist!flat} as previously defined is the same thing as a Dehn twist $\xi$-flat.

Given a connected $V \esssubset S$ with $\xi(V) \le 1$, the marking complex $\MM(V)$ is quasi-isometric to a tree which we denote $\TM(V)$.\index{1aatree@$\TM(V)$, tree QI to $\MM(V)$} When $V$ is an annulus, $\TM(V)$ is isometric to $\reals$, and the Dehn twist about $V$ acts naturally by translation on $\TM(V)$. In the other two cases, where $V$ is a one-holed torus or four-holed sphere, $\TM(V)$ is isometric to the dual tree of the usual modular diagram for $\SL 2 \Z$, on which $\MCG(V)$ acts naturally. Given a geodesic segment $r \subset \TM(V)$ of positive length --- finite, half-infinite, or bi-infinite --- and given an annulus $U \esssubset S$, we say that \emph{$r$ is a twist segment with support $U$}\index{twist!segment} if one of the following holds: $U \isotopic V$ is an annulus; or $V$ is not an annulus, $U \esssubset V$, and $r$ is contained in the axis of the Dehn twist about $U$. 

Consider a subset of $\MM(S)$ formed as follows.  Choose a partial marking~$\mu$ such that the components $W_1,\ldots,W_m$ of $W=\open(\mu)$ satisfy $\xi(W_i) \le 1$.  In each tree $\TM(W_i)$ choose $r_i$ to be a geodesic, finite, half-infinite or bi-infinite (we allow length 0 as well). Let $r=\{r_1,\ldots,r_m\}$. The {\em cube}\index{cube}\index{1aacube@$C(\mu,W,r)$, cube} $C(\mu,W,r)$ is the subset of $\QQ(\mu)$ consisting of markings which, in each $W_i$, restrict to a marking in the geodesic
$r_i$. In other words, under the quasi-isometry
\begin{align*}
\QQ(\mu) &\approx \MM(W_1) \cross \cdots \cross \MM(W_m) \\
   &\approx \TM(W_1) \cross \cdots \cross \TM(W_m) \\
   \intertext{we have}
 C(\mu,W,r) &\approx r_1 \cross \cdots \cross r_m
 \end{align*}

Junctures of cubes can be described in a reasonably
straightforward manner, with careful bookkeeping, in terms of the
description of junctures of product sets given in
Lemma~\ref{juncture}.  Here are the details.

\begin{lemma}{cube junctures}
The family of cubes has junctures. 
\end{lemma}

\begin{proof}
Given cubes $C(\mu,W,r)$ and $C(\nu,V,s)$, we must construct subcubes
which will function as junctures.  Denote the components as
$W=W_1 \union \cdots \union W_m$ and $V = V_1 \union\cdots\union V_n$.

First we describe the essential subsurface $\open(\mu,\nu) = W \essint
V$, whose marking complex parameterizes the junctures of
$\QQ(\mu)$ and $\QQ(\nu)$, by Lemma~\ref{juncture}.

We claim that the components may be reindexed as
\begin{align*}
W &= (W_1 \union \cdots \union W_k) \union (W_{k+1} \union\cdots\union
W_m) \\
V &= (V_1 \union \cdots \union V_k) \union (V_{k+1} \union \cdots
\union V_n)
\end{align*}
where $k \ge 0$, so that the components of $W \essint V$ are $$W
\essint V = \underbrace{(W_1 \essint V_1)}_{U_1} \union \cdots \union
\underbrace{(W_k \essint V_k)}_{U_k}
$$
and so that for each $i = 1,\ldots,k$ one of the following holds:
either $W_i \isotopic V_i \isotopic U_i$; or $U_i$ is an annulus which
is essentially contained in $W_i$ and in $V_i$.

More generally, consider two connected, essential subsurfaces $X,Y$ of
$S$ with $\xi(X), \xi(Y) \le 1$.  If $U = X\essint Y$ is nonempty, it can
only be an annulus or all of $X$ and $Y$. The complement of an annulus
in $X$, if $X$ is not an annulus itself, is either one or two 3-holed
spheres. Now if $Z$ is disjoint from
$Y$, we claim that $X\essint Z$ is empty. For any curve $c$ in
$\Gamma(Z)\intersect\Gamma(X)$ would have to be essential in $X$ and
isotopic to the complement of $Y$ -- hence $U$ would be an annulus and
$c$ isotopic to its core. 
This would make $c$ essential in both $Y$ and $Z$, which
is impossible unless $Y$ and $Z$ are isotopic annuli.

In the context of $W$ and $V$, this implies that the relation $W_i \essint V_j \ne \emptyset$
is a bijection between a subset of the components of $W$ and a subset
of the components of $V$, and the claim immediately follows.

Now we will construct a quasi-isometric  embedding of $\QQ(\mu)$ and
$\QQ(\nu)$ into a product of trees, which will allow us to see their
junctures more clearly. 

For each $i=1,\ldots,k$ the inclusion $U_i \subset W_i$ induces an embedding $\TM(U_i) \inject \TM(W_i)$ whose image is a subtree denoted $\tau_i \subset \TM(W_i)$: either $U_i \isotopic W_i$ and $\tau_i = \TM(W_i)$; or $U_i$ is an annulus and $\tau_i$ is the axis in $\TM(W_i)$ of the Dehn twist about $U_i$.  Similarly, the inclusion $U_i \subset V_i$ induces a quasi-isometric embedding $\MM(U_i) \to \TM(V_i)$ whose image is a subtree $\sigma_i \subset \TM(V_i)$.  

By composing a coarse inverse of the map
$\MM(U_i) \to \tau_i$ with the map $\MM(U_i) \to \sigma_i$, we obtain a
quasi-isometry $g_i \from \tau_i \to \sigma_i$.  Notice that we may
take $g_i$ to be a simplicial isomorphism, as one can verify easily in
either of two cases: if $W_i \isotopic V_i \isotopic U_i$ then these
isotopies induce simplicial isomorphisms of marking complexes; and
otherwise $\tau_i$ and $\sigma_i$ are the axes in the trees $\TM(W_i)$
and $\TM(V_i)$, respectively, of the Dehn twist about $U_i$, and we
can take $g_i$ to be a simplicial isomorphism between these two axes.
Let $X_i$ be the tree obtained from the disjoint union of the trees
$\TM(W_i)$ and $\TM(V_i)$ by gluing $\tau_i$ to $\sigma_i$
isometrically using the map $g_i$.  

Let 
$$\Upsilon = \R \times
\prod_{i=1}^k X_i \times \prod_{i=k+1}^m \TM(W_i) \times
\prod_{i=k+1}^n \TM(V_i).
$$
This is a product of trees on which we can put the $\ell^1$ metric. 

Now for $i=k+1,\ldots,m$, let $p_i =
\pi_{\MM(W_i)}(\nu)$ and 
and note that in fact a bounded
neighborhood of $p_i$ contains all of 
$\pi_{\MM(W_i)}(\QQ(\nu))$. Similarly, for $j=k+1,\ldots,n$ let
$q_j = \pi_{\MM(V_j)}(\mu)$ which approximates
$\pi_{\MM(V_j)}(\QQ(\mu)).$

The product structure of $\QQ(\mu)$ (Proposition \ref{Q product structure}) 
now gives us a quasi-isometric embedding 
$$\theta^\mu \from \QQ(\mu) \to \Upsilon
$$
which is the identity on the $\TM(W_i)$ factors (including
those embedded in the $X_i$), and maps to the constant $q_j$ on each
$\TM(V_j)$, $j=k+1,\ldots,n$, and to 0 in the $\R$ factor. 
Similarly we have 
$$\theta^\nu \from \QQ(\nu) \to \Upsilon,$$ 
which is the identity on the $\TM(V_i)$ factors (including
those embedded in the $X_i$), and maps to the constant $p_j$ on each
$\TM(W_j)$, $j=k+1,\ldots,m$, and to $D$ in the $\R$ factor, where 
$D$ is the Hausdorff distance between $E(\mu,\nu)$ and $E(\nu,\mu)$. 

Note, by Lemma \ref{juncture},  that the images $\theta^\mu(E(\mu,\nu))$
and $\theta^\nu(E(\nu,\mu))$ are parallel products of subtrees, namely
$$
\{0\}\times \prod_1^k\sigma_i \times
\prod_{k+1}^m\{p_i\}\times\prod_{k+1}^n\{q_i\}
$$
and
$$
\{D\}\times \prod_1^k\tau_i \times
\prod_{k+1}^m\{p_i\}\times\prod_{k+1}^n\{q_i\}
$$
recalling that $\sigma_i$ and $\tau_i$ are identified in $X_i$. 

Moreover we note that, by the distance formula (3) in
Lemma~\ref{juncture} (and its interpretation in terms of projections
in (\ref{divergence term 1b},\ref{divergence term 2b},\ref{divergence term
  3},\ref{divergence term 4})), 
$\theta^\mu$ and $\theta^\nu$ actually combine to
give us a quasi-isometric embedding of the union $\QQ(\mu)\union \QQ(\nu)$
into $\Upsilon$, which we will call $\theta$. 

In particular $\theta(E(\mu,\nu))$ and $\theta(E(\nu,\mu))$ are junctures for 
$\theta(\QQ(\mu))$ and
$\theta(\QQ(\nu))$ in this product of trees. This is a special case of
the following easy fact:

\begin{lemma}{junctures in tree products}
Let $T = T_1 \times \cdots\times T_N$ be a product of complete trees with the
$\ell^1$ metric. Then the family of products of closed subtrees has
junctures. Moreover the  approximations in the definition of junctures
are all exact.
\end{lemma}
\begin{proof}
For a single tree this is easily checked: 
Any two subtrees either intersect, in which case the junctures are 
(two copies of) their common subtree, or are disjoint, in which case
the junctures are the unique points of closest approach of each tree
to the other. For a product of subtrees in a product of trees, the
junctures are the products of junctures in the factors, and the
distance formulas in the factors sum to give the desired outcome. 
\end{proof}

Now it is easy to understand how the cubes $C(\mu,W,r)$ and
$C(\nu,V,s)$ are situated by considering their $\theta$-images.
$\theta(C(\mu,W,r))$ is a product of lines and points in the factors of
$\Upsilon$, with first coordinate $0$, and $\theta(C(\nu,V,s))$ is a similar
product with first coordinate $D$. Lemma \ref{junctures in tree
  products} implies that the junctures of the images are again
products of subintervals, and we conclude that the $\theta$-preimages,
which are subcubes of the original cubes, are also junctures.
\end{proof}

The proof of Lemma \ref{cube junctures} gives some more
information about the structure of the junctures of two cubes, which
we record here: 

\begin{lemma}{cube juncture details}
Let $C_1 = C(\mu,W,r)$ and $C_2 = C(\nu,V,s)$. The junctures $C_{ij} = E(C_i,C_j) \subset C_i$ are subcubes of the form $C_{12} = C(\mu,W,r')$ and $C_{21} = C(\nu,V,s')$, where each component of $r'$ or $s'$ is a subinterval or point of the corresponding component of $r$   or $s$. 

After the renumbering in the proof of Lemma \ref{cube junctures}, the components of $r'$ and $s'$ that are not single points come in pairs $r'_i,s'_i$ such that $U_i = W_i\essint V_i \ne \emptyset$, and $r'_i,s'_i$ are images of the same segment of $\TM(U_i)$ under the quasi-isometric embeddings $\TM(U_i) \inject \TM(W_i)$ and $\TM(U_i) \inject \TM(V_i)$. Furthemore:
\begin{enumerate}
\item If $U_i$ is an annulus then $r'_i,s'_i$ are twist segments with support $U_i$.
\item If $U_i \isotopic V_i \isotopic W_i$ then $r'_i = s'_i = r_i \intersect s_i$ in the tree $\TM(U_i) = \TM(V_i) = \TM(W_i)$.
\end{enumerate}
\end{lemma}

\begin{proof}
This is a consequence of the fact that the map $\theta$ in the proof of Lemma~\ref{cube junctures} respects the product structures in its domain and range. The image of $C(\mu,W,r)$  in $\Upsilon$ is a product of line segments in the factors $X_1$, \ldots, $X_k$ and $\TM(W_{k+1})$, \ldots, $\TM(W_m)$, and points in the other factors, whereas $C(\nu,V,s)$ maps to a product of line segments in $X_1,\ldots,X_k$ and in $\TM(V_{k+1}),\ldots, \TM(V_n)$, and points in the rest. Thus, any tree factor in which the juncture factor is a nondegenerate segment is an $X_i$ which corresponds to a pair $W_i,V_i$ that has nontrivial essential intersection $U_i$, this nondegenerate segment is the intersection of the images of $r_i$ and $s_i$ via the quasi-isometric embeddings $\TM(W_i) \inject X_i$ and $\TM(V_i) \inject X_i$, the pullbacks of this segment to $\TM(W_i)$ and to $\TM(V_i)$ are the segments $r'_i, s'_i$ respectively, and the pullbacks of these two segments to $\TM(U_i)$ are the same segment; in the special case when $U_i$ is an annulus item~(1) is an immediate consequence, and when $U_i \isotopic V_i \isotopic W_i$ item~(2) is an immediate consequence. The pullbacks of the junctures by $\theta$ are then subcubes respecting the product structures of the original cubes, and with nondegenerate segments $r'_i, s'_i$ only in the factors $\TM(W_i)$, $\TM(V_i)$ corresponding to the~$X_i$. 
\end{proof}

\subsection{Cubes and junctures in the asymptotic cone}
\label{SectionCubesCone}

From the definition of junctures we can obtain the following statement in the asymptotic cone: Let $\FF$ be a family with junctures in $\MM$, let $(X_n)$ and $(Y_n)$ be sequences in $\FF$ and let $X_\omega$ and $Y_\omega$ be their rescaling ultralimits in the cone $\MM_\omega$.  We find that $E(X,Y)_\omega$ and $E(Y,X)_\omega$ are either disjoint or identical, depending on rate of growth of the Hausdorff distance.  Property~2 in the definition of junctures also implies that
$$X_\omega \intersect Y_\omega = E(X,Y)_\omega \intersect E(Y,X)_\omega,
$$
and hence this intersection is either empty or equal to the limit of
the junctures.

Now given a sequence of cubes $C^n = C(\mu^n,W^n,r^n)$, which we denote $\bbar C = C(\bbar\mu,\bbar W,\bbar r)$, we can take the
cone $ C^\omega(\bbar\mu,\bbar W,\bbar r)$, which is nonempty provided that the distance from the cubes to the basepoint of $\MM(S)$ does not grow too fast. This object has dimension less than or equal to the number of components of $W^n$ for $\omega$-a.e.\ 
$n$.  In fact the limit cube is naturally bilipschitz homeomorphic to
$r^\omega_1\times\cdots\times r^\omega_k$ where each $r^\omega_i$ is an embedded path in the $\reals$-tree $\MM_\omega(\bbar W_i)$ whose length is in $[0,\infinity]$, having positive length if and only if the length of the sequence $(r^n_i)$ grows linearly.  We will continue calling these objects cubes.

Lemma \ref{cube junctures} on junctures for cubes implies, using
the discussion in the beginning of this section, that the intersection
of two cubes in $\MM(S)$ is empty or is a cube, possibly a trivial cube, meaning a single point.  Moreover this cube is described by data closely related to the original cubes.  We will use this in 
Section~\ref{orthant defs} to understand the {\em complex of orthants} in the asymptotic cone.

\subsection{Coarse set theory of Dehn twist $k$-flats}
\label{SectionKFlat}

We have already explained in Lemma~\ref{LemmaCoarseIntJunct} how
junctures are examples of coarse intersection.  In this section we
make a further study of coarse inclusion and coarse equivalence among
Dehn twist $k$-flats in $\MM(S)$.

A \emph{Dehn twist $k$-flat}\index{Dehn twist!$k$-flat} in $\MM(S)$,
$0 \le k \le \xi(S)$, is a subset of the form $\QQ(\mu)$ where $\mu$
is a marking such that $\base(\mu)$ is a pants decomposition with
exactly $k$ unmarked components.  As the terminology suggests, each
Dehn twist $k$-flat is a quasi-isometrically embedded copy of
$\reals^k$ in $\MM(S)$, by Proposition~\ref{Q product structure}~(1).

Throughout the paper, the phrase ``Dehn twist flat'',\index{Dehn
twist!flat} when unadorned by a dimension, will by default refer to
the top dimensional case, namely a Dehn twist $\xi(S)$-flat.
Sometimes we emphasize this by referring to \emph{maximal Dehn twist
flats}.\index{maximal Dehn twist flat}

In a metric space $\MM$, given two subsets $A,B \subset \MM$, we say
that $A$ is \emph{coarsely included} in $B$\index{coarse!inclusion} if
there exists $r \in [0,\infinity)$ such that $A \subset \NN_r(B)$.
We say that $A$ is \emph{coarsely equivalent} to
$B$\index{coarse!equivalence} if there exists $r \in [0,\infinity)$
such that $A \subset \NN_r(B)$ and $B \subset \NN_r(A)$; the infimum
of all such $r \in [0,\infinity]$ is equal to the Hausdorff
distance\index{Hausdorff distance} between $A$ and $B$.

Given a $k-1$ simplex $\sigma$ in $\CC(S)$ and a Dehn twist $k$-flat
$\QQ(\mu)$, we say that $\sigma$ is \emph{represented by $\QQ(\mu)$}
if the system of unmarked curves of $\base(\mu)$ is isotopic
to~$\sigma$.  In particular, vertices are represented by Dehn twist
1-flats and edges by Dehn twist 2-flats.  Note that each Dehn
twist $k$-flat represents a unique $k-1$ simplex. Conversely, each $k-1$ simplex is represented by infinitely many Dehn
twist $k$-flats, except in the case $k=\xi$ where each maximal dimension simplex is represented by a unique maximal Dehn twist flat.

\begin{lemma}{LemmaFlatSetTh}
Given simplices $\sigma_i \subset \CC(S)$, $i=0,1$, and representative
Dehn twist $k_i$-flats $\QQ(\mu_i)$ respectively, we have:
\begin{enumerate}
\item \label{ItemContCCont}
$\QQ(\mu_0)$ is coarsely contained in $\QQ(\mu_1)$ if and only if
$\sigma_0 \subset \sigma_1$.
\item \label{ItemEqCEq}
$\QQ(\mu_0)$ is coarsely equivalent to $\QQ(\mu_1)$ if and only if
$\sigma_0 = \sigma_1$.
\item \label{ItemIntCInt}
Given a simplex $\tau \subset \CC(S)$ and a representative Dehn twist
$l$-flat $\QQ(\nu)$, $\QQ(\nu)$ represents the coarse intersection of
$\QQ(\mu_0)$ and $\QQ(\mu_1)$ if and only if $\tau = \sigma_0
\intersect \sigma_1$.
\end{enumerate}
\end{lemma}

\begin{proof} First note that $\sigma_0 \subset \sigma_1$ if and only
if $\QQ(\mu_0)$ is equal to $E(\QQ(\mu_0),\QQ(\mu_1)) = \QQ(\mu_0
\rfloor \mu_1)$.  The ``if'' direction of item \pref{ItemContCCont} is
then an immediate consequence of Lemma~\ref{juncture} \pref{J
Hausdorff}.  Conversely, if $\sigma_0 \not\subset \sigma_1$ then by
applying Lemma~\ref{Q product structure} \pref{ItemProdQI} it follows
that there exist points of $\QQ(\mu_0)$ which are arbitrarily far from
$\QQ(\mu_0 \rfloor \mu_1)$, and then Lemma~\ref{juncture} \pref{J
divergence} proves that $\QQ(\mu_0)$ is not contained in any finite
radius neighborhood of $\QQ(\mu_1)$.

Item~\pref{ItemEqCEq} follows by symmetric applications of
item~\pref{ItemContCCont}.

To prove item~\pref{ItemIntCInt}, first apply Lemma~\ref{juncture} to
conclude that $\QQ(\mu_0 \rfloor \mu_1)$ and $\QQ(\mu_1 \rfloor
\mu_0)$ are the junctures of $\QQ(\mu_0)$ and $\QQ(\mu_1)$, then apply
Lemma~\ref{LemmaCoarseIntJunct} to conclude that these two junctures
both represent the coarse intersection of $\QQ(\mu_0)$ and
$\QQ(\mu_1)$, and then note that both of these junctures represent the
simplex $\sigma_0 \intersect \sigma_1$.  Now apply
item~\pref{ItemEqCEq} to $Q(\nu)$ and either of the two junctures.
\end{proof}

The proof of Lemma~\ref{LemmaFlatSetTh} \pref{ItemIntCInt} does not
make full use of the uniform control on coarse intersection that is
provided by Lemma~\ref{LemmaCoarseIntJunct}.  The following lemma make
use of this control, which will be needed in
Section~\ref{SectionQIClassificationProof}.

\begin{lemma}{LemmaFlatIntersection} For each Dehn twist $k$-flat $\QQ(\mu)$ there exist two maximal Dehn twist flats $\QQ(\mu_0)$, $\QQ(\mu_1)$ whose junctures $\QQ(\mu_0 \rfloor \mu_1)$ and $\QQ(\mu_1 \rfloor \mu_0)$ are Dehn twist $k$-flats at uniform Hausdorff distance from each other and from~$\QQ(\mu)$, and which uniformly represent the coarse intersection of $\QQ(\mu_0)$ and $\QQ(\mu_1)$. To be precise: 
\begin{align*}
\QQ(\mu) \union \QQ(\mu_0 \rfloor \mu_1) \union \QQ(\mu_1 \rfloor \mu_0) &\subset \NN_R(\QQ(\mu_0)) \intersect \NN_R(\QQ(\mu_1)) \\
 &\subset \NN_C(\QQ(\mu)) \intersect \NN_C(\QQ(\mu_0 \rfloor \mu_1)) \intersect \NN_C(\QQ(\mu_1 \rfloor \mu_0))
\end{align*}
where the constants $C,R \ge 0$ depend only on the topology of $S$.
\end{lemma}

\begin{proof}[Proof of Lemma \ref{LemmaFlatIntersection}] It suffices to prove that the three sets $\QQ(\mu)$, $\QQ(\mu_0 \rfloor \mu_1)$, $\QQ(\mu_1 \rfloor \mu_0)$ are all at uniformly finite Hausdorff dimension from each other, for once this is done we may apply Lemma~\ref{juncture} to obtain that $\QQ(\mu_0 \rfloor \mu_1)$, $\QQ(\mu_1 \rfloor \mu_0)$ are the junctures of $\QQ(\mu_0)$, $\QQ(\mu_1)$, and then we may apply Lemma~\ref{LemmaCoarseIntJunct} to obtain the desired uniform control on coarse intersection.

We first prove the lemma in the special case that $k=0$, so $\mu$ is a marking and $\QQ(\mu)=\{\mu\}$. Note that a pants decomposition $\Delta$  overlaps each component of $\base(\mu)$ if and only if $\base(\mu)$ overlaps each component of~$\Delta$, in which case both $\base(\mu) \rfloor \Delta$ and $\Delta \rfloor \base(\mu)$ are markings of $S$. We shall find such a $\Delta$ so that each of the markings $\base(\mu) \rfloor \Delta$ and $\Delta \rfloor \base(\mu)$ is uniformly close to $\mu$ in $\MM(S)$. 

To find the appropriate $\Delta$, letting $\base(\mu) = \{c_1,\ldots,c_\xi\}$, we shall build up $\Delta = \{d_1,\ldots,d_\xi\}$ one component at a time. Proceeding by induction, choose the subset $\{d_1,\ldots,d_k\}$ so that:
\begin{itemize}
\item $\{d_1,\ldots,d_k,c_{k+1},\ldots,c_\xi\}$ is a pants decomposition.
\item $c_k,d_k$ are connected by an edge in the curve complex of the complexity~$1$ component of $S-(d_1 \union\cdots\union d_{k-1} \union c_{k+1} \union\cdots\union c_\xi)$ that contains them.
\item $\pi_{\CC(c_k)}(d_k)$ is a uniform distance from the $\mu$-transversal of $c_k$.
\end{itemize}
The last item is possible because each orbit of the action on $\CC(c_k)$ of the Dehn twist group $\<\tau_{c_k}\>$ comes uniformly close to each point of $\CC(c_k)$. Application of Proposition~\ref{Q product structure} \pref{ItemProdQI} provides a uniform bound to $d(\mu,\base(\mu) \rfloor \Delta)$. Consider the sequence of markings defined by
\begin{align*}
\mu_0 &= \base(\mu) \rfloor \Delta \\ &=  ((c_1,d_1),\ldots,(c_k,d_k)) \\
\mu_i &= ((d_1,c_1),\ldots,(d_{i},c_{i}),(c_{i+1},d_{i+1}),\ldots,(c_k,d_k)), \quad i=1,\ldots,k-1 \\
\mu_k  &= \Delta \rfloor \base(\mu) \\
&= ((d_1,c_1),\ldots,(d_k,c_k)) \\
\end{align*}
The markings $\mu_i,\mu_{i+1}$ are connected by an edge in $\MM(S)$ and so 
$$d(\base(\mu) \rfloor \Delta,\Delta \rfloor \base(\mu)) \le k
$$

Now we reduce the general case of the lemma to the special case just proved. Consider any Dehn twist $k$-flat $\QQ(\mu)$ with $k \le \xi-1$. Write $\mu = \mu' \union \mu''$ where $\mu'$ consists of the $k$ unmarked curves of $\base(\mu)$ and $\mu''$ consists of the $\xi-k$ marked curves together with their transversals. The subsurface $F = \supp(\mu)$ is the union of the components of $S-\mu'$ that are not 3-holed spheres, and $\mu''$ is a marking of $F$. Applying the special case we obtain a pants decomposition $\Delta$ of $F$ such that $\base(\mu'') \rfloor \Delta$ and $\Delta \rfloor \base(\mu'')$ are markings of $F$ each at uniformly bounded distance from $\mu''$ in $\MM(F)$. Let $\mu_0 = \base(\mu) = \mu' \union \base(\mu'')$ and $\mu_1 = \mu' \union \Delta$. By construction the sets $\QQ(\mu)$, $\QQ(\mu_0 \rfloor \mu_1)$, $\QQ(\mu_1 \rfloor \mu_0)$ are all their own junctures among each other, and Lemma~\ref{juncture} \pref{J Hausdorff} provides a uniform bound to their Hausdorff distances.
\end{proof}

\section{Consistency theorem}
\label{consistency}

\newcommand\wprec{\mathrel{\ll}}

In this section we will derive a coarse characterization of the image
of the \emph{curve complex projections map}\index{curve complex
projections map}\index{1aacurvecomplexprojection@$\Pi$, curve complex
projections map}
$$
\Pi\co\MM(S) \to \prod_{W\subseteq S} \CC(W)
$$
defined by $\Pi(\mu) = (\pi_W(\mu))_W$.

Consider the following {\em Consistency Conditions}\index{consistency
conditions} on a tuple $(x_W)\in\prod_W \CC(W)$, where $c_1$ and $c_2$
are a pair of positive numbers:
\begin{itemize}
\item[C1:] Whenever $W\pitchfork V$,
$$
\min \left( d_W(x_W,\boundary V) , d_V(x_V,\boundary W) \right) < c_1.
$$
\item[C2:] Whenever $V\esssubset W$ and $d_W(x_W,\boundary V) > c_2$,
$$
d_V(x_V,x_W) < c_1.
$$
\end{itemize}
In C2 recall that $d_V(x_V,x_W)$ is shorthand for
$d_V(x_V,\pi_V(x_W))$ and so is only defined if $x_W$ is in the domain
of $\pi_V$, which it may not be if $d_W(x_W,\bdy V)$ is too small.
Notice that if C1 and C2 are satisfied with respect to positive
constants $c_1,c_2$, then they are satisfied with respect to any
larger constants.

The conditions C1-C2, with suitable constants, are satisfied by the
image of~$\Pi$, and moreover
\begin{lemma}{consistency necessary}
Given $K$ there exist $c_1,c_2\ge 1$ such that, if $\mu\in\MM(S)$ and
$(x_W)\in\prod\CC(W)$ such that $d_W(x_W,\mu) \le K$ for all
$W\subseteq S$, then $(x_W)$ satisfies C1 and C2 with constants
$c_1,c_2$.
\end{lemma}

\begin{proof}
The case $K=0$, i.e. $(x_W) = \Pi(\mu)$, follows from Behrstock's
inequality \cite{behrstock:asymptotic}, namely
\begin{lemma}{behrstock inequality}
There exists $m_0$ such that for any marking $\mu\in\MM(S)$ and
subsurfaces $V\pitchfork W$,
$$
\min \left( d_W(\mu,\boundary V) , d_V(\mu,\boundary W) \right) < m_0.
$$
\end{lemma}
This gives condition C1.  Condition C2, with $c_2=1$ and suitable
$c_1$, follows simply because $\pi_V$ is determined by intersections,
so whenever $V\esssubset W$, $\pi_V\circ\pi_W$ is a bounded distance
from $\pi_V$ when both are defined (Lemma \ref{coarse composition}).
 
For $K>0$ we simply observe that (C1-2) are preserved, with suitable
change in constants, when all the coordinates of $(x_W)$ are changed a
bounded amount.
\end{proof}

Our main point here is to show that conditions (C1-2) are also {\em
sufficient} for a point to be close to the image of $\Pi$, namely:

\begin{theorem}{consistency suffices}
Given $c_1$ and $c_2$ there exists $c_3$ such that, if (C1-2) hold
with $c_1$ and $c_2$ for a point $(x_W)$, then there exists
$\mu\in\MM(S)$ such that
$$
d_W(x_W,\mu) < c_3
$$
for all $W\subseteq S$.
\end{theorem}

The proof of this theorem will take up the rest of
Section~\ref{consistency}.

\subsection{Subsurface ordering induced by projections}
\label{SectionProjOrd}

In order to approach the proof of Theorem~\ref{consistency suffices}
we will first study more carefully the structure imposed by (C1) and
(C2).  Recall from \cite{masur-minsky:complex2} that for any two
markings of $S$, there is a natural partial order on the set of
component domains of subsurfaces that occur in a hierarchy between
those two markings.  With this motivation, the structure that we study
will involve similar partial orders on the collection of proper,
connected, essential subsurfaces (up to isotopy).

Let us fix a tuple $(x_W)$ satisfying (C1-2).  Without loss of
generality, we will assume $c_1>\max\{c_2,m_0,B\}$, where $m_0$ is the
constant given by Lemma~\ref{behrstock inequality} and $B$ is the
constant given by Theorem~\ref{bounded geodesic projection}.
 
If $W,V$ are proper, connected, essential subsurfaces of $S$ and
$k\in\N$, define a relation\index{1aaprojectionorderstrong@$W \prec_k
V$, projection ordering, strong}
$$ W\prec_k V
$$
to mean that
$$ W\pitchfork V \ \text{and}\ d_W(x_W,\boundary V) \ge k(c_1+4).
$$
The role of 4 here is that it is twice the maximal diameter of
$\pi_{\CC(Y)}(\gamma)$ for a curve system $\gamma$ --- see
lemmas~\ref{coarse definition} and~\ref{coarse lipschitz} and the
comments thereafter.

We also allow the right hand side of $\prec_k$ to be a marking $\rho$:
define $W \prec_k \rho$ to mean that $d_W(x_W,\rho) \ge k(c_1+4)$.

Although $\prec_k$ is not quite an order relation on the set of
proper, connected, essential subsurfaces, the family of all $\prec_k$
behaves roughly like a partial order in a way we shall now explore.
Let us also define a relation\index{1aaprojectionorderweak@$W \wprec_k
\rho$, projection ordering, weak} $$ W \wprec_k \rho,
$$
where $\rho$ is any partial marking, to mean that $$ W\pitchfork \rho
\ \text{and}\ d_W(x_W,\rho) \ge k(c_1 +4).
$$
Notice that if $\rho$ is a marking then $W \prec_k \rho$ and $W
\wprec_k \rho$ are equivalent, since $W$ and $\rho$ always overlap.

We then define $W\wprec_k V$ to mean $W\wprec_k \boundary
V$.\index{1aaprojectionorderweak@$W \wprec_k V$, projection ordering,
weak} This is a weaker relation than $W \prec_k V$ because
$W\pitchfork \boundary V$ allows the possibility that $V \esssubset
W$, which cannot happen if $W \pitchfork V$.

Note that if $k\ge p$ then $U\wprec_k V \implies U\wprec_p V$ and
$U\prec_k V \implies U\prec_p V$.  Next we point out that property C1
implies the following:
\begin{itemize}
\item $\prec_k$ is antisymmetric in the following sense: if $U\prec_k
V$ holds, then $V\not\prec_1 U$, and hence $V\not\prec_k U$.
\end{itemize}
Of course $\wprec_k$ is antisymmetric as well, since containment is
already antisymmetric.  Now we will prove the following lemma, which
states that the system of relations is transitive in a certain sense.

\begin{lemma}{order properties}
Given an integer $k>1$ we have:
\begin{enumerate}
\item If $U\prec_k V$ and $V\wprec_2 W$ then $U\prec_{k-1} W$.

Also, if $\rho$ is a marking, and if $U \prec_k V$ and $V \wprec_2
\rho$, then $U \prec_{k-1} \rho$.

\item If $U\wprec_k V$ and $V\wprec_2 W$ then $U\wprec_{k-1} W$.

Also, if $\rho$ is a marking, and if $U \wprec_k V$ and $V \wprec_2
\rho$ then $U \wprec_{k-1} \rho$.
 
\item If $U\pitchfork V$ and both $U\wprec_k \rho$ and $V\wprec_k
\rho$ for some partial marking $\rho$, then $U$ and $V$ are
$\prec_{k-1}$-ordered --- that is, either $U\prec_{k-1} V$ or
$V\prec_{k-1} U$.
\end{enumerate}
\end{lemma}

Note that the weak transitivity of parts (1) and (2) tends to
``decay'' ($k$ decreases) each time it is applied, and hence does not
give a partial order.  However part (3) can be used to re-strengthen
the inequalities under appropriate circumstances.

\begin{proof}
Beginning with (1), suppose $U\prec_k V$ and $V\wprec_2 W$.  From
$V\wprec_2 W$ we have
$$ d_V(x_V,\boundary W) \ge 2(c_1+4)
$$
and from $U\prec_k V$ and property C1 we have
$$ d_V(x_V,\boundary U) < c_1.
$$
By the triangle inequality, together with the fact that
$\diam\pi_{\CC(V)}(\gamma) \le 2$ for any disjoint curve system
$\gamma$ (see comments after Lemma~\ref{coarse lipschitz}),
\begin{align*}
d_V(\boundary U,\boundary W) & \geq d_V(x_V,\boundary W) -
d_V(x_V,\boundary U) -\diam_V(\boundary U) - \diam_V(\boundary W) \\
& > 2(c_1+4) - c_1 - 4 = c_1+4.
\end{align*}
In particular $d_V(\boundary U,\boundary W)>2$, so $\boundary U
\pitchfork \bdy W$, and so $U\pitchfork W$.  Now applying Lemma
\ref{behrstock inequality} we also get
$$ d_U(\boundary V,\boundary W) < m_0 \le c_1
$$
and hence, using $U\prec_k V$ and the triangle inequality as above,
$$ d_U(x_U,\boundary W) > (k-1)(c_1+4).
$$
Hence, $U\prec_{k-1} W$, as desired.

Replacing $W$ (and $\bdy W$) by a marking $\rho$, the second clause of
(1) is proved similarly, noting that the question of whether $U
\pitchfork \rho$ is not an issue.

\medskip

The proof of part (2) is similar.  The case not covered by part (1) is
when $V \esssubset U$, and we consider two cases, depending on whether
$W \esssubset V$ or $W \pitchfork V$.

In the first case suppose that $W \esssubset V$.  Since $V
\not\isotopic U$, it follows that $\boundary V$ and $\boundary W$
together form a curve system in $U$, and hence $\diam_U(\boundary
V\union\boundary W) \le 1$.  So by the triangle inequality we have
$$ d_U(x_U,\boundary W) \ge d_U(x_U,\boundary V) - 1
\ge (k-1)(c_1+4)
$$
and we conclude $U\wprec_{k-1} W$.

In the second case suppose that $W\pitchfork V$.  Since $V\wprec_2 W$
we have
$$ d_V(x_V,\boundary W) \ge 2(c_1+4).
$$
Since $V\esssubset U$ we know that $\boundary W\pitchfork U$.  Since
$d_U(x_U,\boundary V) \ge k(c_1+4)>c_2$, by property C2 we have that
$$ d_V(x_V,x_U) < c_1
$$
and hence
$$ d_V(x_U,\boundary W)
\ge d_V(x_V,\boundary W) - d_V(x_V,x_U) -\diam_V(\boundary W) > c_1+4.
$$
But now by Theorem~\ref{bounded geodesic projection}, this implies
that any $\CC(U)$-geodesic $[x_U,\pi_U(\boundary W)]$ must pass within
distance 1 of $\boundary V$, and we conclude
$$ d_U(x_U,\boundary W) \ge d_U(x_U,\boundary V) - 1 -
\diam_U(\boundary V)
 \ge (k-1)(c_1+4)
$$
and again we have $U\wprec_{k-1} W$.

Again replacing $W$ (and $\bdy W$) by a marking $\rho$, the second
clause of (2) is proved similarly, only the second case of the proof
being relevant.

\medskip

Now we prove (3): starting with $U\pitchfork V$ and
$$ d_U(x_U,\rho) \ge k(c_1+4)
$$
and
$$ d_V(x_V,\rho) \ge k(c_1+4),
$$
suppose $U\not\prec_{k-1} V$, so that
$$ d_U(x_U,\boundary V) < (k-1)(c_1+4).
$$
Then by the triangle inequality
$$ d_U(\rho,\boundary V) > c_1
$$
and by Lemma \ref{behrstock inequality}
$$ d_V(\boundary U,\rho) < m_0.
$$
Now by the triangle inequality
$$ d_V(x_V,\boundary U) > (k-1)(c_1+4)
$$
so $V\prec_{k-1} U$, and we are done.
\end{proof}

Let us now define\index{1aaprojectionlarge@$\FF_k(\rho)$, }
$$ \FF_k((x_{W}),\rho) = \{W\subsetneq S: W\pitchfork \rho
\ \text{and}\ d_W(x_W,\rho) \ge k(c_1 +4) \}.
$$
which is the collection of subsurfaces in whose curve complexes $x_W$
and $\rho$ have large distance. If $x_{W}=\pi_{W}(x)$ for 
$x\in\MM(S)$, then we use the notation $\FF_{k}(x,\rho)$. 

Note that 
$\FF_k((x_{W}),\rho) = \{W\subsetneq S: W \wprec_k \rho\}$, and so we 
sometimes simplify the notation by simply writing $\FF_k(\rho)$.
In this vein, when $Z$ is a subsurface, we let
$\FF_k(Z)$ denote $\FF_k(\partial Z)$.

As a corollary of the previous lemma we obtain:
\begin{lemma}{F order}
If $k>2$ then the relation $\prec_{k-1}$ is a partial order on
$\FF_k((x_{W}),\rho)$.
\end{lemma}
\begin{proof}
All that is needed is to prove that $\prec_{k-1}$ is transitive on
$\FF_k(\rho)$ --- antisymmetry is already established.

Suppose $U,V,W \in \FF_k(\rho)$, and $U\prec_{k-1} V$ and $V\prec_{k-1}
W$.  By Lemma \ref{order properties} part (1), this implies
$U\prec_{k-2} W$.  In particular $U\pitchfork W$, so by Lemma
\ref{order properties} part (3) $U$ and $W$ are $\prec_{k-1}$-ordered.
Antisymmetry together with $U\prec_{k-2} W$ implies that $U\prec_{k-1}
W$, as desired.
\end{proof}

We can also obtain a finiteness statement:
\begin{lemma}{F finite}
If $k>2$ then $\FF_k((x_{W}),\rho)$ is finite.
\end{lemma}

\begin{proof}
Suppose that $\FF_k(\rho)$ is infinite and let $\{Y_i\}$ be an infinite,
injective sequence within it.  After extracting a subsequence we may
assume that $\boundary Y_i \to \lambda$ in $\PML(S)$, the projective
measured lamination space of $S$.  Let $U$ be a subsurface filled by a
component of $\lambda$; possibly $U=S$.  Then $\boundary Y_i$ meets
$U$ for all sufficiently large $i$, and $\pi_U(\boundary
Y_i)\to\infty$ in $\CC(U)$ --- that is, $d_U(\boundary
Y_i,q)\to\infty$ for any fixed $q$.  This is a consequence of the
Kobayashi/Luo argument that $\CC(U)$ has infinite diameter, see \cite
[Proposition 3.6] {masur-minsky:complex1}.  Note in the special case
that $U$ is an annulus we are obtaining that the twisting of
$\boundary Y_i$ around $U$ is going to~$\infty$.

Now, $d_U(x_U,\boundary Y_i) \to \infty$ means for any given $p$ that
eventually $U\wprec_p Y_i$.  However we have $Y_i\wprec_k \rho$ by
assumption, so
$$
U\wprec_{p-1} \rho
$$
by Lemma \ref{order properties}, part (2).  However $U$ and $\rho$ are
fixed and $p$ is arbitrary, so this is impossible.  We conclude that
$\FF_k(\rho)$ is finite.
\end{proof}

\subsection{Proof of the consistency theorem}
As stated at the outset of Section~\ref{SectionProjOrd} we have
$(x_W)$ satisfying C1-2 with the same assumptions on $c_1$ and $c_2$.
We will construct $\mu$ by induction.

Consider  $\FF_3(x_S)$, which we recall is shorthand 
for $\FF_3((x_{W}),x_S)$.  If $\FF_3(x_S) = \emptyset$, let $\mu_0 =
x_S$.  Otherwise, by Lemmas \ref{F finite} and \ref{F order}, the set
$\FF_3(x_S)$ is finite and partially ordered by $\prec_2$, and so this
partial order contains minimal elements.  Among these minimal
elements, choose one, $Y$, of maximal complexity $\xi(Y)$, and let
$\mu_0 = \boundary Y$.

Now consider any $Z\esssubset S$ which overlaps $\mu_0$.  We claim
that
\begin{equation}\label{Z hits mu0 bound}
d_Z(x_Z,\mu_0) < 4(c_1+4).
\end{equation}
Suppose otherwise, so $Z \wprec_4 \mu_0$.  If $\FF_3(x_S)=\emptyset$
then $Z\wprec_4 x_S$ which implies $Z\in\FF_4(x_S)\subseteq
\FF_3(x_S)$, a contradiction.  When $\FF_3(x_S)\ne\emptyset$, we would
have $Z\wprec_4 Y \wprec_3 x_S$, and by Lemma \ref{order properties}
part (2), $Z\wprec_3 x_S$.  Hence $Z\in \FF_3(x_S)$.

Now since $Y$ was $\prec_2$-minimal, we can't have $Z\prec_2 Y$ and we
conclude $Y \esssubset Z$.  Now $Z$ cannot be $\prec_2$-minimal
because its complexity is larger than that of $Y$, so there must be
$V\in\FF_3(x_S)$ with $V\prec_2 Z$.  But then Lemma \ref{order
properties} part (1) implies $V\prec_1 Y$.  In particular $V
\pitchfork Y$ and so, arguing as in the proof of Lemma \ref{F order},
we apply Lemma~\ref{order properties} part (3) with $\rho = x_{S}$ to
conclude that $V$ and $Y$ are $\prec_2$ ordered, and since $V \prec_1
Y$ it follows by asymmetry that $V \prec_2 Y$.  Again this is a
contradiction.  We conclude that (\ref{Z hits mu0 bound}) holds.

Now consider the restriction of $(x_W)$ to subsurfaces in $S\setminus
\mu_0$.  In each component $V$ of $S\setminus \mu_0$, the assumptions
on $(x_W)$ still hold, so inductively there is a marking $\mu_V$ in
$\MM(V)$ satisfying
\begin{equation}\label{Z complement mu0 bound}
d_Z(\mu_V,x_Z) < c_3(V)
\end{equation}
for all $Z\subseteq V$.  We append the $\mu_V$ to $\mu_0$ to obtain a
marking $\mu'$ which almost fills the surface except that it has no
transversal data on the curves of $\mu_0$.  By (\ref{Z hits mu0
bound}) and (\ref{Z complement mu0 bound}), it satisfies a bound on
$d_Z(x_Z,\mu')$ for every $Z\subseteq S$ except the annuli whose cores
are components of $\mu_0$.  Let $\mu$ be the enlargement of $\mu'$
obtained by setting the transversal on each $\gamma\in\mu_0$ to be
$x_\gamma$.  Now we obtain a bound on $d_Z(x_Z,\mu)$ for all $Z$, so
$\mu$ is the desired marking and the proof is complete.

\section{$\Sigma$-hulls}
\label{hulls}

A $\Sigma$-hull of a finite set in $\MM(S)$ (and then, taking
limits, of a finite set in $\AM(S)$) is a substitute for convex hull
which is well adapted to the presence in $\MM(S)$ of both
hyperbolicity and product structure. In ``hyperbolic directions'' it
looks like a hyperbolic convex hull, and in product regions the hull
of two points can be a rectangle. In general it is a hybrid of these. 

In this section we focus on the ``coarse $\Sigma$-hull'' of a finite set $A \subset \MM(S)$, a parameterized family of sets $\Sigma_\epsilon(A)$ which are coarsely well-defined for sufficiently large $\epsilon$ (see Lemma~\ref{coarse Sigma properties}~(3)). Our main goal is Proposition~\ref{hull retraction}, in which we show that $\Sigma$-hulls admit coarse retractions. 

In Section~\ref{contractibility} we will apply this to $\Sigma$-hulls in the asymptotic cone, showing that they are contractible and vary continuously with their extreme points.

\subsection{Hulls in hyperbolic spaces}
\label{S:hyperbolic hulls}
If $A\subset X$ is a finite subset of a $\delta$-hyperbolic geodesic space $X$, let
$\hull_X(A)$\index{1aahull@$\hull_X$, hyperbolic hull} denote the union of geodesics $[a,a']$ with $a,a'\in A$, the \emph{hyperbolic hull}\index{hyperbolic hull} of $A$ in~$X$. We will need the following properties of this construction, which are easy exercises.

\begin{lemma}{hyperbolic hulls}
The sets $\hull_X(A)$ satisfy the following properties, with implicit constants depending only on the hyperbolicity constant of $X$ and the cardinality of $A$:
\begin{enumerate}
\item $\hull_X(A)$ is quasi-convex. 
\item If $x\in X$ and $y\in \hull_X(A)$ is the point nearest to $x$, and if $y'\in \hull_X(A)$, then $d(y,y') \prec d(x,y') - d(x,y)$.   
\item The map $A \mapsto \hull_X(A)$ is coarsely Lipschitz in the Hausdorff metric.
\end{enumerate}
Also, for points $x \in X$ and closed, quasiconvex subsets $B \subset X$:
\begin{enumeratecontinue}
\item The nearest point retraction, which takes the pair $(x,B)$ to a point of $B$ closest to $x$, is coarsely Lipschitz in both $x$ and $B$ (in terms of Hausdorff distance for~$B$), with implicit constants depending only on the hyperbolicity constant of $X$ and the quasiconvexity constants of $B$.
\end{enumeratecontinue}
\qed\end{lemma}

We will apply these properties to curve complexes of surfaces and their subsurfaces, as
follows: if $A$ is a finite subset of $\MM(S)$ then we let $\hull_S(A)$ denote
$\hull_{\CC(S)}(A')$, where $A'$ is the set of curves in the bases of the
markings in~$A$.  Similarly, if $W\subset S$ we let $\hull_W(A)$ denote
$\hull_{\CC(W)}(\pi_W(A))$ where $\pi_W\co\MM(S) \to \CC(W)$ is the usual
subsurface projection.

\subsection{$\Sigma$-hulls and their projections}
If $A$ is a finite subset of $\MM(S)$, and $\ep>0$, we define\index{1aasigmahull@$\Sigma_\ep(A)$, $\Sigma$-hull (coarse version)}
$$
\Sigma_\ep(A) = \{\mu\in\MM(S): d_W(\mu,\hull_W(A)) \le \ep,\ \  \forall W\subseteq
S\}.
$$
Here $W$ varies over all essential subsurfaces of $S$ (including $S$) and $\hull_W(A)$ is the hyperbolic hull as defined in \S\ref{S:hyperbolic hulls}. These sets, one for each $\epsilon \ge 0$, are called the \emph{coarse $\Sigma$-hulls} or just the \emph{$\Sigma$-hulls} of $A$\index{sigmahull@$\Sigma$-hull} in $\MM(S)$. Usually we will assume $\epsilon$ is large enough so that the conclusions of Propositions~\ref{hull retraction} and~\ref{coarse Sigma properties} apply.

It is clear that $A\subset \Sigma_\ep(A)$ and that $\Sigma_\epsilon(A) \subset \Sigma_{\epsilon'}(A)$ if $\epsilon \le \epsilon'$, but a priori not much else. (For the reader familiar with the constructions in Masur-Minsky \cite{masur-minsky:complex2}, we note one of our motivations for this definition: there exists $\ep_0$ such that, 
if $\ep > \ep_0$ then $\Sigma_\ep(A)$ contains every hierarchy path 
between points $a,a'\in A$.)

In order to understand $\Sigma$-hulls better we will need
a family of coarse retractions. 

\begin{proposition}{hull retraction}
There exists $\epsilon_0 \ge 0$ depending only on $\xi(S)$ such that for each $\epsilon \ge \epsilon_0$ the following hold. Given a finite set $A\subset \MM(S)$ there exists a map\index{1aacoarseretraction@$p_A$, coarse retraction to $\Sigma$-hull}
$$ p_A \co \MM \to \Sigma_\ep(A)
$$
which is a coarse retraction. That is, 
\begin{enumerate}
\item
 $p_A|_{\Sigma_\ep(A)}$ is uniformly close to the identity.
\item
$p_A(x)$ is coarse-Lipschitz not just in $x$, but jointly
in $x$ and in $A$ (using the Hausdorff metric on $A$).
\item 
For each $W\subseteq S$, let $y_W$ be a nearest point on
    $\hull_W(A)$ to $\pi_W(x)$. Then 
$$ d_W(p_A(x), y_W) 
$$
is uniformly bounded. 
\end{enumerate}
The implicit constants depend only on $\ep$, $\xi(S)$, and
the cardinality of $A$.
\end{proposition}

\begin{proof}
The proof will be an application of the Consistency Theorem
\ref{consistency suffices}.  Given $ x \in\MM(S)$, for each
$W\subseteq S$ let $y_W=y_W( x ,A)$ be a nearest point to $\pi_W( x )$
on $\hull_W(A)$.

\begin{lemma}{Sigma projection consistent}
For any $x \in \MM(S)$ and finite $A \subset \MM(S)$, the tuple $(y_W( x ,A))_W$ satisfies the consistency conditions (C1-2) of \S\ref{consistency}, with constants
$c_1,c_2$ depending only on the cardinality of $A$.
\end{lemma}

\begin{proof}
To prove (C1), let $U\pitchfork V$.
First, by Lemma~\ref{behrstock inequality} we have
$$
\min\left(d_U( x ,\boundary V),d_V( x ,\boundary U)\right) < m_0
$$
Consider the case where $d_U( x ,\boundary V)<m_0$; the other case is similar. Now 
if $d_U( x ,y_U) < 2m_0+2$, we are done, because $d_U(y_U,\boundary V) < 3m_0+4$. We have used here that $\diam_U(x)$ is bounded above by $2$ --- see comments after Lemma~\ref{coarse lipschitz} --- and we shall use the same bound several times below. 

If $d_U(x,y_U) \ge 2m_0+2$ then 
$d_U( x ,\hull_U(A)) \ge 2m_0+2$ since $y_U$ was the nearest point to
$\pi_U( x )$, and we conclude by the triangle inequality that
\begin{equation}\label{bdV far}
d_U(\boundary V,\hull_U(A)) \ge m_0.
\end{equation}
Now since $y_V\in \hull_V(A)$, there must exist $a,b\in A$ and $a'\in\pi_V(a)$, $b'\in\pi_V(b)$ such that $y_V\in[a',b']$. Now $d_U(\boundary V,a)$ and
$d_U(\boundary V,b)$ are $\ge m_0$ by (\ref{bdV far}), and it follows by Lemma~\ref{behrstock inequality} that $d_V(\boundary U,a)<m_0$ and $d_V(\boundary U,b) < m_0$, and so $d_V(a',b') < 2m_0+6$. Since $y_V\in [a',b']$ we conclude that one of $d_V(y_V,a')$, $d_V(y_V,b')$ is $< m_0+3$, implying that $d_V(y_V,\boundary U) < 2m_0+5$, and again we are done.

That is, we have shown that (C1) holds with $c_1=3m_0+5$.

\medskip

It remains to prove (C2). Let $V\esssubset W$, and suppose that
$d_W(y_W,\boundary V) >4$.
We will bound  $d_V(y_V,y_W)$. 

Suppose first that $d_V(y_W, x ) \ge m_0$. Then by Theorem~\ref{bounded geodesic 
projection}, the $\CC(W)$-geodesic $[y_W,\pi_W( x )]$ must pass
through a point $t$ within 1 of $\boundary V$. By the assumption that  
$d_W(y_W,\boundary V) >4$, it follows that $d_W(t,y_W)>3$ 
and hence  $d_W(t,x)< d_W( x ,y_W)-3$.  
Now let $\gamma$ be a $\CC(W)$-geodesic $[\pi_W(a),\pi_W(b)]$ for
$a,b\in A$. If $\gamma$ were to pass within $1$ of $\boundary V$ then
it would pass within $2$ of $t$, so there would be a point of $\gamma$
which is within $d_W( x ,y_W) - 1$ of $\pi_W( x )$. 
This contradicts the choice of $y_W$ as a closest point to
$\pi_W( x )$. We conclude, by Theorem~\ref{bounded geodesic 
projection} that
$\diam_V(\gamma)< m_0$, and hence $\diam_V(A) < m_0$. 

Moreover, since $y_W$ itself is on such a geodesic,
$d_V(y_W,A)<m_0$. Since $y_V\in\hull_V(A)$ we also have
$d_V(y_V,A)<m_0$ and we conclude $d_V(y_V,y_W) < 3m_0$. 

Now suppose that $d_V(y_W, x ) < m_0$. 
Let $a,b\in A$ be such that $y_W \in [\pi_W(a),\pi_W(b)]$. Now, by our 
assumption that $d_W(y_W,\boundary V) >4$, we have that $\pi_W(\boundary V)$ 
may be within distance 1
of either subsegment $[\pi_W(a),y_W]$ or $[y_W,\pi_W(b)]$, but not
both. Suppose the former. Then by Theorem~\ref{bounded geodesic 
projection} we have $d_V(y_W,b) < m_0$. This yields that
$d_V( x ,A) < 2m_0$, and hence that the closest point $y_V$ to
$\pi_V( x )$ is within $3m_0$ of $\pi_V(y_W)$. 

Hence we have proved (C1) and (C2) both hold 
 with constants $c_1 = 3m_0+5$ and $c_2=4$. 
\end{proof}

We turn now to the proof of Proposition~\ref{hull retraction}. Using Lemma~\ref{Sigma projection consistent}, the definition and properties of $p_A$ follow directly from
Theorem~\ref{consistency suffices}: 
given $x\in\MM(S)$ and $(y_W)$ as in part (3), 
Lemma~\ref{Sigma projection consistent} tells us that $(y_W)$ satisfies
conditions (C1-2) with uniform constants, and hence by 
Theorem~\ref{consistency suffices}  there exists $\mu\in\MM(S)$ with
$d_W(\mu,y_W) < c_3$ for uniform $c_3$ and all $W\subseteq S$. 
We take $\epsilon_0 = c_3$, for any $\epsilon \ge \epsilon_0$ we define $p_A(x) \equiv \mu \in \Sigma_\epsilon(A)$, and clearly (3) holds. 

\medskip

Finally, let us show that the rest of the proposition follows from
(3). To see (1), let $x\in \Sigma_\ep(A)$, i.e.,  for all $W$ we have
$d_W(x,\hull_W(A)) \le \ep$, and so $d_W(x,y_W) \le \ep$.
Now if $\mu = p_A(x)$ we have from (3)
that $d_W(\mu,y_W)$ is uniformly bounded, and hence we have a uniform
bound on $d_W(x,\mu)$. Corollary~\ref{C bounds M} of the quasidistance formula
now gives us a bound on $d_\MM(x,p_A(x))$.

To prove (2), suppose that we have $d_{\MM}(x,x') < b$ and $d_H(A,A')
< b$, where $d_H$ is Hausdorff distance in $\MM$. The coarse-Lipschitz
property of $\pi_{W}$ (Lemma \ref{coarse lipschitz}) implies that for any $W$ we have bounds of
the form $d_{W}(x,x') < b'$, and $d_{H,\CC(W)}(A,A') < b'.$
The latter implies a Hausdorff distance bound 
$$d_{H,\CC(W)}(\hull_W(A),\hull_W(A')) \le b''
$$
by Lemma~\ref{hyperbolic hulls}(3). If $y'_W$ is a nearest point to $\pi_W(x')$ in $\hull_W(A')$ then, by Lemma~\ref{hyperbolic hulls}(4) we obtain a uniform bound on $d_{W}(y_{W},y'_{W})$.

But (3) now implies that $d_W(p_A(x),p_{A'}(x'))$ is uniformly
bounded for all $W$. Again the quasidistance formula gives us a uniform bound on
$d_{\MM}(p_A(x),p_{A'}(x'))$. 
\end{proof}

As a consequence of Proposition \ref{hull retraction} we obtain the
following facts: 

\begin{lemma}{coarse Sigma properties} There exists $\epsilon_0 \ge 0$ depending only on $\xi(S)$ such that for all $\epsilon,\epsilon' \ge \epsilon_0$ and all $I$ there exist $K$, $C$, and $\ep''$ such that if $A ,A' \subset \MM(S)$ each have cardinality $\le I$ then
\begin{enumerate}
\item $\diam(\Sigma_\ep(A)) \le K \, \diam(A)+C$
\item If $A'\subset \Sigma_\ep(A)$ then $\Sigma_\ep(A') \subset
  \Sigma_{\ep''}(A)$.
\item $d_H(\Sigma_\epsilon(A'),\Sigma_\epsilon(A)) \le K d_H(A',A)+C$.
\item $d_H(\Sigma_\epsilon(A), \Sigma_{\epsilon'}(A)) \le C$
\end{enumerate}
\end{lemma}

\begin{proof}
Parts (1) and (2) follow from the definition of $\Sigma_\ep(A)$ and the
quasidistance formula.

To prove (3), note by Proposition \ref{hull retraction}~(2) the maps $p_A$ and $p_{A'}$ differ by at most $K \, d_H(A',A) + C''$, for some $K,C''$ (with the appropriate dependence). Since $p_{A'}$ is uniformly close to the identity on $\Sigma_\ep(A')$, the restriction of $p_A$ to $\Sigma_\ep(A')$ must be within $K \, d_H(A',A) + C'$ of the identity, for some $C'$. It follows that $\Sigma_\ep(A')$ is within  $K \, d_H(A',A) + C$ of $\Sigma_\ep(A)$, for some $C$. The opposite inclusion is obtained in the same way. 

To prove (4), we may assume that $\epsilon \le \epsilon'$ from which it immediately follows that $\Sigma_\epsilon(A) \subset \Sigma_{\epsilon'}(A)$. Consider the two projection maps $p_{A,\epsilon} \from \MM(S) \to \Sigma_\epsilon(A)$, $p_{A,\epsilon'} \from \MM(S) \to \Sigma_{\epsilon'}(A)$. Given $\mu \in \Sigma_{\epsilon'}(A)$, by applying Proposition \ref{hull retraction}~(3) we obtain a uniform bound on $d_W(y_W,p_{A,\epsilon'}(\mu))$ and on $d_W(y_W,p_{A,\epsilon}(\mu))$ and so also on $d_W(p_{A,\epsilon'}(\mu),p_{A,\epsilon}(\mu))$, over all essential subsurfaces $W \subset S$. Corollary~\ref{C bounds M} then gives a bound on $d(p_{A,\epsilon'}(\mu),p_{A,\epsilon}(\mu))$. Also, Proposition~\ref{hull retraction}~(1) gives a bound on $d(\mu,p_{A,\epsilon'}(\mu))$, and since $p_{A,\epsilon}(\mu) \in \Sigma_\epsilon(A)$ we are done.
\end{proof}

We shall also have use for the following lemma, where $m_0$ is the constant in Lemma~\ref{behrstock inequality}.

\begin{lemma}{hull projections} There exists $m_1 \ge 0$ depending only on $\xi(S)$ such that if $A\subset\MM(S)$ is any subset and $W\subset S$ is any essential subsurface satisfying 
$$\diam_{\CC(W)}(A) > m_1,
$$
then for all essential subsurfaces $U\subseteq S$  with
$U\pitchfork\boundary W$ we have
$$d_{\CC( U)}(\partial W,\hull_U(A))\leq m_{0}
$$
\end{lemma}

\begin{proof} Since $U\pitchfork\boundary W$ we either have $W \esssubset U$ and $W \not\isotopic U$, or $W\pitchfork U$. We treat these two cases separately. 
      
First, if $W\subset U$ and $W \not\isotopic U$ then Theorem~\ref{bounded geodesic projection} immediately implies that 
$$d_{\CC( U)}(\boundary W,\hull_U(A))\leq 1
$$
as long as $\diam_{\CC(W)}(A) > B$.
      
Now, consider the case that $W \pitchfork U$. As long as $\diam_{\CC( W)}(A)>2m_{0}+2$, it follows that there exists $a\in A$ for which $d_{\CC(W)} (a,\partial U)\geq m_{0}$. Then by Lemma~\ref{behrstock inequality} it follows that $d_{\CC( U)}(a,\partial W)< m_{0}$.
\end{proof}

\newcommand{\lra}{\to}
\def\Si{\Sigma}
\def\defeq{\mathrel{:=}}
\newcommand{\nerve}{\operatorname{Nerve}}
\newcommand{\im}{\operatorname{Im}}


\section{Contractibility and homology}
\label{contractibility}

In this section we prove contractibility of $\Sigma$-hulls in 
$\AM(S)$ (Lemma \ref{hull contractibility}), and use this to develop
\emph{$\Sigma$-compatible} chains in the cone. This has 
applications in Sections~\ref{separation} and \ref{finiteness2}, as 
well as providing another proof of Hamenst\"adt's theorem on the homological dimension of $\AM(S)$.

If $\seq A$ is a finite set in $\AM$ represented by a sequence
$(A_n)$, the $\Sigma$-hull of~$\seq A$, denoted $\Sigma(\seq A)$\index{1aasigmahull@$\Sigma(\seq A)$, $\Sigma$-hull (asymptotic version)}\index{sigmahull@$\Sigma$-hull}, is the ultralimit of the coarse $\Sigma$-hulls $\Sigma_\ep(A_n)$, where $\ep$ is a fixed constant chosen sufficiently large so that the the lemmas in Section~\ref{hulls} apply.  Note that changing $\ep$ does not change $\Sigma(\seq A)$, by Lemma~\ref{coarse Sigma properties}~(4), nor does changing the representatives, by Lemma~\ref{coarse Sigma properties}~(3).  In fact Lemma~\ref{coarse Sigma properties} applied in the limit gives:

\begin{lemma}{Sigma properties} For all $I \ge 0$ there exists a constant $K$, depending also on $\xi(S)$, such that if $\seq A \subset \AM$ is a set of cardinality $\le I$ then:
\begin{enumerate}
\item $\diam(\Sigma(\seq A)) \le K \, \diam(\seq A)$ 
\item If $\seq A'\subset \Sigma(\seq A)$ then $\Sigma(\seq A') \subset
\Sigma(\seq A)$
\end{enumerate}
\end{lemma}

The retractions $p_{A_n}$ of 
Proposition \ref{hull retraction} ultraconverge to a Lipschitz retraction
$$p_{\subseq A} : \AM \to \Sigma(\seq A)
$$
whose Lipschitz constant depends only on $\xi(S)$ and the cardinality of $\seq A$. Moreover, Proposition \ref{hull retraction} implies that $p_{\subseq A}$ is jointly continuous in its arguments and in the points of
$\seq A$. With this we can establish: 

\begin{lemma}{hull contractibility}
$\Sigma(\seq A)$ is contractible. 
\end{lemma}

\begin{proof}
First we note that $\Sigma(\seq A)$ is path-connected: $\AM(S)$ is
path-connected since it is the asymptotic cone of a path-metric space. Hence
given $\seq a, \seq b\in\Sigma(\seq A)$, let $\gamma(t)$ be a path connecting
them and note that $p_{\subseq A}\circ \gamma$ is a path in $\Sigma(\seq A)$
connecting them. 

Now write $\seq A = \{\seq a_0,\seq a_1,\ldots,\seq a_k\}$, and
for $j=1,\ldots,k$  let $\seq a_j(t)$ be a path in $\Sigma(\seq A)$ from $\seq a_0$ to
$\seq a_j$, where $\seq a_j(0)=\seq a_0$ and $\seq a_j(1) = \seq a_j$. 

Let $\seq A_t = \{\seq a_0, \seq a_1(t),\ldots,\seq a_k(t)\}$ for
$t\in[0,1]$, and let $p_t$ be the retraction from $\AM$ to
$\Sigma(\seq A_t)$. $p_t$ varies continuously in $t$, takes values
within $\Sigma(\seq A)$, and we note that $p_1$ restricted to
$\Sigma(\seq A)$ is the identity while $p_0$ is a constant. Hence
$\Sigma(\seq A)$ is contractible. 
\end{proof}

\subsection{$\Sigma$-compatible chains and homological dimension}

In this subsection we use $\Sigma$-hulls as a device to control singular chains in $\AM(S)$, in terms of what we call \emph{$\Sigma$-compatible} chains. With these we compute the homological dimension from the result of \cite{behrstock-minsky:rank} that the 
topological (covering) dimension of compact subsets of $\AM(S)$ is bounded by $\xi(S)$. We also recall a local homology theorem of Kleiner-Leeb \cite{KleinerLeeb:buildings} (Theorem~\ref{local homology} below) and its corollary~\ref{chain control},  which we will use in 
Sections~\ref{separation} and \ref{finiteness2} to control the support of embedded top-dimensional manifolds in $\AM(S)$ in terms of $\Sigma$-compatible chains.

A \emph{polyhedron}\index{polyhedron} is a finite simplicial complex, and a \emph{polyhedral pair} $(P,Q)$ consists of a polyhedron $P$ and a subcomplex $Q$. A \emph{polyhedral $n$-chain}\index{polyhedral chain} in a space $X$ is a continuous map from an $n$-dimensional polyhedron to $X$, where the domain is equipped with a specified orientation and coefficient on each $n$-simplex. A \emph{polyhedral $n$-cycle}\index{polyhedral cycle} is a polyhedral $n$-chain such that, in the simplicial chain complex of the domain, the linear combination of the $n$-simplices has zero boundary. Polyhedral chains and cycles in $X$ represent singular chains and cycles in the traditional sense, by restricting the map to individual simplices and taking the indicated formal linear combination. Every homology class in $X$ can be represented by a polyhedral cycle.

A continuous map $f\co P\lra \AM$ from a polyhedron to 
$\AM$ is \emph{$\Si$-compatible}\index{Sigma@$\Sigma$-compatible} if for each face $\tau\subset P$,
$$
f(\tau)\subset \Si(f(\tau^{(0)})),
$$ 
where $\tau^{(0)}$ denotes the $0$-skeleton of $\tau$.
By applying Lemma~\ref{coarse Sigma properties}~(1) it follows that if $f\co P\lra \AM$ is $\Si$-compatible then for every
face $\tau\subset P$, 
\begin{equation}
\label{eqn:compatiblediameter}
\diam(f(\tau))\leq \diam(\Si(f(\tau^{(0)})))\leq C\;\diam(f(\tau^{(0)})),
\end{equation}
where the constant $C=C(\dim\,\tau)$ depends explicitly on $\dim\,\tau = \#\text{vertices}(\tau)-1$, and also depends implicitly on 

\begin{lemma}{lem:sicompatibleextension}
Suppose $(P,Q)$ is a finite dimensional polyhedral pair, where the zero skeleton
of $Q$ coincides with the zero skeleton of $P$.
Then any $\Si$-compatible map $f_0\co Q\lra \AM$ can be extended to 
a $\Si$-compatible map $f\co P\lra \AM$.
\end{lemma}

\begin{proof}
The map $f$ may be constructed by induction on the relative $k$-skeleton
using the contractibility of hulls.
\end{proof}

\begin{lemma}{lem:perturb}
If $\ep>0$ and $f_0\co P\to \AM$ is a map from a finite polyhedron to $\AM$, 
then there is a  map $f_1\co P\to\AM$
such that

1. $f_1$ factors through a polyhedron of dimension $\leq \xi(S)$.

2.  $d(f_0,f_1)<\ep$.

\end{lemma}

\noindent (Here $d(f,g) = \sup_{x\in P} d(f(x),g(x))$).

\begin{proof}
Pick $\rho>0$.

Let $Y\defeq f_0(P)\subset \AM$.  Since the topological
dimension of $Y$ is $\leq \xi(S)$ by \cite{behrstock-minsky:rank},
there is an open cover 
$\U=\{U_i\}_{i\in I}$ of $Y$ such that $P'\defeq \nerve(\UU)$
has dimension at most $\xi(S)$, and $\diam(U_i)<\rho$ for all $i\in I$.
Let $\{\phi_i \co Y\to [0,1]\}_{i\in I}$ be a partition of unity subordinate
to $\UU$, and  $\phi\co Y\to P'$ be the map with barycentric
coordinates given by the $\phi_i$'s.

Next, for each $i\in I$, pick $x_i\in U_i$, and using  
Lemma~\ref{lem:sicompatibleextension} construct 
a $\Si$-compatible map $\alpha\co P'\to \AM$
with the property that $\alpha(U_i)=x_i$ (recall that the vertex
set of $\nerve(\UU)$ consists of elements of $\UU$).  

Set $f_1\defeq \alpha\circ\phi\circ f_0\to\AM$.  

We now estimate $d(f_0,f_1)$.

Pick $x\in P$.  If $\phi\circ f_0(x)$ lies
in an open face $\tau\subset P'$ whose vertices
are $U_{i_1},\ldots,U_{i_k}$, then $f_0(x)\in U_{i_1}\cap\ldots\cap U_{i_k}$,
and
$$
f_1(x)\in \Si(\{x_{i_1},\ldots,x_{i_k}\}).
$$
Therefore for a constant $C$ depending only on $\xi(S)$ we have: 
\begin{equation}
\begin{aligned}
d(f_0(x),f_1(x))&\leq d(f_0(x),x_{i_1})+d(x_{i_1},f_1(x))\\
&\leq \rho+C\diam(\{x_{i_1},\ldots,x_{i_k}\})\\
&\leq \rho+2\;C\rho.
\end{aligned}
\end{equation}
So when $\rho< \frac{\ep}{1+2C}$ we will have $d(f_0,f_1)<\ep$.
\end{proof}

\begin{lemma}{lem:controlledhomotopy}
Let $P$ be a finite polyhedron.
Given a pair of maps $f_0,f_1\co P\to \AM$, 
there is a homotopy $\{f_t\}_{t\in [0,1]}$
from $f_0$ to $f_1$ whose tracks have diameter $<C\;d(f_0,f_1)$,
where $C=C(\dim P)$. 
\end{lemma}

\begin{proof}
Pick $\rho>0$. By subdividing $P$ we may assume without loss of generality that for $i\in \{0,1\}$ and every face $\tau$ of $P$, 
$$\diam(f_i(\tau))<\rho.
$$
Let $P=P_1, P_2,\ldots,P_k,\ldots$ be a sequence of successive  barycentric
 subdivisions of $P$, so the mesh size tends to zero.  For $i\in \{0,1\}$, $k\in \Z_+$, 
let $f_{i,k}\co P_k\to \AM$ be a $\Si$-compatible map agreeing with $f_i$
on the $0$-skeleton of $P_k$.  Since $f_i$ is uniformly continuous, the 
diameter estimate (\ref{eqn:compatiblediameter}) implies that $f_{i,k}$
converges uniformly to $f_i$ as $k\to\infty$.

We will construct the homotopy from $f_0$ to $f_1$ as an infinite concatenation
of homotopies 
$$
f_0\ldots  \stackrel{H_{0,3}}{\sim}f_{0,3}\stackrel{H_{0,2}}{\sim}f_{0,2}\stackrel{H_{0,1}}{\sim} f_{0,1}
\stackrel{H}{\sim} f_{1,1}\stackrel{H_{1,1}}{\sim}f_{1,2}\stackrel{H_{1,2}}{\sim}f_{1,3}\ldots  f_1.
$$  

\noindent
The homotopy $f_{i,k}\stackrel{H_{i,k}}{\sim}f_{i,k+1}$ is constructed
as follows.  Triangulate $P\times [0,1]$ such that the $0$-skeleton lies
in $P\times \{0,1\}$, and the induced triangulation of $P\times \{j\}$ agrees with $P_{k+j}$,  for $j\in \{0,1\}$. Now apply Lemma \ref{lem:sicompatibleextension} to get a homotopy from $f_{i,k}$ to $f_{i,k+1}$. The homotopy $H$ is constructed similarly.

Now consider the track of the point $x\in P$ during the homotopy 
$H_{i,k}$. The point $x\in P$ lies in some open simplex 
$\tau_k$ of $P_k$. Let $\sd\tau_k$ be the barycentric subdivision of $\tau_k$, a subcomplex of $P_{k+1}$. By the uniform continuity of $f_i$ and the fact that $f_{i,k}\ra f_i$ uniformly, it follows that the diameter of $f_{i,k}(\tau_k)\cup f_{i,k+1}(\sd\tau_k)$ tends to zero as $k\ra \infty$.
Now for every $t\in[0,1]$, the point $(x,t)\in P\times [0,1]$ lies in
a face  of the subdivision of $P\times[0,1]$ used to construct $H_{i,k}$,
and this face has vertices in $\left(\tau_k\times\{0\}\right) \union \left(\sd\tau\times\{1\}\right)$.
By Lemma \ref{lem:sicompatibleextension}, we get 
$d(H_{i,k}(x,t),f_i(x))< \delta_k$, where $\delta_k<\rho$
and $\delta_k\to 0$
as $k\to\infty$.  It follows that the concatenation of 
$$
H_{i,j},H_{i,j+1},\ldots
$$
has tracks of diameter tending to zero as $j\to \infty$, yielding
a homotopy $f_{i,1}\sim f_i$ whose tracks have diameter $<C_1\rho$.

Similar estimates imply that the tracks of $H$ have diameter 
$<C_2(d(f_0,f_1)+\rho)$.  So if $\rho$ is sufficiently small, we 
obtain the desired homotopy.
\end{proof}

We now give some corollaries of Lemmas \ref{lem:perturb} and
\ref{lem:controlledhomotopy}.

Our first corollary produces acyclic sets in $\AM(S)$. A set $X\subset \AM(S)$ is {\em $\Sigma$-convex}\index{sigmaconvex@$\Sigma$-convex} if, for any finite set $\seq A\subset X$ the hull $\Sigma(\seq A)$ is in $X$ as well.  

\begin{corollary}{AcyclicHulls}
If $X \subset \AM(S)$ is open and $\Sigma$-convex then $X$ is acyclic.
\end{corollary}

Notice that, as a consequence, $\AM(S)$ is itself acyclic. But contractibility of $\AM(S)$ was already known, as a consequence of the fact that $\MCG(S)$ is automatic \cite{Mosher:automatic}, by using a folk theorem which says that the combing lines of an automatic structure induce a contraction in the asymptotic cone.

\begin{proof}[Proof of Corollary \ref{AcyclicHulls}] We claim that any polyhedral cycle $f\colon P\to X$ can be refined and then approximated by a $\Sigma$-compatible cycle, and openness allows us to do this within~$X$. To see how, refine $P$ until the mesh size is sufficiently small, apply Lemma \ref{lem:sicompatibleextension} to the pair $(Q,Q^0)$ where
$Q$ is the refined polyhedron and $Q^0$ is its 0-skeleton, and then
use Lemma \ref{Sigma properties} to show that the new map $f'\colon Q\to\AM(S)$ is sufficiently close to the old to still be in $X$. Moreover Lemma \ref{lem:controlledhomotopy} tells us the new and old maps are homotopic with similar control, so that the homotopy may be made to lie in $X$. The old and new cycles therefore represent the same homology class.

Now let $\CC(Q)$ be the cone of $Q$, the join of $Q$ with a single point $p$. Extend $f'$ to $Q \union p$ by taking $f'(p)$ to be in $f'(Q^0)$. Applying Lemma \ref{lem:sicompatibleextension} extend $f'$ to a $\Sigma$-compatible map on $\CC(Q)$. By 
$\Sigma$-convexity this polyhedral chain lies in $X$, so our original
cycle bounds in $X$, and therefore $X$ is acyclic.
\end{proof}

Our next corollary computes the homological dimension of $\AM(S)$. Note that in Hamenst\"adt's
approach \cite{hamenstadt:III} the homological statement comes
directly.

\begin{corollary}{homology dimension}
If $(U,V)$ is an open pair in $\AM$, then $H_k(U,V)=\{0\}$ for all 
$k>\xi(S)$.
\end{corollary}

\begin{proof}
Pick $[c]\in H_k(U,V)$.  Then there is a finite polyhedral pair
$(P,Q)$ and a continuous map of pairs $f_0\co (P,Q)\to (U,V)$ such that
$$
[c]\in \im\bigl(H_k(P,Q)\stackrel{f_{0*}}{\lra}H_k(U,V)\bigr). 
$$

Pick $\ep>0$. 
Applying Lemma \ref{lem:perturb}, we obtain a continuous map $f_1\co P\to \AM$
with $d(f_0,f_1)<\ep$, such that $f_1$ factors through a polyhedron $P'$ of 
dimension
at most $\xi(S)$.   By Lemma \ref{lem:controlledhomotopy} there is a homotopy
$\{f_t\}_{t\in [0,1]}$ whose tracks have diameter $<C\ep$, where $C=C(\dim P)$.

If $\ep$ is sufficiently small, then $f_1$ will induce  a map of pairs
$(P,Q)\to (U,V)$, and the homotopy $\{f_t\}$ will be a homotopy of maps
of pairs, so that $f_0$ and $f_1$ induce the same map $H_k(P,Q)\to H_k(U,V)$.
But since $f_1$ factors through a polyhedron $P'$ of dimension $\leq 
\xi(S)$, by subdividing $P'$  if necessary we can arrange that
$f_1$ factors as $(P,Q)\to (P',Q')\to (U,V)$, where $(P',Q')$ is a
polyhedral pair of dimension $\leq \xi(S)$ .  This implies that $f_{1*}=0$.
Hence $[c]=0$.
\end{proof}

We are now in a position to apply the following local homology results of Kleiner-Leeb \cite{KleinerLeeb:buildings}, which we will be using in the proof of Theorems \ref{Lambda misses flat} and \ref{local cube finiteness}. 

\begin{theorem}{local homology}
Let $X$ be a contractible metric space and suppose $H_k(U,V)=0$ for any open pair $V\subset U\subset X$ and $k>n$. If $M\subset X$ is an embedded $n$-manifold then 
$$ H_n(M,M-p) \to H_n(X,X-p)
$$
is injective for any $p\in M$. 
\end{theorem}

\begin{corollary}{chain control}
Let $X$ be a contractible metric space and suppose $H_k(U,V)=0$ for any open pair $(U,V)$ and $k>n$. Let $M\subset X$ be an
oriented compact $n$-manifold with boundary, and let $C$ be a singular chain in $X$, such that $\boundary C = \boundary M$. Then $M \subset C$.
\end{corollary}

By Corollary \ref{homology dimension}, we will be able to apply Corollary
\ref{chain control} in the setting of $\xi(S)$-dimensional manifolds
in $\AM(S)$. 

\newcommand{\jetslice}{\LL}
\newcommand{\jetsliceom}{\jetslice_\omega}

\section{Separation properties}
\label{separation}

In this section we develop the notion of {\em jets}, which are local
structures in the cone corresponding to sequences of geodesics in
subsurface complexes. Projections to a jet serve to control separation
properties in the cone. The two main results of the section are
Theorem~\ref{projection separation} and Theorem~\ref{Lambda acyclic}, 
which are concerned with separation properties
of {\em microscopic} and {\em macroscopic} jets, respectively. Much of
the technical work is done in Lemma~\ref{quantitative projection}. 
Section~\ref{pieces} provides a brief digression, where we deduce 
information about the tree-graded structure of $\AM(S)$ as an 
application of microscopic jets.

\subsection{Jets}

Recall the following definition from \MMTwo. Consider a finite type
surface $Y$. If $\xi(Y)\geq 2$, a \emph{tight geodesic}\index{tight
geodesic} in $\CC(Y)$ is a sequence of simplices $\sigma =
(w_0,\ldots,w_n)$ such that any selection of vertices $v_i\in w_i$
yields a geodesic in the 1-skeleton of $\CC(Y)$, and such that for
each $1\le i \le n-1$, the curve system $w_i$ is the boundary of the
subsurface filled by $w_{i-1}$ and $w_{i+1}$. If $\xi(Y)=1$, every
geodesic in $\CC(Y)$ is considered to be tight. If $Y \subset S$ is
an essential annulus, then every geodesic in $\CC(Y)$ is considered
to be tight as long as it satisfies a technical finiteness condition
on the endpoints of arcs representing the vertices. It is shown in
\MMTwo\ that any two vertices in $\CC(Y)$ can be joined by a tight
geodesic, and there are only finitely many possibilities. When the
sequence $(w_0,\ldots,w_n)$ is understood, we use the shorthand
notation $[w_0,w_n]$, and we also refer to
$[w_i,w_j]=(w_i,\ldots,w_j)$ as a \emph{subsegment} of $[w_0,w_n]$.

Let $a,b\in\MM(S)$, let $W\subseteq S$ be a connected essential
subsurface, and let $g$ be a tight geodesic in $\CC(W)$ from an
element of $\pi_{\CC(W)}(a)$ to an element of $\pi_{\CC(W)}(b)$. If
$\sigma=[\alpha,\beta]$ is a subsegment of $g$, we call
$(\sigma,a,b)$ a \emph{tight triple} supported in~$W$. Let $|\sigma|$
denote the length of $\sigma$ in $\CC(W)$. Although we have
suppressed the geodesic $g$ from the notation for a tight triple,
when we need to refer to it we shall call it the \emph{ambient
geodesic}.

We also associate to the triple $(\sigma, a, b)$ the following pair of
points in $\MM(W)$: $\iota(\sigma,a,b) = \alpha\rfloor
\pi_{\MM(W)}(a)$,\index{1aainitialmarking@$\iota(\cdot)$, initial
marking} the \emph{initial marking}\index{initial marking} of the
triple; and $\tau(\sigma,a,b) = \beta\rfloor
\pi_{\MM(W)}(b)$,\index{1aaterminalmarking@$\tau(\cdot)$, terminal
marking} the \emph{terminal marking}.\index{terminal marking} Up to
the usual bounded ambiguity one can think of $\iota$ as $\alpha\union
\pi_{\MM(W_\alpha)}(a)$, where $W_\alpha$ denotes the union of
$W\setminus \alpha$ with the annuli whose cores are the curves of
$\alpha$; and similarly for $\tau$. 

We define
$$\| \sigma \|_{(W,a,b)} = \dist_{\MM(W)}(\iota(\sigma, a,
b),\tau(\sigma,
a, b)).
$$

Using Theorem~\ref{distance formula} we can establish the following
properties of this notion of size:

\begin{lemma}{triple size}
Given $\sigma=[\alpha,\beta], W, a, b$ as above,  letting
$\Phi(\sigma)$
denote the set of subsurfaces $Y\esssubset W$ that do not overlap
some simplex of $\sigma$, we have
\begin{enumerate}
\item $\displaystyle \|\sigma\|_{(W,a,b)} \ceq{}
d_{\CC(W)}(\alpha,\beta) + \sum_{Y\in \Phi(\sigma)}\Tsh
A{d_{\CC(Y)}(a,b)}
$
\item If $\sigma$ is written as a concatenation of subintervals
  $\sigma_1*\ldots *\sigma_k$, then
$$
\|\sigma\|_{(W,a,b)} \ceq{} \sum_i \|\sigma_i\|_{(W,a,b)}
$$
\end{enumerate}
where the constant $A$ and the constants of approximation depend only
on $\xi(W)$.  
\end{lemma}

\begin{proof}
Let $\iota=\iota(\sigma,a,b)$ and
$\tau=\tau(\sigma,a,b)$.
To prove (1), first recall that Theorem~\ref{distance formula}
gives us, for large enough $A$ and uniform constants of approximation
(depending on $A$), that
$$ d_{\MM(W)}(\iota,\tau) \ceq{} \sum_{Y\esssubset W} 
\Tsh A{d_{\CC(Y)}(\iota,\tau)}
$$
On the other hand,  Theorem~\ref{bounded geodesic projection} gives 
a constant $B$ such that, if $Y \esssubset W$, $Y \not\isotopic W$,
and if $Y$ overlaps every simplex of $\sigma$, then
$\diam_{\CC(Y)}(\sigma) \le B$; in particular since $\alpha, \beta$
are in $\sigma$, and since $\alpha,\beta$ are contained in $\iota$
and $\tau$ respectively, we get
$$d_{\CC(Y)}(\iota,\tau) \le B
$$
Thus, if the threshold constant $A$ is raised above $B$ all of these
terms drop out of the sum, leaving the $Y=W$ term, and what is almost
the summation in (1), indexed by $Y \in \Phi(\sigma)$, except 
with $d_{\CC(Y)}(\iota,\tau)$ in place of $d_{\CC(Y)}(a,b)$.  

Now consider $Y \esssubset W$ which does not overlap some simplex
of~$\sigma$. By tightness of the ambient geodesic $g$
containing~$\sigma$, the set of simplices in $g$ not overlapping $Y$
is a contiguous sequence of at most 3 simplices. If $Y$ overlaps
$\beta$ then the rest of $g$ between $\beta$ and $\pi_{\CC(W)}(b)$
consists of simplices overlapping $Y$ and so, since $\tau$ contains
$\beta$, Theorem~\ref{bounded geodesic projection} implies
$$ d_{\CC(Y)}(\tau,b) \le B.
$$
If $Y$ does not overlap $\beta$ then $Y\subset W_\beta$. By
definition the restriction of $\tau$ to $W_\beta$ is
$\pi_{\MM(W_\beta)}(b)$, and again we have a uniform bound on
$d_{\CC(Y)}(\tau,b)$. The same logic yields a uniform bound on
$d_{\CC(Y)}(\iota,a)$. Thus, at the cost of again raising the
threshold, we can replace $\iota,\tau$ in the sum by $a,b$ --- thus
completing the proof of (1).

To prove (2), we simply apply the approximation of (1) to each
$\sigma_i$ separately and sum, noting that any  $Y\subset W$ there are
at most 3 (successive) simplices disjoint from it, and hence it can
be in at most
$4$ different $\Phi(\sigma_i)$. This bounds the overcounting by a
factor of 4, and gives the estimate. 
\end{proof}

\medskip

A {\em jet}, denoted $J$, is a quadruple  $(\bbar\sigma, \bbar W,
\bbar a,
\bbar b)$, where $(\sigma_n,a_n,b_n)$ are tight triples with
$\sigma_n$ supported in $W_n$, and we assume that
$\bbar a$ and $\bbar b$ have ultralimits in $\AM(S)$ (i.e., that they
do not go to $\infty$ faster than linearly).  We refer to $\bbar W$ as
the \emph{support surface} of the jet $J$. The sequence of initial
points $\iota_n=\iota(\sigma_n,a_n,b_n)$ defines a point
$\iotaomega(\bbar\sigma, \bbar a, \bbar b) \in\AM(\bbar W)$, which we
will call the {\em basepoint} of the jet and denote $\seq\iota(J)$
or just $\seq\iota$ when $J$ is understood. Similarly one obtains
$\seq\tau(J)$.  

Call a jet {\em microscopic} if $\|\sigma_n\|_{(W_{n},a_{n},b_{n})}$
grows sublinearly --- that is, if
$$ \frac1{s_n} \|\sigma_n\|_{(W_{n},a_{n},b_{n})} \to_\omega 0.
$$
A jet $J$ is {\em macroscopic} if it is not microscopic, which occurs
if and only if $\seq\iota(J) \ne \seq\tau(J)$. Often we write
$\|\sigma_n\|_{J}$\index{1aaanorm@$\parallel\sigma\parallel_J$} 
to denote $\|\sigma_n\|_{(W_{n},a_{n},b_{n})}$.

\subsection{Projection and separation properties of microscopic jets}
Let $J= (\bbar \sigma, \bbar W, \bbar a, \bbar b)$ be a microscopic
jet with basepoint $\seq\iota\in\AM(\bbar W)$.  As in
Section~\ref{SectionCubes} we have product regions
$$
\QQ(\boundary W_n) \homeo \MM(W_n) \times \MM(W_n^c)
$$
which give rise in the cone to
$$
\QQomega(\boundary\bbar W) \homeo \AM(\bbar W) \times \AM(\bbar W^c).
$$
We let $\jetslice_n(J)$ denote the slice
$\QQ(\iota_n\union\boundary W_n)$, which by Proposition~\ref{Q
product structure} can be identified
with $\{\iota_n\}\times\MM(W_n^c)$. In the cone we get
$$\jetsliceom(J) = 
\QQomega(\bbar \iota \union \boundary\bbar W)\homeo
\{\seq \iota\}\times\AM(\bbar W^c).
$$
Applying Lemma~\ref{DimensionLemma}, the locally compact dimension of
$\jetsliceom(J)$ equals $\delta(\bbar W^c)$ which, by applying the
codimension formula, Proposition~\ref{CodimensionProp}, equals
$\xi(S) - \delta(\bbar W) = \xi(S) - \xi(\bbar W)$.

Denote by $\pi_{\sigma_n} \from \MM(S) \to \sigma_n$ the composition
of projection $\MM(S) \to \CC(W_n)$ with closest point projection
$\CC(W_n) \to \sigma_n$.

\subsubsection*{Projection equivalence}
In terms of the jet $J$ we can define a relation on sequences $(x_n)$
in $\MM(S)$ as follows. Say that $(x_n) \sim_{\bbar\sigma} (x'_n)$ if 
$$ d_{\CC(W_n)}(\pi_{\sigma_n}(x_n),\pi_{\sigma_n}(x'_n))
$$
is bounded for $\omega$-a.e.\ $n$. It is immediate that this is an
equivalence relation on the ultraproduct $\MM(\bbar S)$. We will
deduce the following stronger result.

\begin{theorem}{projection separation}
For any microscopic jet $J$, the relation $\sim_{\bbar\sigma}$
descends to an equivalence relation on $\AM(S) \setminus
\jetsliceom(J)$.  Moreover, every equivalence class is open.
\end{theorem}

This theorem is a consequence of the following more quantitative
statement:

\begin{lemma}{quantitative projection inequivalent}
There exists $C>0$ such that for any microscopic jet $J$, if $(\xi_n)$
and $(\xi'_n)$ are sequences in $\MM(S)$ representing $\seq\xi,
\seq\xi' \in \MM_\omega(S)$, and if $(\xi_n)$, $(\xi'_n)$ are {\em
inequivalent} under $\sim_{\bbar\sigma}$, then
$$ d(\seq\xi,\seq\xi') \ge C \, d(\seq\xi,\jetsliceom(J)).
$$
\end{lemma}

We now show how the lemma implies the theorem.

\begin{proof}[Proof of Theorem~\ref{projection separation}]
If $\seq\xi = \seq\xi'$ then Lemma~\ref{quantitative projection
inequivalent} implies either $(\xi_n)\sim_{\bbar\sigma}(\xi'_n)$ or
$\seq\xi\in\jetsliceom(J)$.  Hence in the complement of
$\jetsliceom(J)$ the equivalence relation $\sim_{\bbar\sigma}$
descends to an equivalence relation in the asymptotic cone.

Further, if $\seq\xi\notin\jetsliceom(J)$ then
Lemma~\ref{quantitative projection inequivalent} implies that there
is a positive radius neighborhood of $\seq\xi$ consisting of points
represented by sequences which are $\sim_{\bbar\sigma}$ equivalent to
$\seq\xi$.  Hence equivalence classes are open. 
\end{proof}

Lemma~\ref{quantitative projection inequivalent} is an immediate
consequence of the following stronger statement, which will have
other applications in what follows:

\begin{lemma}{quantitative projection}
There exist $K,C>0$ such that for any microscopic jet $J$, if
$(\xi_n)$ and $(\xi'_n)$ are sequences in $\MM(S)$ representing
points $\seq\xi, \seq\xi' \in \MM_\omega(S)$, and if for \oaen\  we
have 
$d_{\CC(W_n)}(\pi_{\sigma_n}(\xi_n),\pi_{\sigma_n}(\xi'_n))>K$, then 
$$
d(\seq\xi,\seq\xi') \ge C \, d(\seq\xi,\jetsliceom(J)).
$$
\end{lemma}

\begin{proof}
Proposition~\ref{Q product structure} (2) gives us the
following estimate on distance to $\jetslice_n(J)$:
\begin{align*}
d_{\MM(S)}(\xi_{n},\jetslice_n(J) ) & \ceq{}
\sum_{Y\pitchfork(\iota_n\union\boundary W_n)}
\Tsh A{d_{\CC(Y)}(\xi_n,\iota_n\union\boundary W_n)}\\
& \ceq{} 
\sum_{Y\esssubset W_{n}} \Tsh 
A{d_{\CC(Y)}(\xi_{n},\iota_n)} 
+ \sum_{Y\pitchfork\partial W_{n}}
\Tsh A{d_{\CC(Y)}(\xi_{n},\partial W_{n})}
\end{align*}
where $A$ is any sufficiently large threshold and the approximation
constants depend only on $A$. 

For each $Y$ indexing this sum, we will show an inequality of the form
\begin{equation}\label{termwise Y in Wn}
d_{\CC(Y)}(\xi'_n,\xi_n) \ge d_{\CC(Y)}(\xi_n,\iota_n) -
d_{\CC(Y)}(\iota_n,\tau_n) - q
\end{equation}
if $Y\esssubset W_n$, and of the form
\begin{equation}\label{termwise Y intersect Wn}
d_{\CC(Y)}(\xi'_n,\xi_n) \ge d_{\CC(Y)}(\xi_n,\boundary W_n) - q
\end{equation}
if $Y\pitchfork\boundary W_n$, where $q$ is a 
uniform constant. Since the left hand sides of these inequalities are
terms in the quasidistance formula for $d_{\MM(S)}(\xi_n,\xi'_n)$, we
will obtain (with the usual threshold adjustment)
\begin{equation}\label{bound with iota tau}
d_{\MM(S)}(\xi_n,\xi'_n) \ge p' d_{\MM(S)}(\xi_n,\jetslice_n(J)) -
p''d_{\MM(W_n)}(\iota_n,\tau_n) - q'
\end{equation}
where $p',p'',q''$ are additional constants.
This will be sufficient, since by assumption
$\ulim d_{\MM(S)}(\iota_n,\tau_n)/s_n = \ulim \|\sigma_n\|_J/s_n =
0$, and hence the second
term disappears in the asymptotic cone.  We proceed to establish
(\ref{termwise Y in Wn}) and~(\ref{termwise Y intersect Wn}). 

\medskip

Let 
$$ x_n = \pi_{\CC(W_n)}(\xi_n), \qquad x'_n = \pi_{\CC(W_n)}(\xi'_n),
$$
and
$$ z_n=\pi_{\sigma_n}(\xi_n), \qquad z'_n = \pi_{\sigma_n}(\xi'_n).
$$ 
Let $h_n=[x_n,z_n]$ and $h'_n = [x'_n,z'_n]$ be $\CC(W_n)$-geodesic
segments.
Because $z_n$ is a nearest point to $x_n$ on $\sigma_n$ (and similarly
for $z'_n$ and $x'_n$), and $\CC(W_n)$ is $\delta$-hyperbolic, there
is a constant $K_\delta$ such that, if $d(z_n,z'_n) > K_\delta$, 
the union $T_n=\sigma_n\union h_n\union h'_n$ can
be considered as a finite tree, and the distance function of
$\CC(W_n)$ restricted to $T_n$ is approximated by the distance
function along the tree, up to some additive error $\delta'$. 

In the case that $Y \isotopic W_n$, we immediately find that
\begin{align}\label{Y is Wn}
d_{\CC(W_n)}(x'_n,x_n) &\ge d_{\CC(W_n)}(x_n,z_n) - \delta' \notag \\
   & \ge d_{\CC(W_n)}(x_n,\iota_n) - \diam_{\CC(W_n)}(\sigma_n) -
\delta'.
\end{align}
which is (\ref{termwise Y in Wn}) in this case. 

Consider next the case that $Y \esssubset W_n$ and $Y \not\isotopic
W_n$. Let $B$ be the bound in Theorem~\ref{bounded geodesic
projection}. Suppose first that $\boundary Y$ is disjoint from a
radius~1 neighborhood of $h_n$ in $\CC(W_n)$. Then $z_n \pitchfork
Y$, and $d_{\CC(Y)}(x_n,z_n) \le B$. Moreover, $\boundary Y$ can only
be disjoint from simplices on one side of $z_n$ in $\sigma_n$ (not
both) since $\sigma_n$ is a tight geodesic. It follows that
$\pi_Y(z_n)$  is within $B$ of either $\pi_Y(\iota_n)$ or
$\pi_Y(\tau_n)$, and hence
$$
\min\{d_{\CC(Y)}(x_n,\iota_n),d_{\CC(Y)}(x_n,\tau_n)\} \le 2B.
$$
It follows that
$$
d_{\CC(Y)}(x'_n,x_n) \ge 0 \ge d_{\CC(Y)}(x_n,\iota_n) -
d_{\CC(Y)}(\iota_n,\tau_n) - 2B - 3,
$$
which again gives (\ref{termwise Y in Wn}).

Now suppose that $\boundary Y$ intersects a radius~2 neighborhood of
$h_n$. Assuming $K > \max(K_\delta,2\delta' + 4)$, and using the
remark above about the tree $T_n$, it follows that $\boundary Y$ is
disjoint from a radius~2 neighborhood of $h'_n$. Hence the same
argument as above gives
$$ \min\{d_{\CC(Y)}(x'_n,\iota_n),d_{\CC(Y)}(x'_n,\tau_n)\} \le 2B.
$$
If $d_{\CC(Y)}(x'_n,\iota_n)\le 2B$ then the triangle inequality gives
$$ d_{\CC(Y)}(x_n,x'_n) \ge d_{\CC(Y)}(x_n,\iota_n) - 2B - 3
$$
and
if $d_{\CC(Y)}(x'_n,\tau_n)\le 2B$ then the triangle inequality gives
$$
d_{\CC(Y)}(x_n,x'_n) \ge d_{\CC(Y)}(x_n,\iota_n) -
d_{\CC(Y)}(\iota_n,\tau_n) - 2B - 6.
$$
Either way this again gives us (\ref{termwise Y in Wn}).

\medskip

Now consider the case when $Y\pitchfork\boundary W_n$. We may assume
$d_{\CC(Y)}(\xi_n,\boundary W_n) > A$.  As in
Section \ref{consistency}, for essential subsurfaces $U$ and partial
markings $\gamma$ of $S$, for each $k>0$ define the relation $U
\wprec_k \gamma$ by
$$
U\pitchfork \gamma \ \text{and}\  
d_{\CC(U)}(\xi_n,\gamma) \ge k \, c_{1}+4
$$
where $c_{1}$ is a constant which is chosen as follows. As noted in 
Section~\ref{consistency}, the consistency conditions hold for 
$(\pi_U(\xi_n))$ for any sufficiently large constants
$c_1,c_2$. 
Without loss of
generality, we will assume $c_1>\max\{c_2,m_0,B\}$, where $m_0$ is the
constant given by Lemma~\ref{behrstock inequality} and $B$ is the
constant given by Theorem~\ref{bounded geodesic projection}.

We may assume $A > 4(c_{1}+4)$,  and setting $k= \lfloor
d_{\CC(Y)}(\xi_n,\boundary W_n)/(c_{1}+4)\rfloor \ge 4$, we have 
$$ Y\wprec_k W_n. 
$$
Moreover we have $d_{\CC(W_n)}(\xi_n,\xi'_n) > d_{\CC(W_n)}(z_n,z'_n)
- \delta'
> K-\delta'$, so assuming $K-\delta'>2(c_{1}+4)$ we get
$$ W_n \wprec_2 \xi'_n.
$$
Now Lemma \ref{order properties} implies that 
$$ Y \wprec_{k-1} \xi'_n
$$
so in particular 
\begin{equation}\label{Y intersect Wn}
d_{\CC(Y)}(\xi_n,\xi'_n) \ge (\lfloor d_{\CC(Y)}(\xi_n,\boundary W_n)
/(c_{1}+4)\rfloor-1) (c_{1}+4) \geq 
d_{\CC(Y)}(\xi_n,\boundary W_n) - 2(c_{1}+4). 
\end{equation}
This gives us (\ref{termwise Y intersect Wn}). 
\end{proof}

\subsubsection*{Finding microscopic jets} The next lemma constructs
microscopic jets having properties that we will utilize in
Sections~\ref{pieces} and~\ref{finiteness2}.

\begin{lemma}{building microscopic jets} Let $\aomega, \bomega\in
\AM(\bbar S)$ be represented by $\bbar a, \bbar b \in \MM(\bbar S)$,
let $(W_n)$ be a sequence of connected essential subsurfaces, and
suppose  that $d_{\CC(W_n)}(a_n,b_n) \to_\omega \infty$.  Then there
exists a microscopic jet $J=(\bbar\sigma,\bbar W,\bbar a,\bbar b)$
such that
$$ \bbar a\not\sim_{\bbar\sigma} \bbar b
$$
\end{lemma}

\begin{proof}
Let $\mu_n = \pi_{\MM(W_n)}(a_n)$ and $\nu_n = \pi_{\MM(W_n)}(b_n)$.
Let $\ell_n$ be a tight $\CC(W_n)$ geodesic between $\mu_n$ and
$\nu_n$. 

We use a counting argument to produce a sequence of subsegments
$\sigma_n$ of $\ell_n$ with $\abs{\sigma_n} \to \infinity$ but $\|
\sigma_n \|_{(W_n,a_n,b_n)}$ growing sublinearly. Let $f(n)$ be an
integer valued function going to $+\infty$ more slowly than
$|\ell_n|$, meaning that $f(n)\to+\infty$ but $|\ell_n|/f(n) \to
+\infty$. Divide $\ell_n$ into $f(n)$ subsegments, each of length
between $\left( |\ell_n|/f(n) \right) -1$ and $\left( |\ell_n|/f(n)
\right) + 1$. The sum of the $\|\cdot\|_{(W_n, a_n, b_n)}$-sizes of
these subsegments is approximated by $\|\ell_n\|_{(W_n, a_n, b_n)}$
up to bounded multiple by part (2) of Lemma \ref{triple size}, and
this in turn is bounded by a multiple of $d_{\MM(S)}(\mu_n,\nu_n)$ by
part (1)  of Lemma \ref{triple size}, and hence by a multiple of
$s_n$. Therefore
there must be a fixed $C$ such that there is, for \oaen, a subsegment
$\sigma_n$ with $\|\sigma_n\|_{(W_n, a_n, b_n)} \le C s_n/f(n)$. 

Sublinear growth of $\|\sigma_n\|_{(W_n,a_n,b_n)}$ implies that
$J=(\bbar\sigma,\bbar W,\bbar a,\bbar b)$ is a microscopic jet. 
By construction, $a_n$ and $b_n$ project to opposite ends of
$\sigma_n$, and therefore $\bbar a$, $\bbar b$ are inequivalent under
$\sim_{\bbar\sigma}$.
\end{proof}

\subsection{Linear/sublinear decomposition of macroscopic jets}

For macrosopic jets, the way in which the linear growth happens turns
out to be important. 

Consider a macrosopic jet $J = (\bbar\sigma, \bbar W, \bbar a, \bbar
b)$ with $\sigma_n = [\alpha_n, \beta_n]$, and with initial and
terminal markings $\iota_n, \tau_n \in \MM(W_n)$. We will say that
the jet $J$ has {\em sudden growth}\index{growth!sudden} if there
exist simplices $y_n$ and $z_n$ on $\sigma_n$ such that 
\begin{itemize}
\item $\|[\alpha_n,y_n]\|_{J}$ grows sublinearly,
\item $\|[y_n,z_n]\|_{J}$ grows linearly, and
\item $d_{\CC(W_n)}(y_n,z_n)$ is bounded for $\omega$-a.e.\ $n$.
\end{itemize}
We say that $J$ has {\em gradual growth}\index{growth!gradual} if it
does not have sudden growth. 

Given $\bbar z \in \MM(\bbar S)$ we say that $\bbar z$ \emph{escapes
linearly along $J$} if 
$$\|[\alpha_n,\pi_{\sigma_n}(z_n)]\|_{J}
$$
has linear growth, otherwise $\bbar z$ {\em escapes sublinearly}. 

Although linear and sublinear escape are only defined in the
ultraproduct $\MM(\bbar S)$, the following lemma says that they are
defined in the asymptotic cone $\AM(S)$ as long as the jet in
question has gradual growth.

\begin{lemma}{escape defined in cone}
Let $J$ be a macroscopic jet with the gradual growth property. 
The linear/sublinear escape properties for sequences descend to the
ultralimits when these lie in $\AM(S) \setminus \jetsliceom(J)$. 
In other words, 
we can decompose $\AM(S)  \setminus \jetsliceom(J)$ as a disjoint
union, 
$$
\AM(S)  \setminus\jetsliceom(J)= \Omega_J \union
\Lambda_J
$$ 
so that $\seq z \in \Lambda_J$ implies that any  sequence
$(z_n)$ representing $\seq z$ escapes linearly along $J$, and
$\seq z\in\Omega_J$ implies any $(z_n)$ escapes sublinearly. 
\end{lemma}

We will establish the following.

\begin{theorem}{Lambda acyclic}
Let $J$ be a macroscopic jet with the gradual growth property. 
Then $\Lambda_J$ and $\Omega_J$ are both open. Moreover, 
$\Lambda_J$ is $\Sigma$-convex, and therefore is acyclic. 
\end{theorem}

Note that $\Omega_J$ is not acyclic at all  --- indeed it is not
even connected. It breaks up into uncountably many connected
components as an application of Theorem~\ref{projection separation}.

We will establish both the theorem and the lemma as consequences of 
the following more quantitative fact: 

\begin{lemma}{projection to Lambda}
There exists $C>0$ such that the following holds for any macroscopic
jet $J$ with gradual growth. Suppose that $\bbar \xi, \bbar\xi' \in
\MM(\overline S)$ represent $\seq\xi, \seq\xi' \in \AM(S)$ and that 
\begin{enumerate}
\item $\seq \xi \notin  \jetsliceom(J)$
\item $d_{\AM(S)}(\seq\xi,\seq\xi') < C \, d(\seq \xi,
\jetsliceom(J))$
\end{enumerate}
Then either $\xi_n$ and $\xi'_n$ both escape linearly along $J$, or
both escape sublinearly.
\end{lemma}

\begin{proof}[Proof of Lemma \ref{projection to Lambda}:] 
Write $J=(\bbar\sigma,\bbar W,\bbar a,\bbar b)$, let $C$ be the
constant in Lemma \ref{quantitative projection},  and
suppose, by way of contradiction, that (1) and (2) hold but one of
$\bbar \xi,
\bbar\xi'$ escapes sublinearly and the other escapes linearly. After
renaming the one that escapes linearly $\bbar\zeta$ and the one that
escapes sublinearly $\bbar \eta$, we find that
$\pi_{\sigma_n}(\eta_n)$ must precede $\pi_{\sigma_n}(\zeta_n)$ along
$\sigma_n$ for \oaen,  and that
$\|[\alpha_n,\pi_{\sigma_n}(\eta_n)]\|_J$ grows sublinearly while
$\|[\pi_{\sigma_n}(\eta_n),\pi_{\sigma_n}(\zeta_n)]\|_J$ grows
linearly (this uses the additivity property (2) in Lemma \ref{triple
size}).

Since $\bbar\eta$ escapes sublinearly, the restricted jet 
$$
J' = (\bbar\sigma',\bbar W,\bbar a,\bbar b),
$$
defined by letting $\sigma'_n = [\alpha_n,\pi_{\sigma_n}(\eta_n)]$,
is microscopic. 
Moreover, gradual growth implies, as above, that if we enlarge
$\sigma'_{n}$ in the forward direction by an amount which is bounded
for 
\oaen\ then we still obtain a microscopic jet. Thus, we may produce a
new microscopic jet 
extending $\sigma'_n$ along $\sigma_n$ by any
bounded amount which is larger than the constant, $K$, needed to
apply Lemma~\ref{quantitative 
projection}. In this new jet $J''=(\bbar\sigma'',\bbar W,\bbar
a,\bbar b)$ we find that $\pi_{\sigma''_n}(\eta_n) =
\pi_{\sigma_n}(\eta_n)$, while $\pi_{\sigma''_n}(\zeta_n) $ equals (up
to bounded error) the
forward endpoint of $\sigma''_n$, and hence
$$
d_{\CC(W_{n})}(\pi_{\sigma''_n}(\xi'_n), 
\pi_{\sigma''_n}(\xi_n))>K.
$$
Thus 
Lemma~\ref{quantitative projection} implies that 
$d(\seq\xi,\seq\xi') \ge C \, d(\seq\xi,\jetsliceom(J))$. (Note that
$\seq\iota(J)=\seq\iota(J'')$, so $\jetsliceom(J)=\jetsliceom(J'')$.)
This contradicts our  
hypothesis that $d(\seq\xi,\seq\xi') < cd(\seq\xi, \jetsliceom(J))$.
Hence 
it must hold that $\bbar\xi$ escapes 
linearly if and only if $\bbar\xi'$ does.
\end{proof}

Lemma \ref{escape defined in cone} follows immediately from  
Lemma~\ref{projection to Lambda} by considering the 
case $\seq\xi = \seq\xi'$.

\begin{proof}[Proof of Theorem~\ref{Lambda acyclic}] The openness of
$\Lambda_J$ and $\Omega_J$ is an easy consequence of
Lemma~\ref{projection to Lambda}. It remains to prove that
$\Lambda_J$ is $\Sigma$-convex, for we can then apply
Corollary~\ref{AcyclicHulls} to conclude that $\Lambda_J$ is acyclic. 

Let $\seq A \subset \Lambda_J$ be finite and let $\bbar A$ represent
it. Then each $\bbar a\in\bbar A$ has
projections to $\sigma_n$ which escape linearly. The projection of
$\Sigma_\ep(A_n)$ to $\sigma_n$ is, up to bounded error,  the
projection in $\CC(W_n)$ of $\hull_{W_n}(A_n)$ to $\sigma_n$,
and hyperbolicity implies that this is contained (up to bounded error)
in the hull along $\sigma_n$ of the projections of $A_n$. Hence any
point in the hull has projections that escape linearly, and so is in
$\Lambda_J$. This proves $\Sigma$-convexity. 
\end{proof}

\subsection{Classification of pieces}
\label{pieces}

We now record an application of Theorem~\ref{Lambda acyclic}, 
which classifies the maximal subsets of $\AM$ which can not be 
separated by a point. This result will help to motivate the statement
of Theorem~\ref{simultaneous trim}.

First, let us recall the notions of \emph{pieces} and
\emph{tree-graded} spaces as defined by Dru\c{t}u--Sapir
\cite{drutu-sapir:treegraded}.

A complete geodesic metric space $X$ is called \emph{tree-graded} if
it there exists a collection of proper closed convex subsets, $\PP$,
called \emph{pieces}, which pairwise intersect in at most one point
and such that every non-trivial simple geodesic triangle in $X$ is
contained in one piece. It is an easy observation that if $X$
contains a point 
whose removal disconnects it, then $X$ is tree-graded. Further, in
any 
tree-graded space $X$ there exists a unique \emph{finest} way to
write $X$ as a union of pieces none of which can be separated by a
point \cite{drutu-sapir:treegraded}.

In terms of $\CC(S)$--distance, we now provide a complete criterion
for when two points in $\AM(S)$ can be globally separated by a 
point. In particular, the following result describes the pieces in 
the finest decomposition of $\AM(S)$ as a tree-graded space. We note 
that by results of \cite{Behrstock-Drutu-Mosher:thick} such pieces
can not be realized as asymptotic cones of subgroups of $\MM(S)$.

\begin{theorem}{piece condition} Suppose that $\xi(S)\geq 2$. For any
pair of points $\seq \mu, \seq \nu\in\AM(S)$, the following are
equivalent:
\begin{enumerate}
\item No point of $\AM(S)$ separates $\seq \mu$ from $\seq \nu$.

\item In any neighborhoods of $\seq\mu, \seq\nu$, respectively, there
exist points $\seq \mu',\seq \nu'$ with representative sequences
$(\mu'_{n}),(\nu'_n)$ such that  
$$\ulim d_{\CC(S)}(\mu'_{n},\nu'_{n})<\infty
$$
\end{enumerate}
\end{theorem}

\subsubsection*{Remark.} Theorem~\ref{piece condition} can be thought
of as a first step in working out the results of
Section~\ref{finiteness2}, which provide conditions under which a
finite subset of $\AM(S)$ is separated by a product region.
Nevertheless, Theorem~\ref{piece condition} will not be used in the
rest of the paper.

\begin{proof}
    We begin by showing (2) implies (1). Suppose first that
$\seq\mu$ and $\seq\nu$ have
    representative sequences $(\mu_{n})$ and 
    $(\nu_{n})$ for which $\ulim
    d_{\CC(S)}(\mu_{n},\nu_{n})<\infty$.  Hence there is a fixed $m\ge
    0$ such that $d_{\CC(S)}(\mu_{n},\nu_{n}) = m$ for
    $\omega$-a.e. $n$, and we can let 
    $v_{n,0},\ldots, v_{n,m}$ denote the simplices of a tight geodesic
 in  $\CC(S)$ connecting $v_{n,0}\in\base(\mu_{n})$ to $v_{n,m}\in
 \base(\nu_{n})$. For a fixed $i$ let $\bbar v_i = \uprod{v_{n,i}}$. 
The regions $\QQ(v_{n,i})$ have the structure described in Lemma
\ref{Q
  product structure},  and in particular the cone
$\QQ_\omega(\bbar v_i)$ is nontrivial (not a singleton) and
connected. 

Since $v_{n,i}$ and $v_{n,i+1}$ are disjoint (here we use $\xi(S)\ge
2$), we have 
$v_{n,i}\rfloor v_{n,i+1} = v_{n,i}\union v_{n,i+1}$, so 
by Lemma \ref{juncture} 
the intersection  $\QQ_\omega(\bbar v_i) \intersect \QQ_\omega(\bbar
v_{i+1})$ is equal to $\QQ_\omega(\bbar v_i \union \bbar
v_{i+1})$. 
This again is not a singleton, and it follows that the union
$$\QQ_\omega(\bbar v_0)\union\cdots\union \QQ_\omega(\bbar v_m)$$
cannot be disconnected by a point. Since $\seq \mu \in
\QQ_\omega(\bbar v_0)$ and $\seq\nu\in\QQ_\omega(\bbar v_m)$, this
gives property (1) in this case where $\ulim
d_{\CC(S)}(\mu_{n},\nu_{n})<\infty$.  

Now for general $\seq\mu$ and $\seq\nu$ satisfying (2), 
the above argument implies that $\seq\mu$ and $\seq\nu$ can be
approximated arbitrarily closely by $\seq\mu'$ and $\seq\nu'$ which
cannot be separated by a point. 
Since maximal subsets without 
cutpoints are closed \cite{drutu-sapir:treegraded}, this completes the
proof that (2) implies (1).

We now establish that (1) implies (2), by 
proving the contrapositive. Namely suppose that (2) fails to hold for
$\seq\mu$ and $\seq\nu$, so that there exists $r>0$ such that
whenever $d(\seq\mu,\seq{\mu'})\le r$ and
$d(\seq\nu,\seq{\nu'})\le r$,  we have $d_{\CC(S)}(\mu'_n,\nu'_n)\to
\infty$ for any representative sequences. 
We can assume
$r<d_{\AM}(\seq\mu,\seq{\nu})/2$.

Note that Proposition \ref{hull retraction} implies that 
$\Sigma_\ep(\mu_n,\nu_n)$ is coarsely connected. In particular, 
by projecting a continuous path in $\MM(S)$
from $\mu_n$ to $\nu_n$ into $\Sigma_\ep(\mu_n,\nu_n)$, we can obtain
points
$\mu'_n,\nu'_n\in\Sigma_\ep(\mu_n,\nu_n)$ such that  
$d_{\MM(S)}(\mu_n,\mu'_n)$ and $d_{\MM(S)}(\nu_n,\nu'_n)$ are in the
interval $[\half r s_n,r s_n]$ for all sufficiently large $n$. Fix 
such a pair of sequences $(\mu'_{n})$ and $(\nu'_{n})$. It follows
that $d_{\CC(S)}(\mu'_n,\nu'_n) \to \infty$, and by
Lemma~\ref{building microscopic jets} there exists a microscopic jet
$J=(\bbar \sigma,\bbar S,\bbar\mu',\bbar\nu')$ such that 
$\bbar\mu'\not\sim_{\bbar\sigma}\bbar\nu'$. 

Since $\mu'_n$ and $\nu'_n$ are in $\Sigma_\ep(\mu_n,\nu_n)$, the
segments $\sigma_n$ must be within a bounded distance of any 
$\CC(S)$-geodesic between $\base(\mu_n)$  and $\base(\nu_n)$. 
It follows that $\pi_{\sigma_n}(\mu_n)$ and
$\pi_{\sigma_n}(\nu_n)$ are within bounded distance of 
$\pi_{\sigma_n}(\mu'_n)$ and $\pi_{\sigma_n}(\nu'_n)$, respectively. 
Hence we also have $\bbar\mu\not\sim_{\bbar\sigma}\bbar\nu$. 

Now $\jetsliceom(J) = \{\seq\iota(J)\}$ since $J$ is built on the main
surface $S$. We claim that $\seq\mu\ne\seq\iota(J)$ and
$\seq\nu\neq\seq\iota(J)$. This follows from two facts about
$\Sigma$-hulls:

First, for any $a,b\in\MM(S)$ we claim that, if
$a',b'\in\Sigma_\ep(a,b)$, then  
$$
d(a,\Sigma_\ep(a',b')) \gtrsim d(a,\{a',b'\})
$$
 (with uniform constants).
In the projection to each $\CC(W)$, the $\Sigma$-hulls map to 
coarse intervals, and the corresponding inequality
is simply the fact that if two
intervals are nested then the endpoints of the inner one 
separate its interior from the endpoints of the outer one.
The statement then follows from 
the quasidistance formula, Theorem \ref{distance formula}.

Second, if $v$ is a vertex on a tight $\CC(S)$--geodesic $g$ from a
vertex of $\base(a)$ to a vertex of $\base(b)$, then $v\rfloor a$ is
in $\Sigma_\ep(a,b)$, for a uniform $\ep$. This will follow by
showing, for all $W\subseteq S$, that $\pi_{\CC(W)}(v\rfloor a)$ is
uniformly close to $\hull_W(a,b)$. If $W\not\pitchfork v$, the
projections
of $v\rfloor a$ and $a$ are by definition close. If $W\pitchfork v$
and $W\ne S$, then tightness of $g$ implies that either the
subsegment from $v$ to $a$ or the one from $v$ to $b$ consists of
simplices overlapping $W$, and so Theorem~\ref{bounded geodesic
projection} implies that one of $d_W(v,a)$ or $d_W(v,b)$ is uniformly
bounded.  If $W=S$ then $v$ is already in $\hull_S(a,b)$.

Applying the first fact and the choice of $r,\mu'_n$ and $\nu'_n$, we
see that
$$
d_{\MM(S)}(\mu_n,\Sigma_\ep(\mu'_n,\nu'_n)) \gtrsim \half rs_n.
$$ 
From the second fact and the definition of $\iota_n(J)$,  we see that
$\iota_n(J) \in \Sigma_\ep(\mu'_n,\nu'_n)$, 
for a uniform $\ep$. Hence $\seq\mu\ne\seq\iota(J)$. The same applies
to $\seq\nu$.

We have shown that $\bbar\mu\not\sim_{\bbar\sigma}\bbar\nu$, and that
$\seq\mu$ and $\seq\nu$ are different from $\seq\iota(J)$. Hence by
Theorem~\ref{projection separation}, $\seq\mu$ and $\seq \nu$ are
separated by the point $\seq\iota(J)$. 

\end{proof}

\subsection{Manifolds and jets}
As an application of the linear/sublinear decomposition associated to
a macroscopic jet with the gradual growth property, together with the
homology result Theorem~\ref{local homology} and Corollary~\ref{chain
control}, we can obtain:

\begin{theorem}{Lambda misses flat}
Let $E$ be a $\xi(S)$--dimensional connected manifold in $\AM(S)$ 
and let $J$ be a macroscopic jet with the gradual growth property.
Suppose the supporting subsurface $\bbar W$ of $J$ has $\xi(\bbar W)
> 1$.  Then if
$E\intersect \jetsliceom(J) \ne \emptyset$, we conclude
$$
E\intersect \Lambda_J = \emptyset.
$$
\end{theorem}

\begin{proof}
Let $\seq q\in \jetsliceom(J)\intersect E$. Suppose on the contrary
that $E\intersect\Lambda_J \ne \emptyset$. Now $\jetsliceom(J)$ has
codimension at least~2 since $\xi(\bbar W) > 1$. Applying
Theorem~\ref{NonseparationLemma} the set $\jetsliceom(J)$ cannot
separate $E$ which has dimension $\xi(S)$. However, $\jetsliceom(J)$
does separate $\Lambda_J$ from $\Omega_J$ in $\AM$, by
Theorem~\ref{Lambda acyclic}. We conclude that $E \setminus
\jetsliceom(J)$ is contained in $\Lambda_J$. 

Now let $B$ be a ball in $E$ containing $\seq q$ in its interior.
Since $ \jetsliceom(J)\intersect\boundary B$ is a compact set of
codimension $\ge 1$ in $\boundary B$, for any $\ep>0$ there exists a
triangulation of $\boundary B$ with vertices outside $
\jetsliceom(J)$ and mesh size $\ep$.  Using Lemmas \ref{Sigma
properties}, \ref{lem:sicompatibleextension}, and
\ref{lem:controlledhomotopy}, as in the proof of Theorem \ref{Lambda
  acyclic}, 
$B$ can be
deformed to a $\Sigma$-compatible chain $C$, such that every point
moves at
most $c\ep$ (with $c$ a uniform constant) and the 0-skeleton does not
move at all. Since the $0$-skeleton is contained in $\Lambda_J$,
by Theorem~\ref{Lambda acyclic} all of $C$ is
contained in $\Lambda_J$ as well. 
Let $U$ be the $r$-chain giving the homotopy of $\boundary B$ to $C$,
i.e. $\boundary U = \boundary B - C$ and $U$ is supported in a
$c\epsilon$ neighborhood of $\boundary B$. 

Since $C$ sits within $\Lambda_J$, Theorem~\ref{Lambda acyclic}
also implies that it bounds an $r$-chain $B'$ in $\Lambda_J$.

Corollary~\ref{chain control} now implies, since $B$ is embedded and
$\boundary B = \boundary (B'+U)$, that $B\subset B'+U$. Assuming we
have chosen $\ep$ so that $c\ep < d(\seq q,\boundary B)/2$, we find
that $\seq q$ cannot 
be in $U$. Hence $\seq q \in B' \subset \Lambda_J$. This is a
contradiction. 
\end{proof}

\section{Local finiteness for manifolds}
\label{finiteness2}

Our main goal in this section is Theorem \ref{local cube finiteness},
which says that any top-dimensional submanifold of $\AM(S)$ is locally
contained in a union of finitely many cubes.

This will be a consequence of Theorem \ref{Sigma cubes}, in which we
will consider the $\Sigma$-hull of a finite number of points in a
connected top-dimensional manifold in $\AM(S)$, and show that it is
always contained in a finite complex made of cubes of the appropriate
dimension.  In order to do this we will prove Theorem
\ref{simultaneous trim}, which will show that points in the manifold
can be represented by sequences of markings whose projections to all
but the simplest subsurfaces remain bounded.  This in turn will be
possible because of the separation theorems established in Section
\ref{separation}.

\subsection{Trimming theorem}
\label{SectionTrimming}

In Theorem~\ref{piece condition} we showed that for any $\seq \mu,
\seq\nu \in \AM(S)$, if no point in $\AM(S)$ separates
$\seq\mu,\seq\nu$ then --- after perturbation --- the curve complex
distance $d_{\CC(S)}(\mu_n,\nu_n)$ is $\omega$-a.e.\ bounded.

Theorem~\ref{piece condition} can be regarded as a baby version of
Theorem~\ref{simultaneous trim}.  Given a top-dimensional manifold $E$
in $\AM(S)$, the Alexander duality argument given in
Lemma~\ref{NonseparationLemma} shows that $E$ \emph{cannot} be
separated by any product region of codimension $\ge 2$ in $\AM(S)$;
such product regions are associated to sequences of connected,
essential subsurfaces $\bbar W$ such that $\xi(\bbar W) > 1$.  From
this, Theorem~\ref{simultaneous trim} will conclude that the curve
complex diameters of finite subsets of~$E$ --- after trimming --- are
$\omega$-a.e.\ bounded, and this will be true in the curve complexes
of \emph{all} $\bbar W$ such that $\xi(\bbar W) > 1$.

A manifold in $\AM(S)$ is \emph{top dimensional} if it is of dimension
$\xi(S)$, equal to the locally compact dimension of $\AM(S)$.

\begin{theorem}{simultaneous trim}
  Let $\bbar A$ be a finite set of elements in $\MM(\bbar S)$.
  Suppose $A_\omega$ is contained in a connected top-dimensional
  manifold $E\subset \AM(S)$.  There exist $\epsilon$ and $k_0$, a new
  set, $\bbar A' $, and an onto map $\tau\co\bbar A \to \bbar A'$ with
  the following properties.
\begin{enumerate}
\item $\tau(\bbar a)_\omega = a_\omega$ for each $\bbar a \in \bbar A$
\item $A'_n \subset \Sigma_\ep(A_n)$ for $\omega$-a.e.\ $n$.  \item
For any $\bbar W$ with $\xi(\bbar W) > 1$,
$$\diam_{\CC(W_n)}(A'_n) < k_0
$$
for $\omega$-a.e.\ $n$.
\end{enumerate}
\end{theorem}

\emph{Notation:} because $\bbar A$ is finite we can think of it, up to
ultraproduct equivalence, as a sequence of finite sets $A_n$, and we
can think of the trimming map $\tau$ as a sequence of maps from $A_n$
to $A'_n$.  With slight abuse of notation we will also use $\tau$ to
denote these maps, thus writing for example $\tau(a_n)\in \tau(A_n) =
A'_n$.

\begin{proof}
We will argue by induction on the cardinality of $\bbar A$.  The case
of cardinality 1 is trivial, so let us consider the case that $\bbar
A$ has two points.

The fundamental step of the proof is the following lemma, which
``trims'' $\bbar A$ to reduce its projections to a given subsurface
sequence $\bbar W$:

\begin{lemma}{trim W}
Suppose $\bbar A$ has two elements and $A_\omega$ is contained in a
connected top-dimensional manifold $E\subset \AM(S)$.  Let $\bbar W$
be represented by a sequence $(W_n)$ of connected, essential
subsurfaces with $\xi(\bbar W) > 1$.  There exists a map $\tau\co
\bbar A \to \Sigma_\ep(\bbar A)$ such that
$$
(\tau\bbar x)_\omega =x_\omega
$$
for each $\bbar x\in\bbar A$, and
$$
\diam_{\CC(W_n)} \left( \tau(A_n) \right)
$$
is bounded, for $\omega$-a.e.\ $n$.  The constant $\ep$ depends only
on the topological type of $S$.
\end{lemma}

\begin{proof} We may assume that $\diam_{\CC(W_n)}(A_n) \to_\omega
\infinity$, for if not then $\diam_{\CC(W_n)}(A_n)$ is $\omega$-a.e.\
bounded, and we may simply take $\tau$ to be inclusion of $\bbar A$
into $\Sigma_\epsilon(\bbar A)$.  All we need is that
$\diam_{\CC(W_n)}(A_n) \to_\omega \infinity$ is greater than the
constant $m_1$ of Lemma~\ref{hull projections} for $\omega$-a.e.\ $n$,
for we may then apply that lemma, concluding that for each connected,
essential subsurface $U$ of $S$,
\begin{equation}\label{HullBound}
\text{if $U \pitchfork W_n$ for $\omega$-a.e.\ $n$ then $d_{\CC
(U)}(\boundary W_n,\hull_U(A_n)) \le m_{0}$ for $\omega$-a.e.\ $n$}
\end{equation}
This is fact is used repeatedly below, although only once with
explicit details, in the proof of~(3) below.

We find the trimming map $\tau$ in stages.  First let $\tau_1(\bbar
x)$, for $\bbar x\in\bbar A$, be defined as follows:
\begin{equation}\label{tau1}
\tau_1(\bbar x) = \begin{cases} \bbar x & x_\omega \notin
\QQomega(\boundary \bbar W) \\
  \pi_{\QQ(\boundary \bbar W)}(\bbar x) & x_\omega \in
  \QQomega(\boundary\bbar W).
		  \end{cases}
\end{equation}
The notation $\pi_{\QQ(\mu)}\co\MM(S)\to\QQ(\mu)$ denotes the map
$\nu\mapsto \mu\rfloor\nu$ from Section \ref{SectionProducts}.  In
particular the sequence $\pi_{\QQ(\boundary W_n)}$ gives rise to a map
$\pi_{\QQ(\boundary \bbar W)}\co \MM(\bbar S) \to \QQ(\boundary \bbar
W)$.

We claim that $\tau_1(\bbar x)$ has the properties
\begin{enumerate}
\item $\tau_1(\bbar x)_\omega = x_\omega$, \item Either $\tau_1(\bbar
x)_\omega \not\in \QQomega(\boundary\bbar W)$, or $\tau_1(x_n) \in
\QQ(\boundary W_n)$ for $\omega$-a.e.\ $n$.  \item $\tau_1(\bbar x)
\in \Sigma_\ep(\bbar A)$ for suitable $\ep$.
\end{enumerate}
Property (1) is immediate from the definition and the fact that (by
Proposition \ref{Q product structure} (3)) $\pi_{\QQ(\boundary W)}(x)$
realizes, within bounded factor, the distance from $x$ to
$\QQ(\boundary W)$.  Property (2) similarly follows from the
definition.

To see Property (3), note it is obvious when $\tau_1(\bbar x) = \bbar
x$.  Hence, assume that $\tau_1(\bbar x) = \pi_{\QQ(\boundary \bbar
W)}(\bbar x)\in \QQ(\boundary \bbar W)$.  By definition of
$\Sigma_\ep$, it suffices to bound the distance from $\pi_{\CC
(U)}(\tau_1(x_n))$ to $\hull_U(A_n)$, for every $U\subseteq S$ and
$\omega$-a.e.\ $n$.  We treat separately the cases $U \pitchfork \bdy
W_n$ and $U \not\pitchfork \bdy W_n$.

If $U \not\pitchfork \boundary W_n$ for $\omega$-a.e.\ $n$, then
$d_{\CC (U)}(\tau_1(x_n),x_n)$ is uniformly bounded: this can be seen
easily from the definition of the projection $\pi_{\QQ(\boundary
W_n)}$ and the coarse composition properties of subsurface projection
maps, Lemma \ref{coarse composition}.  Since $x_n\in
A_n\subset\Sigma_{\epsilon}(A_{n})$, this gives the desired bound
for~(3).

If $U \pitchfork \boundary W_n$ for $\omega$-a.e.\ $n$, then $d_{\CC
(U)}(\tau_1(x_n),\boundary W_n)$ is bounded by Lemma~\ref{coarse
lipschitz}, and by applying~(\ref{HullBound}) we obtain~(3).

Now for notational simplicity let us assume that we have replaced
$\bbar A$ by $\tau_1(\bbar A)$.  For each $\bbar x \in \bbar A$
properties (1) and (3) become trivial, and property (2) holds with the
consequence that either $x_\omega \not\in \QQ_\omega(\bdy \bbar W)$,
or $\bbar x \in \QQ(\bdy \bbar W)$.  Recalling that $\bbar A$ has two
elements, write $\bbar A = \{\bbar a,\bbar b\}$, and the discussion
separates into two cases: neither of $a_\omega$, $b_\omega$ is in
$\QQ_\omega(\bdy \bbar W)$; or one of $\bbar a, \bbar b$ is in
$\QQ(\bdy \bbar W)$.

\subsubsection*{Case 1:} Neither $\aomega$ nor $\bomega$ lies in
$\QQomega(\boundary \bbar W)$.

We claim that already $d_{\CC (W_n)}(a_n,b_n)$ is bounded for
$\omega$-a.e.\ $n$, and hence there is nothing left to do in this
case.  Suppose otherwise, that $d_{\CC (W_n)}(a_n,b_n) \to_\omega
\infty$.  Then Lemma~\ref{building microscopic jets} yields a
microscopic jet $J = (\bbar\sigma, \bbar W, \bbar a, \bbar b)$, built
from $\CC(W_n)$-geodesics $\sigma_{n}$ for which
$(a_{n})\not\sim_{\bbar\sigma} (b_{n})$.  Moreover, by our assumption
that $\aomega,\bomega \not\in\QQomega(\boundary \bbar W)$, we know
that neither of them is contained in $\jetsliceom(J) \subset
\QQ_\omega(\bdy \bbar W)$.  It follows from Lemma~\ref{projection
separation} that $\jetsliceom(J)$ separates $\aomega$ from $\bomega$.

However, $\jetsliceom(J)$ is homeomorphic to $\AM(\bbar W^c)$, which
by Proposition~\ref{CodimensionProp} has codimension at least 2 since
$\xi(\bbar W) > 1$, so by Lemma~\ref{NonseparationLemma} it cannot
separate $E$.  This contradiction implies that in fact
$d_{\CC(W_n)}(a_n,b_n)$ is bounded $\omega$-a.s.

\subsubsection*{Case 2:} At least one of $\bbar a, \bbar b$, say
$\bbar a$, lies in $\QQ(\boundary\bbar W)$.

Now we consider the projections of $\aomega$ and $\bomega$ to the
factor $\AM(\bbar W)$ of $\QQomega(\boundary\bbar W) \approx
\MM_\omega(\bbar W) \cross \MM_\omega (\bbar W^c)$.  Using the coarse
product structure $\QQ(\bdy W_n) \approx \MM(\bbar W_n) \cross
\MM(\bbar W_n^c)$, let
\begin{align*}
a_n &\approx \left( \pi_{\MM(\bbar W_n)}(a_n), \pi_{\MM(\bbar
W_n^c)}(a_n) \right) \\
       & = (\alpha_n,\beta_n) \\
\pi_{\QQ(\bdy W_n)}(b_n) &\approx \left( \pi_{\MM(\bbar W_n)}(b_n),
\pi_{\MM(\bbar W_n^c)}(b_n) \right) \\
       &= (\gamma_n, \delta_n)
\end{align*}

\subsubsection*{Case 2a:} $\alpha_\omega = \gamma_\omega$.

In this case, we can simply adjust $\bbar a$ so that the projections,
before rescaling, are a bounded distance apart: replace $a_n$ by
$$\tau_2(a_n) =(\gamma_n,\beta_n).
$$
We need to check, as before, that $\tau_2(a_n) \in \Sigma_\ep(A_n)$
for a fixed $\ep$.  This is again done by considering projections to
all connected, essential subsurfaces $U \subset S$, treating
separately the cases $U \esssubset W_n$, $U \esssubset W_n^c$, or $U
\pitchfork \bdy W_n$, for $\omega$-a.e.\ $n$; in the latter
case~(\ref{HullBound}) is again applied, and remaining details are
left to the reader.  Similarly $\tau_2(\bbar a)_\omega = \aomega$, and
of course $d_{\CC(W_n)}(\tau_2(a_n),b_n)$ is now bounded $\omega$-a.s.

Note that this argument works whether or not $\bbar b\in \QQ(\boundary
\bbar W)$.  If it is, then the roles of $\bbar a$ and $\bbar b$ can be
reversed.

\subsubsection*{Case 2b:}
$\alpha_\omega \ne \gamma_\omega$.

In this case, consider the jet $J = (\bbar \sigma,\bbar W,\bbar
a,\bbar b)$, where $\sigma_n = [x_n,y_n]$ with
$x_n\in\pi_{\CC(W_n)}(a_n)$ and $y_n\in\pi_{\CC(W_n)}(b_n)$.  Notice
that the initial marking of $J$ is $\bbar\iota(J) = \bbar \alpha$,
because $\iota_n(J) = \iota(\sigma_n,a_n,b_n) = x_n\rfloor
\pi_{\MM(W_n)}(a_n) = x_n \rfloor \alpha_n = \alpha_n$ where the last
equation follows since $x_n$ is a vertex of $\alpha_n$.  Similarly the
terminal marking of $J$ is $\bbar\tau(J) = \bbar \gamma$.  Since
$\alpha_\omega \ne \gamma_\omega$ it follows that the jet $J$ is
macroscopic.

In the arguments to follow we abbreviate the notation for the jet norm
$\| \sigma_n \|_J = \| \sigma_n \|_{(W_n,a_n,b_n)}$ to $\| \cdot \|$,
and similarly for jets built on subsegments of $\sigma_n$, because in
all cases the jet norm is equal to $\| \cdot \|_{(W_n,a_n,b_n)}$ which
is independent of how the subsegments were chosen.

Since $\aomega$ and $\bomega$ are contained in a connected
top-dimensional manifold $E$, we can apply Theorem~\ref{Lambda misses
flat} to conclude that $J$ cannot have the gradual growth property,
for $\jetsliceom(J) = \QQomega(\bbar\alpha\union \boundary\bbar W)$
contains $\aomega$, and if $J$ had gradual growth then $\bomega$ would
be in $\Lambda_J$, but then Theorem~\ref{Lambda misses flat} would
forbid $\bomega$ from being in~$E$.

Since $J$ has sudden growth, the following must occur: in $\sigma_n$
there must be points $p_n$ and $q_n$ such that $\|[x_n,p_n]\|$ grows
sublinearly and $\|[p_n,q_n]\|$ grows linearly, while
$d_{\CC(W_n)}(p_n,q_n)$ stays bounded.

Let $\tau_3(a_n)$ be the marking obtained by projecting $a_n$ to
$\QQ(p_n\union \boundary W_n)$.  Since $a_n$ is already in $\QQ(\bdy
W_n)$ it follows that $\tau_3(a_n)$ is coarsely equal to the
projection of $a_n$ to $\QQ(p_n)$.  We have
\begin{align*}
d_{\MM(S)}(a_n,\tau_3(a_n)) &\approx d_{\MM(S)}(a_n, \QQ(p_n)) \\
   &\approx d_{\MM(W_n)}(\pi_{\MM(W_n)}(a_n),\QQ_{W_n}(p_n)) \\
   \intertext{the last equation following from Proposition~\ref{Q
   product structure}, where $\QQ_{W_n}(p_n)$ denotes the product
   region in $\MM(W_n)$ associated to the partial marking $p_n$.  This
   product region contains the projection to $\MM(W_n)$ of $p_n
   \rfloor \gamma_n$ and so}
d_{\MM(W_n)}(\pi_{\MM(W_n)}(a_n),\QQ_{W_n}(p_n)) &\prec
d_{\MM(W_n)}(\pi_{\MM(W_n)}(a_n),p_n \rfloor \gamma_n) \\
   &\approx d_{\MM(W_n)}(x_n \rfloor \alpha_n,p_n \rfloor \gamma_n) \\
   & \approx \| [x_n,p_n] \|
\end{align*}
The quantity $d_{\MM(S)}(a_n,\tau_3(a_n))$ therefore grows
sublinearly, and so $\tau_3(\bbar a)_\omega = \aomega$.  As before we
can show that $\tau_3(\bbar a) \in \Sigma_\ep(\bbar A)$.

In the case that $\bbar b\in \QQ(\boundary \bbar W)$, by reversing the
direction of $J$ so that its initial marking is $\bbar \gamma$, the
same argument as above shows that the reversed jet has sudden growth,
and we obtain a path sequence $[v_n,u_n]$ in $\sigma_n$ such that $\|
[v_n,u_n] \|$ grows linearly and $\| [u_n,y_n] \|$ grows sublinearly
while $d_{\CC(W_n)}(v_n,u_n)$ stays bounded $\omega$-a.e. As with
$\bbar a$, in this case we define $\tau_3(b_n) =
\pi_{\QQ(u_n\union\boundary W_n)}(b_n)$, and so $\tau_3(\bbar
b)_\omega = b_\omega$.  In the case that $\bomega \notin
\QQomega(\boundary\bbar W)$, then we simply let $\tau_3(\bbar b) =
\bbar b$, and let $u_n=v_n=y_n$.  In $\sigma_n$ we have a sequence
$x_n \cdots p_n \cdots q_n \cdots v_n \cdots u_n \cdots y_n$, where we
can ensure that $q_n$ precedes $v_n$ by shortening the segments
$[p_n,q_n]$ and $[v_n,u_n]$ if necessary, maintaining the property
that $\| [p_n,q_n] \|$ and $\| [v_n, u_n] \|$ each grow linearly.

We claim now that $d_{\CC(W_n)}(\tau_3(a_n),\tau_3(b_n)) \approx
d_{\CC(W_n)}(p_n,u_n) $ is bounded $\omega$-a.s. For if it were not
then neither would $d_{\CC(W_n)}(q_n,v_n)$ be bounded $\omega$-a.s.,
and so as in case~(1) we could extract a microscopic jet $J' = (\bbar
\sigma',\bbar W,\bbar a, \bbar b)$ with $\sigma'_n$ a subsegment of
$[q_n,v_n]$.  Neither of the points $\aomega$ and $\bomega$ can be in
$\jetsliceom(J')$.  For $\aomega$, this follows from the fact that
$\|[p_n,q_n]\|$ grows linearly and hence insulates $a_n$ from
$\sigma'_n$ --- that is, by Lemma \ref{triple size} (1) and the
quasidistance formula we obtain, term-by-term, a linearly growing
lower bound for $d_{\MM(W_n)}(\alpha_n,\iota_n(J'))$.  For $\bomega$
this is the same argument if $\bbar b\in \QQ(\boundary\bbar W)$, and
if not it is even easier for $\bomega$ is not even in
$\QQomega(\boundary \bbar W)$.

Hence, Lemma \ref{projection separation} would imply that
$\jetsliceom(J')$ separates $E$, and applying
Lemma~\ref{NonseparationLemma} this would again contradict the
assumption that $\xi(\bbar W) > 1$.

We conclude that, in case (2), we can find $\tau_3(\bbar A)$ such that
$\diam_{\CC(W_n)}(\tau_3(\bbar A))$ is bounded $\omega$-a.s. This
concludes the proof of Lemma~\ref{trim W}, where $\tau$ is the
composition of the appropriate $\tau_i$.
\end{proof}

\subsubsection*{Hierarchies of geodesics}
Before we can continue the proof of Theorem \ref{simultaneous trim} we
must recall a few of the details of the construction of hierarchies of
tight geodesics from \cite{masur-minsky:complex2}.  A hierarchy
$H=H(a,b)$ is associated to any $a,b\in\MM(S)$, and is a certain
collection of tight geodesics in curve complexes of connected,
essential subsurfaces of $S$.  The subsurface whose complex contains a
geodesic $h$ is called its {\em domain} $D(h)$.  The properties
relevant to us are the following:

\begin{theorem}{hierarchy properties}
Let $a,b\in \MM(S)$ and $H(a,b)$ a hierarchy between them.
\begin{enumerate}
\item \label{H:main} There is a unique {\em main geodesic} $g_H$ with
$D(g_H)=S$, whose endpoints lie on $\base(a)$ and $\base(b)$.  \item
\label{H:component domain} For any geodesic $h\in H$ other than $g_H$,
there exists another geodesic $k\in H$ such that, for some simplex $v$
in $k$, $D(h)$ is either a component of $D(k)\setminus v$, or an
annulus whose core is a component of $v$.  We say that $D(h)$ is a
{\em component domain} of $k$.  \item \label{H:unique} A subsurface in
$S$ can occur as the domain of at most one geodesic in~$H$.  \item
\label{H:endpoints} For each $h\in H$, the endpoints of $h$ are within
uniformly bounded distance of $\pi_{D(h)}(a)$ and $\pi_{D(h)}(b)$.
\item \label{H:large} For each connected, essential subsurface $W
\subset S$, if $d_{\CC(W)}(a,b) > m_0$, then there exists $h\in
H(a,b)$ with $D(h)=W$.
\end{enumerate}
\end{theorem}

Define the complexity of a geodesic $g \in H(a,b)$ to be $\xi(D(g))$.
The following counting argument allows us, under the appropriate
circumstances, to bound the cardinality of the set of geodesics with a
given lower bound on complexity.

\begin{lemma}{inductive domain count} For all $t \ge 1$, all $a,b \in
\MM(S)$, and all $k \ge 1$, if $d_{\CC(W)}(a,b) \le k$ for all
subsurfaces $W$ with $\xi(W) \ge t$, then the hierarchy $H(a,b)$
contains at most $O(k^{\xi(S)-t})$ geodesics of complexity $\xi=t-1$.
\end{lemma}

\begin{proof}
The proof is by induction, using the properties listed in Theorem
\ref{hierarchy properties}.  Every subsurface of complexity $\xi=s$ in
$H(a,b)$ appears as a component domain in some geodesic of complexity
$> s$.  Hence the number of $\xi=s$ geodesics is bounded by the number
of $\xi>s$ geodesics times the length bound on those geodesics.
\end{proof}

\medskip

We now turn to the proof of Theorem \ref{simultaneous trim} in the
case that $\bbar A$ has two elements.  We first apply Lemma \ref{trim
W} with $\bbar W = \bbar S$.  Thus we obtain $\tau(\bbar A)$, such
that $\diam_{\CC(S)}(\tau(A_n))$ is $\omega$-a.s.\ bounded.  Again for
notational convenience we replace replace $\bbar A$ by $\tau(\bbar A)$
and continue.

Writing $\bbar A = \{\bbar a,\bbar b\}$ as before, we consider
hierarchies $H_n=H(a_n,b_n)$.  By property \pref{H:large}, for any
$\bbar W$ with $\diam_{\CC(W_n)}(A_n)\to_\omega\infty$, $W_n$ must be
a domain in $H_n$ for $\omega$-a.e.\ $n$.  The main geodesics
$g_{H_n}$ have bounded length for $\omega$-a.e.\ $n$, by property
\pref{H:main} and the bound on $\diam_{\CC(S)}(A_n)$.  Applying
Lemma~\ref{inductive domain count} with $t = \xi(S)$ we obtain an
$\omega$-a.s.\ bound on the number of geodesics of complexity
$\xi(S)-1$ in $H_n$, because the domain of such a geodesic must be a
component domain of the main geodesic $g_{H_n}$.  We have therefore
bounded how many $\bbar W$ exist with $\xi(\bbar W) = \xi(S)-1$ and
with $W_n$ a domain in $H_n$ for $\omega$-a.e.\ $n$ --- we use here
the general fact that the ultraproduct of a sequence of sets $X_n$ of
finite cardinality $\le k$ has cardinality $\le k$.  For each such
$\bbar W$ successively, use Lemma~\ref{trim W} again to find $\tau
(\bbar A)$ such that $\diam_{\CC(W_n)}(\tau(\bbar A)_n)$ is bounded,
and again replace $\bbar A$ by $\tau(\bbar A)$ and continue.

Every time we apply Lemma~\ref{trim W}, we maintain the boundedness
that we had for $\diam_{\CC(\bbar U)}$ for any previous $\bbar U$.
This is because $\tau(\bbar A)$ always lies in $\Sigma_\ep(\bbar A)$,
so in the projections to $\CC(U_n)$, it follows that
$\pi_{\CC(U_n)}(\tau(A_n))$ lies uniformly near the hull of
$\pi_{\CC(U_n)}(A_n)$ which is bounded.  Hence after finitely many
steps we have diameter bounds for all $\bbar W$ with $\xi(\bbar W) =
\xi(S)-1$.

This procedure repeats $\xi(S)$ times.  At the $k^{\text{th}}$ step we
have bounds on the lengths of all geodesics of complexity $\ge
\xi(S)-k+1$ that occur in the hierarchy, and by applying
Lemma~\ref{inductive domain count} we bound the {\em number} of
geodesics of complexity $\ge \xi(S)-k$.  A finite number of
applications of Lemma \ref{trim W} renders bounded the projections to
those surfaces without spoiling the previous ones.

The procedure ends when all projections to surfaces of $\xi>1$ are
bounded.  The final set, which we might denote $\tau^N(\bbar A)$ (for
some $N$ which grows with $\xi(S)$ and the bounds at each level), lies
in $\Sigma_{\ep'}(\bbar A)$ (where $\ep'$ depends on $\ep$ and $N$),
and each $\tau^N(\bbar x)$ defines the same point in the cone as
$\bbar x$.

\medskip

This concludes the proof of Theorem \ref{simultaneous trim} when
$\bbar A$ has two elements.  We are now ready for the inductive step,
where we write $\bbar A $ as $\{\bbar a\} \union \bbar B$, and we
assume that there is already a bound on $\diam_{\CC(W_n)}(B_n)$ for
$\omega$-a.e.\ $n$, whenever $\xi(\bbar W) > 1$.

We wish to prove an analogue of Lemma \ref{trim W}, and there is a
similar breakup into cases.  Let $\bbar W$ be such that $\xi(\bbar
W)>1$ and $\diam_{\CC(W_n)}(A_n)\to_\omega \infty$.  First we note as
in the proof of Lemma \ref{trim W} that we may assume (after a first
trimming operation $\tau_1$) that each element $\bbar x\in\bbar A$
either satisfies
$$
x_\omega \notin \QQomega(\boundary \bbar W)
$$
or satisfies
$$
x_n \in \QQ(\boundary W_n)
$$
for $\omega$-a.e.\ $n$ (or as we wrote above, $\bbar x\in
\QQ(\boundary \bbar W)$).

\subsubsection*{Case 1a$'$:} Suppose that there is at least one
element $\bbar b\in\bbar B$ with $\bomega \notin \QQomega(\boundary
\bbar W)$, and that $a_\omega \not\in \QQ_\omega(\bdy \bbar W)$.  Then
the same argument as Case 1 of Lemma~\ref{trim W} shows that
$\diam_{\CC(W_n)}(a_n,b_n)$ is bounded.  Since $\diam_{\CC(W_n)}(B_n)$
was already bounded, this gives us the desired bound for $A_n$.

\subsubsection*{Case 1b$'$:}
Suppose that there is at least one element $\bbar b\in\bbar B$ with
$\bomega \notin \QQomega(\boundary \bbar W)$, but that $\bbar a\in
\QQ(\boundary \bbar W)$.  Depending on whether or not
$\pi_{\MM_\omega(\bbar W)}(a_\omega) = \pi_{\MM_\omega(\bbar
W)}(b_\omega)$, we can apply the argument of Cases 2a and 2b of
Lemma~\ref{trim W}, concluding that $a_n$ can be replaced by
$\tau_2(a_n)$, for which $d_{\CC(W_n)}(\tau_2(a_n),b_n)$ is bounded.
Again since $\diam_{\CC(W_n)}(B_n)$ is assumed bounded we are done.

\subsubsection*{Case 2a$'$:}
Suppose that $\bbar b\in \QQ(\boundary\bbar W)$ for each $\bbar b\in
\bbar B$, and suppose that also $\pi_{\AM(\bbar W)}(B_\omega)$ is a
single point.

In this case, choose one element $\bbar b_0 \in\bbar B$.  Now apply
the argument of Case 2a and 2b in Lemma~\ref{trim W} to $\bbar b_0$
and $\bbar a$.  Note that here $\bbar b_0$ plays the role that $\bbar
a$ played in 2a and 2b, whereas $\bbar a$ itself may or may not be in
$\QQ(\boundary \bbar W)$.  This step produces $\tau_3(\bbar b_0)$
which possibly modifies the $\MM(\bbar W)$ component of $\bbar b_0$
(and similarly for $\bbar a$), so that afterwards their $\CC(W_n)$
distance is $\omega$-a.s.\ bounded.  Define $\tau_3$ on the remaining
elements of $\bbar B$ by making their $\MM(\bbar W)$ components equal
to that of $\tau_3(\bbar b_0)$.  This is a sublinear change which as
before produces points in $\Sigma_\ep(B_n)$.  We now have the desired
bound on $\diam_{\CC(W_n)}(\tau_3(A_n))$.

\subsubsection*{Case 2b$'$:}
Again suppose that $\bbar b\in \QQ(\boundary\bbar W)$ for each $\bbar
b\in\bbar B$, but now suppose that $\pi_{\AM(\bbar W)}(B_\omega)$
contains at least 2 distinct points.  Let $\bbar b_1,\bbar b_2 \in
\bbar B$ have distinct projections to $\AM(\bbar W)$.

If $\bbar a\in \QQ(\boundary \bbar W)$, and $\pi_{\AM(\bbar
W)}(a_\omega) = \pi_{\AM(\bbar W)}({b_i}_\omega)$ for $i=1$ or $i=2$,
then as in Case 2a of Lemma \ref{trim W}, we can replace the
$\MM(\bbar W)$ component of $\bbar a$ to agree with that of $\bbar
b_1$ or $\bbar b_2$, respectively, and are done.

If $\bbar a\in \QQ(\boundary \bbar W)$ but $\pi_{\AM(\bbar
W)}(a_\omega)$ is different from both $ \pi_{\AM(\bbar
W)}({b_1}_\omega)$ and $ \pi_{\AM(\bbar W)}({b_2}_\omega)$, or if
$a_\omega \notin \QQ_\omega(\boundary \bbar W)$, then we work with
$\bbar b_1$ and $\bbar a$ as follows.

If $a_\omega \notin \QQ_\omega(\boundary \bbar W)$ then let
$\tau_4(\bbar a) = \bbar a$.  If $\bbar a\in \QQ(\boundary \bbar W)$,
we argue as in Case 2b of Lemma \ref{trim W}, first to show that a jet
from $\bbar a$ to $\bbar b_1$ cannot have gradual growth, and then to
modify $\bbar a$: along the geodesic from $x_n\in\pi_{\CC(W_n)}(a_n)$
to $y_n\in\pi_{\CC(W_n)}({b_1}_n)$, we find $p_n$ and $q_n$ such that
$\|[x_n,p_n]\|$ grows sublinearly, $\|[p_n,q_n]\|$ grows linearly, and
$d_{\CC(W_n)}(p_n,q_n)$ is $\omega$-a.s.\ bounded.  We then let
$\tau_4(\bbar a) = (p_n\union\boundary W_n)\rfloor a_n$.

Unlike Case 2b of Lemma~\ref{trim W}, we do not attempt to modify
$\bbar b_1$.  Now if $d_{\CC(W_n)}(\tau_4(a_n),{b_1}_n)$ is still
unbounded, we find a microscopic jet $J'$ built from subgeodesics
$\sigma_{n}$ of $[q_n,y_n]$, so that $a_\omega\notin\jetsliceom(J')$
by the same argument at Case 2b.  The points $a_n$ and ${b_1}_n$
project to opposite sides of $\sigma_n$ so $\bbar a \not\sim_{\bbar
\sigma}\bbar b_1$.  Hence if ${b_1}_\omega\notin\jetsliceom(J')$, then
we are done, because $\jetsliceom(J')$ then separates $a_\omega$ from
${b_1}_\omega$ and hence separates $E$, which is a contradiction.  But
if ${b_1}_\omega\in\jetsliceom(J')$ then we must have
${b_2}_\omega\notin\jetsliceom(J')$, because $(b_1)_\omega$ and
$(b_2)_\omega$ have distinct images in $\AM(\bbar W)$.  Since
$d_{\CC(W_n)}({b_1}_n,{b_2}_n)$ is $\omega$-a.s.\ bounded, we also
have $\bbar a \not\sim_{\bbar \sigma} \bbar b_2$, and hence
$\jetsliceom(J')$ separates ${b_2}_\omega$ from $a_\omega$, and we
still have a contradiction.

We conclude that $d_{\CC(W_n)}(\tau_4(a_n),{b_1}_n)$ is $\omega$-a.s.\
bounded, which is what we wanted to show.

This gives the analogue of Lemma~\ref{trim W} for $\bbar A = \bbar
B\union\{\bbar a\}$.  Now we finish the proof as we did before: we
repeatedly apply this result, bounding first the lengths of the main
geodesics in hierarchies between elements of $\bbar A$, and then
inducting downward to bound the lengths of geodesics of lower
complexities, until only domains of complexity~1 are left with
unbounded diameters.
\end{proof}

\subsection{Finitely many cubes}
As a consequence of Theorem \ref{simultaneous trim}, we will show that
the $\Sigma$-hull of a finite number of points in a connected
top-dimensional manifold is composed of finitely many cubes (in the
sense of Section~\ref{SectionProducts}).  From this we'll get the
statement on finitely many orthants in a neighborhood of a point.

\begin{theorem}{Sigma cubes}
If $\seq A$ is a finite subset of a connected top-dimensional manifold
$E$ in $\AM(S)$, then $\Sigma(\seq A)$ is contained in a finite union
of cubes.
\end{theorem}

The first step towards establishing this theorem is the following,
here we are using $\NN_r$ to denote a radius $r$ neighborhood.

\begin{lemma}{G cover}
For each integer $N$, there exists a constant $k_2$, such that for any
finite set $A\subset \MM(S)$ with $\# A=N$ and
$\diam_{\CC(S)}(A_n)<k_{0}$ the following holds for each $a \in A$:
the set $\Sigma_\ep(A)$ is contained in the union of sets
$$
G'(A,\UU,a) = \NN_{k_2}(G(A,\UU,a))
$$
where $\UU$ varies over all sets of the form $\UU(\mu,a)$ for
$\mu\in\Sigma_\ep(A)$.
\end{lemma}

\begin{proof}
Fix $\ep$ large enough for Proposition \ref{hull retraction} (on
retractions of $\Sigma$-hulls) to apply.  For later use fix another
constant $k_0 > \max\{3(m_0+4),2m_0+\ep,B-\epsilon\}$, where $m_0$ is
the constant of Lemma~\ref{behrstock inequality} and $B$ is the
constant of Theorem~\ref{bounded geodesic projection}.

Fix $a \in A$, and now consider any $\mu\in\Sigma_\ep(A)$.  Following
Section~\ref{consistency}, we use $\mu$ and $a$ to define a partial
order among certain subsurfaces of $S$.  
As noted in Section~\ref{consistency}, the projections
$(\pi_W(\mu))$ satisfy the consistency conditions 
with any $c_1 \ge m_0$.  Now define
$\wprec_k$ and $\prec_k$ as in section~\ref{consistency}, that is,
$V\wprec_k W$ if and only if $V\pitchfork \boundary W$ and
\begin{equation}\label{mu prec inequality}
d_{\CC(V)}(\mu,\boundary W) > k(c_1+4),
\end{equation}
whereas $V\prec_k W$ if and only if $V\wprec_k W$ and $V\pitchfork W$.
We choose $c_1$ so that $k_0+\ep=3(c_1+4)$.  In particular
$\FF_3(\mu,a)=\{W:W\wprec_3 a\}$ is the set
$$
\{ W\subsetneq S: d_{\CC(W)}(\mu,a) > k_0+\ep\}.
$$
Lemma \ref{F order} now tells us that $\prec_2$ is a partial order on
$\FF_3(\mu,a)$.  Moreover, by Lemma \ref{order properties}, if $V,W
\in \FF_3(\mu,a)$ and if $V\pitchfork W$ then $V,W$ are
$\prec_2$-ordered.

The set $\FF_3(\mu,a)$ is finite --- using the quasidistance formula
for example, or Lemma \ref{F finite} --- so we can let
$\VV=\VV(\mu,a)$ be the set of $\prec_2$-minimal elements.  Any two
elements of $\VV$ are disjoint or nested in $S$, so let
$\UU=\UU(\mu,a)$ be the subset of $\VV$ consisting of elements maximal
with respect to containment in $S$.  Hence $\UU$ enumerates the
components of an essential subsurface of $S$, which we abuse notation
by also calling $\UU$.  Recall that $\QQ(\boundary\UU)$ has a natural
product structure $\MM(\UU^c) \times\MM(\UU)$.  We claim that $\mu$ is
within uniformly bounded distance of a subset of $\QQ(\boundary \UU)$
of the form
\begin{equation}\label{almost tree product}
G(A,\UU,a) = \{\pi_{\MM(\UU^c)}(a)\} \times \prod_{U\in\abs{\UU}}
\Sigma_{\ep'}(\pi_{\MM(U)}(A))
\end{equation}
where $\abs{\UU}$ is the set of components of $\UU$, and
$\Sigma_{\ep'}$ is defined within $\MM(U)$ just as it was in $\MM(S)$.
The constant $\ep'$ depends only on $\ep$ and $\xi(S)$.

To prove this, we first bound $d_{\MM(S)}(\mu,\QQ(\boundary\UU))$.  By
Proposition \ref{Q product structure}, we just need to establish a
bound on $d_{\CC(W)}(\mu,\boundary \UU)$ for all $W$ that overlap
$\boundary \UU$.  By hypothesis, $\diam_{\CC(S)}(A_n)<k_{0}$ which
implies by Theorem~\ref{bounded geodesic projection} a uniform bound
on $d_{\CC(S)}(\mu,\boundary \UU)$; hence we now assume $W\subsetneq
S$.  Suppose that $d_{\CC(W)}(\mu,\boundary \UU) > 4(c_1+4)$.  In
particular $W\wprec_4 U$ for some $U\in\UU$ such that $W\pitchfork
\boundary U$.  Since $U\wprec_3 a$, by Lemma \ref{order properties}
(2) we have $W\wprec_3 a$, so that $W\in\FF_3(\mu,a)$.

If $W\pitchfork U$ then $W\prec_3 U$ and in particular $W\prec_2 U$,
contradicting the minimality of $U$.  Hence $W$ must contain $U$.
However, by choice of $\UU$ this means $W$ cannot be
$\prec_2$-minimal, so there exists $Z\in\FF_3(\mu,a)$ such that
$Z\prec_2 W$.  By Lemma \ref{order properties} (1), $Z\prec_2
W\wprec_3 U$ implies that $Z\prec_1 U$.  But in particular this means
$Z\pitchfork U$ so they are $\prec_2$-ordered.  $U\prec_2 Z$ would
contradict $Z\prec_1 U$, so we must have $Z\prec_2 U$, but this
contradicts again the minimality of $U$.

We conclude that for all $W$ such that $W\pitchfork \boundary \UU$,
$d_{\CC(W)}(\mu,\boundary\UU) \le 4(c_1+4)$, and this gives a bound of
the form
$$
d_{\MM(S)}(\mu,\QQ(\boundary\UU)) \le k_1
$$
for some $k_1$ depending on $c_1$ (and hence on $m_0$ and $k_0$).

Next we claim that $\pi_{\MM(\UU^c)}(\mu)$ is uniformly close to
$\pi_{\MM(\UU^c)}(a)$.  For this, by the quasidistance formula we need
to bound $d_{\CC(W)}(\mu,a)$ for all $W\subset\UU^c$.  Suppose that
$d_{\CC(W)}(\mu,a) > 3(c_1+4)$, and so $W\in\FF_3(\mu,a)$.  Since $W$
is disjoint from all components of $\UU$ and hence of $\VV$, it is not
$\prec_2$-ordered with or isotopic to any of them.  $W$ cannot be
$\prec_2$-minimal as then it would have to be one of~$\VV$.  Hence
there is some $W'\prec_2 W$ which is $\prec_2$-minimal --- but then
$W'$ is in $\VV$, and again we have a contradiction.

Finally we consider $\pi_{\MM(\UU)}(\mu)$.  Since
$\mu\in\Sigma_\ep(A)$, for each connected subsurface $W \esssubset
\UU$ we have $\pi_W(\mu)\in\NN_\ep(\hull_W(A))$, where
$\NN_\epsilon(\cdot)$\index{1aaneighborhood@$\NN_\epsilon$,
neighborhood of radius $\epsilon$} denotes the radius $\epsilon$
neighborhood (this involves an abuse of distance notation, as
explained under the heading ``Subsurface projections'' in
Section~\ref{SectionBasicDefs}).  But $\hull_W(A)$ is within uniformly
bounded distance of $\hull_W(\pi_{\MM(\UU)}(A))$ by the coarse
composition property of projections (Lemma \ref{coarse composition}).
Hence $\pi_{\MM(\UU)}(\mu)\in\Sigma_{\ep'}(\pi_{\MM(\UU)}(A))$ for
some $\ep'$ depending on $\ep$ and $\xi(S)$.

This establishes that $\mu$ is within uniform distance of the set
$G(A,\UU,a)$ described in (\ref{almost tree product}), proving the
Lemma.
\end{proof}

\begin{proof}[Proof of Theorem~\ref{Sigma cubes}]
From Theorem \ref{simultaneous trim} we may assume that $\seq A$ is
represented by $(A_n)$ such that, for $\omega$-a.e.\ $n$,
$\diam_{\CC(W)}(A_n)$ is bounded by some fixed $k_0$ whenever
$\xi(W)>1$.  Let us consider an arbitrary $A\subset \MM(S)$, of fixed
cardinality $\# A = \#\seq A$, satisfying this condition.

Let $k_{2}$ denote the constant given by Lemma~\ref{G cover} for $\#
A$.  Now we would like to bound the number of $\UU$ that can occur in
the output of Lemma~\ref{G cover}.

If $W$ occurs as a component of $\UU(\mu,a)$ for some $\mu$, then
$d_{\CC(W)}(\mu,a) > k_0 +\ep$.  Since
$\pi_W(\mu)\in\NN_\ep(\hull_W(A))$,
$$
\diam_{\CC(W)}(A) > k_0.
$$
By our assumptions about $A$, this means $\xi(W)\le 1$.  Now define
\begin{align*}
\SSS_1 &= \{ U\subsetneq S: \xi(U) =1 \text{ and } \diam_{U}(A) >
k_0\},\\
\SSS_0 &= \{ U\subsetneq S: \xi(U) =0 \text{ and } \diam_{U}(A) > k_0
\}.
\end{align*}
By Theorem \ref{hierarchy properties} (\ref{H:large}), every element
in $\SSS_0\union \SSS_1$ must be the domain of some geodesic in
$H(a,b)$ for some $a,b\in A$.  Hence the counting argument,
Lemma~\ref{inductive domain count}, directly gives a bound on the
cardinality of $\SSS_1$.

There is no uniform bound for the cardinality of $\SSS_0$, but we can
control the number of annuli $U \in \SSS_0$ which are components of
$\UU(\mu,a)$ for some $\mu \in \Sigma_\epsilon(A)$.  The main reason
for this is the following statement:

\begin{itemize}
\item[$(*)$] For each connected, essential $W \subset S$, each
essential annulus $U \esssubset W$, each $a \in A$, and each $\mu \in
\Sigma_\epsilon(A)$, if $\diam_{\CC(W)}(A) > k_0$, and if $U$ is a
component of $\UU(\mu,a)$, then $d_{\CC(W)}(\boundary U,A) \le k_3$
for $k_3$ depending on $k_0$ and $\# A$.
\end{itemize}

To prove $(*)$, define $\FF_3(\mu,a)$ using $\mu$ as before, and
consider two cases depending on whether $W \in \FF_3(\mu,a)$.

If $W\notin\FF_3(\mu,a)$ then $d_{\CC(W)}(\mu,a) \le k_0+\ep$.  Since
$d_{\CC(U)}(\mu,a) > k_0+\ep>B$, any $\CC(W)$-geodesic from
$\pi_W(\mu)$ to $\pi_W(a)$ must pass within distance 1 in $\CC(W)$ of
$\boundary U$, by Theorem \ref{bounded geodesic projection}.  It
follows that $d_{\CC(W)}(\boundary U,a) \le d_{\CC(W)}(\mu,a) + 1 \le
k_0 + \ep + 1$.

If $W\in\FF_3(\mu,a)$ then, since $U\esssubset W$, the surface $W$
cannot be a $\prec_2$-minimal element of $\FF_3(\mu,a)$, because then
$W$ would have been included in $\UU(\mu,a)$ instead of~$U$.  Hence
there is some element $Y \in \VV(\mu,a)$ such that $Y\prec_2 W$.

We claim that $d_{\CC(W)}(\boundary Y, b)$ is bounded for some $b\in
A$.  The argument is similar to the partial-order arguments in Section
\ref{consistency}.  Since $\mu\in\Sigma_\ep(A)$, we have $\pi_Y(\mu)
\in\NN_\ep(\hull_Y(A))$.  Also, since $d_{\CC(Y)}(\mu,a) > k_0 + \ep$
there must be $b\in A$ such that $d_{\CC(Y)}(a,b) \ge k_0$.  Now
$Y\prec_2 W$ implies that $d_{\CC(Y)}(\mu,\boundary W) \ge 2(c_1+4)>
m_0$ so that $d_{\CC(W)}(\mu,\boundary Y) < m_0$ by Lemma
\ref{behrstock inequality}.  Further, since $d_{\CC(W)}(a,\mu) \ge
k_0$, we have $d_{\CC(W)}(\boundary Y,a) > k_0 - m_0 - 2>m_0$.  Again
by Lemma \ref{behrstock inequality}, we have $d_{\CC(Y)}(\boundary
W,a) < m_0$.  Now since $d_{\CC(Y)}(a,b) \ge k_0$ we have
$d_{\CC(Y)}(\boundary W,b) > k_0 - m_0 -2 > m_0$, so applying
Lemma~\ref{behrstock inequality} one more time we get
$d_{\CC(W)}(\boundary Y,b)<m_0$.

Since $Y$ and $U$ are disjoint, we conclude $d_{\CC(W)}(\boundary U,b)
< m_0 + 1$.  This finishes the proof of $(*)$.

Now we can control the number of elements in $\SSS_0$ which occur as
components of $\UU(\mu,a)$ for $\mu\in\Sigma_\ep(A)$.  Given such a
$U$, there exists $b\in A$ such that $d_{\CC(U)}(a,b) \ge (k_0 -
\epsilon) / 2 > m_0$, and so by Theorem~\ref{hierarchy properties}
there exists a geodesic $h \in H(a,b)$ with domain $D(h)=W$ such that
$\xi(W)\ge 1$ and $U$ is a component domain of $h$.  Noting that
either $h$ has length $\le k_0$ or $\pi_{\CC(W)}(A)$ has diameter $>
k_0$, by applying $(*)$ it follows that there exists $c \in A$ such
that
$$d_{\CC(W)}(\bdy U,c) \le k_4=\max\{k_0,k_3\}
$$
This restricts $\bdy U$, for each $c$, to a segment of length at most
$2k_4$ in $h$.  Now since the number of hierarchies involved is
controlled in terms of $\#A$, and the number of $\xi\ge 1$ surfaces
appearing is controlled in terms of $\# A$ and $k_0$ by
Lemma~\ref{inductive domain count}, this gives us a bound on the total
number of components of the $\UU(\mu,a)$ as $\mu$ varies over
$\Sigma_\ep(A)$.

\medskip

We now apply this result to the sets $A_n$ in the sequence $\bbar A$.
Each one is covered by the uniformly bounded number of sets
$G'(A,\UU,a)$.  Taking rescalings, we obtain in the asymptotic cone
the statement that $\Sigma(\seq A)$ is contained in a finite union of
asymptotic cones of sequences $G(A_n,\UU_n,a_n)$, which by
(\ref{almost tree product}) must be sets of the form
$$
\{x_\omega\} \times \prod_{\bbar U\in|\bbar \UU|} T_{\bbar U}
$$
in $\QQ_\omega(\boundary\bbar\UU)$, where
$x_\omega\in\AM(\bbar\VV^c)$, each component $\bbar U$ of $\bbar \UU$
has $\xi\le 1$, and each $T_{\bbar U}$ is the convex hull of a finite
set in the $\R$-tree $\AM(\bbar U)$.  Hence each $T_{\bbar U}$ is a
finite tree, so after breaking each tree into a finite union of
segments, we obtain the desired finite union of cubes.
\end{proof}

\subsection{Local finiteness}
The main application of Theorem~\ref{Sigma cubes} is the following:
\begin{theorem}{local cube finiteness}
If $E\subset\AM(S)$ is a connected top-dimensional manifold, then any
compact subset of $E$ is contained in a finite union of cubes.
\end{theorem}

\begin{proof}
It suffices to show that a ball $B\subset E$ is contained in finitely
many cubes.

Let $B\subset \text{int}(B')$ where $B'$ is a larger ball.
Triangulate $\boundary B'$ with simplices of diameter smaller than
$r$, where $r$ will be chosen shortly.  Let $f_0\colon B'\to E$ be the
identity, and let $f_1\colon\boundary B'\to \AM(S)$ be a
$\Sigma$-compatible map with respect to the triangulation, which
agrees with $f_0$ on the 0-skeleton; the existence of $f_1$ follows
from Lemma~\ref{lem:sicompatibleextension}.  By Lemma~\ref{Sigma
properties} we have $d(f_0,f_1) < Cr$, and by
Lemma~\ref{lem:controlledhomotopy} there exists a homotopy
$h\colon\boundary B'\times[0,1]\to\AM(S)$ with track diameters at most
$C'r$, for uniform constants $C,C'$.

Choose $r$ small enough that $C' r < \half d(B,\boundary B')$.  Then
we find that the image of $h$ is disjoint from $B$.

Extend the triangulation of $\boundary B'$ to one of $B'$ without
adding any vertices.  Then using Lemma~\ref{lem:sicompatibleextension}
again, $f_1$ can be extended to a $\Sigma$-compatible map $F\colon
B'\to \AM(S)$ with respect to this triangulation.  Let $K$ be the
chain which is the sum of $F$ and $h$ -- then we note that $\boundary
K = \boundary B'$.  By Corollary \ref{chain control}, we conclude that
$B'\subset K$.  Since $B$ is disjoint from $h$, we have
$$
B\subset F.
$$
Now $F$ is contained in the $\Sigma$-hulls of a finite collection of
finite subsets of $E$.  By Theorem \ref{Sigma cubes}, it must
therefore be contained in a finite union of cubes.
\end{proof}

\section{Germs and orthants}
\label{orthant defs}

In this section, we study the local structure of the set of
top-dimensional manifolds passing through a point $\seq x\in\AM(S)$,
by considering the {\em germs} of such manifolds, and using the Local
Finiteness Theorem \ref{local cube finiteness} to relate this to the
{\em complex of orthants} through $\seq x$.

We then apply this to a study of \emph{Dehn twist flats} in the
asymptotic cone $\AM(S)$: given a sequence $\bbar\nu$ of pants
decompositions of $S$, we obtain a sequence $Q(\bbar\nu)$ of Dehn
twist flats in $\MM(S)$, and passing to the rescaled ultralimit we
obtain $Q_\omega(\bbar\nu)$ which, if nonempty, is by definition a
Dehn twist flat in $\AM(S)$.  By Proposition~\ref{Q product
structure}, $Q_\omega(\bbar\nu)$ is a bilipschitz embedded copy of
$\reals^{\xi(S)}$, although we shall not make use of this property.

The main result is Corollary \ref{top char of twist germs}, which
states that germs of Dehn twist flats passing through $\seq x$ admit a
purely topological characterization.  This will be applied in Section
\ref{endgame} in the proof that Dehn twist flats in $\AM(S)$ are
preserved by homeomorphisms, and Dehn twist flats in $\MM(S)$ are
coarsely preserved by quasi-isometries.

\subsection{Poset of Germs}

Consider  the set of closed subsets of $\AM$ containing ${\seq x}$, modulo the equivalence $C \sim C'$ if there exists an open neighborhood $U$ of ${\seq x}$ such that $C\intersect U = C'\intersect U$. The equivalence classes are called  {\em germs through ${\seq x}$},\index{germ} and we let $\gamma(C)$ denote the germ of $C$ through $\seq x$. Note that finite intersection and union yield well-defined operations on the set of germs, and the subset relation is well-defined as well.  Let $\GG$ denote the poset (partially ordered set) of germs at ${\seq x}$; this is a \emph{lattice}, meaning that least upper bounds and greatest lower bounds exist for all pairs $C,C'\in\GG$, namely $C\union C'$ and $C\intersect C'$.

A property $P$ of germs at points of $\AM(S)$ is \emph{topologically characterizable} if any local homeomorphism of $\AM(S)$ taking $\seq x$ to $\seq y$ takes germs at $\seq x$ satisfying $P$ to germs at $\seq y$ satisfying $P$. For example, germs of manifolds are topologically characterizable, as is the dimension function on such germs.

\subsection{Structure of orthants}

In this section we fix $\seq x \in \AM(S)$ and study the set of germs of cubes in $\AM(S)$ for which $\seq x$ is a {\em corner}. These germs will be called \emph{orthants at $\seq x$}. The goal in this section is Lemma~\ref{orthant complex} which shows that the poset of nontrivial orthants at $\seq x$ has the structure of a simplicial flag complex. For this purpose we need to study the relations of equality, subset and intersection of orthants. The main complication that arises is that an orthant has many different representations by sequences of cubes in $\MM(S)$, so these relations are not trivial to detect.

Recall from Section~\ref{SectionCubejunctures} that if $V \esssubset S$ satisfies $\xi(V) \le 1$ then $\TM(V)$ denotes a particular tree quasi-isometric to the marking complex $\MM(V)$, which in the annulus case is a line.

A \emph{cube with distinguished corner} is a cube $C=C(\mu,W,r)$ for which each geodesic $r_i \subset \TM(W_i)$ has a distinguished endpoint $r_i(0)$. The corner of $C$ is, by definition, the marking $\kappa(C) = \{\mu\}\times\prod r_i(0)$, where the right side is interpreted as usual  within $\QQ(\boundary W) \homeo \MM(W^c) \times \prod\MM(W_i)$. Given a sequence $C(\bbar\mu,\bbar W,\bbar r)$ of cubes distinguished corners, we obtain in $\AM(S)$ a cube $C^\omega=C^\omega(\bbar\mu,\bbar W,\bbar r)$ with corner $\kappa^\omega$. 

We recall a few features of the notation for cubes. First, up to ultraproduct equivalence, the sequence $\bbar W = (W^n)$ of essential subsurfaces can be identified with a finite set of sequences $\bbar W_1 = (W^n_1),\ldots,\bbar W_k = (W^n_k)$ of connected essential subsurfaces such that (for $\omega$-a.e.\ $n$) the components of $W^n$ are $W^n_1,\ldots,W^n_k$. Second, the sequence $\bbar r$ can be identified with a finite set of sequences $\bbar r_1 = (r^n_1),\ldots, \bbar r_k=(r^n_k)$ such that $r^n_i$ is a geodesic segment or ray in the tree $\TM(W^n_k)$ with initial point $r^n_i(0)$. Third, recall from Section \ref{SectionCubesCone} that the dimension of the asymptotic cube $C^\omega(\bbar\mu,\bbar W,\bbar r)$ is equal to the number of components $\bbar r_i$ of $\bbar r$ such that the length of the limiting segment $\seq r_i \in \AM(\bbar W)$ is positive --- equivalently the lengths $l(r^n_i)$ grow linearly. Finally, given two essential subsurface sequences $\bbar W, \bbar V$, up to reindexing we may assume that for each $i,j$, if $W^n_i \isotopic V^n_j$ for $\omega$-a.e.\ $n$ then $i=j$; reindexing in this manner is implicit, for example, in the statement of Lemma~\ref{orthant complex}.

Define an \emph{orthant at $\seq x$}\index{orthant} to be the germ $O=\gamma(C^\omega)$ of an asymptotic cube $C^\omega$ with distinguished corner $\kappa^\omega = \seq x$. A \emph{$k$-orthant} is the germ of an asymptotic cube of dimension~$k$. If an orthant~$O$ can be expressed as $O = \gamma(C^\omega(\bbar\mu, \bbar W, \bbar r))$ where all components of $\bbar W$ are annuli, then we say that $O$ is a \emph{Dehn twist orthant}. 

Note that for any orthant $O = \gamma(C^\omega(\bbar\mu,\bbar W, \bbar r))$ at $\seq x$, the asymptotic partial marking $\seq\mu \in \AM(\bbar W^c)$ is determined by the subsurface sequence $\bbar W$ and the corner $\seq x$, being the projection of $\seq x$ to $\AM(\bbar W^c)$. 

For example, the germ of every Dehn twist flat is a union of $2^{\xi(S)}$ Dehn twist $\xi(S)$-orthants. To be precise, consider a Dehn twist flat $Q_\omega(\bbar\nu)$ through $\seq x = x_\omega \in \AM(S)$, where $\base(x_i) = \nu_i$. Let $\bbar W$ be the sequence of annulus neighborhoods of the pants decomposition sequence $\bbar \nu$. The surface $\bbar W^c$ is empty, so we take $\bbar\mu$ to be the empty partial marking. The projection of $\bbar x$ to the line $\TM(\bbar W_i)$ has two directions denoted $\bbar r^1_i, \bbar r^2_i$ (using the usual identification between an ultraproduct of finite sets and a finite subset of an ultraproduct). We may therefore express the germ at $\seq x$ of $Q_\omega(\bbar\nu)$ as the union of the Dehn twist orthants $\gamma(C^\omega(\bbar\mu, \bbar W, \bbar r^{\hat\jmath}_i))$, where the multi-index $\hat\jmath$ varies over $\{1,2\}^{\xi(S)}$. 

First we describe a way to normalize representations of orthants by 
requiring that Dehn twists be transparently represented. Given a cube sequence $C(\bbar \mu, \bbar W, \bbar r)$ and a component $\bbar W_i$ we say that \emph{$r^\omega_i$ is a twist direction}\index{twist!direction} in the $\reals$--tree $\AM(\bbar W_i)$ if $r^\omega_i$ has positive length and there exists a sequence of linearly growing twist segments $\bbar s_i$ in $\TM(\bbar W_i)$ such that $r^\omega_i$ and $s^\omega_i$ have the same germ in $\AM(\bbar W_i)$; equivalently, one can truncate sublinearly growing initial segments of $r^n_i$ so that what is left has linearly growing initial subsegments that are twist segments. A cube sequence $C(\bbar \mu, \bbar W, \bbar r)$ is said to be \emph{twist normalized}\index{twist!normalized} if for each~$i$ the segment $r^\omega_i$ has positive length, and $r^\omega_i$ is a twist direction if and only if $\bbar W_i$ is an annulus. The dimension of a twist normalized asymptotic cube $C^\omega(\bbar \mu, \bbar W, \bbar r)$ equals the number of components of $\bbar W$, which is therefore well-defined independent of the choice of a twist normalization. Lemma~\ref{orthant complex}~(1) shows more, namely that $\bbar W$ itself is well defined up to ultraproduct equivalence.

We claim that every orthant $O$ can be represented by a twist normalized cube sequence. To see why, consider an arbitrary representation $O = \gamma(C^\omega(\bbar \mu, \bbar W, \bbar r))$. First, for each $i$ such that $r^\omega_i$ has zero length, we can extend $\bbar \mu$ by the partial marking $\bbar r_i(0)$, and then we can drop the components $\bbar W_i$ and $\bbar r_i$ from the notation, obtaining a new cube sequence representing~$O$. Next, for each $i$ such that $r^\omega_i$ is not a twist direction, it is already true that $\bbar W_i$ is not an annulus sequence. Finally, for each $i$ such that $r^\omega_i$ is a twist direction, if $\bbar W_i$ is not already an annulus sequence then we can replace $\bbar W_i$ by an annulus sequence $\bbar V_i \esssubset \bbar W_i$, and we can replace $\bbar r_i$ by a linearly growing segment sequence $\bbar s_i$ in $\TM(\bbar V_i)$, such that the image of $s^\omega_i$ under the embedding $\AM(\bbar V) \inject \AM(\bbar W)$ has the same germ as $r^\omega_i$. The result of these replacements is a twist normalized cube sequence still representing the orthant~$O$.

Define $\OO$, the \emph{orthant complex} at $\seq x$,
\index{1aaorthantO@$\OO$, orthant complex (as a poset)}
to be the poset of all nontrivial orthants at $\seq x$ --- all orthants except for the singleton $\{\seq x\}$ --- with respect to the subset relation. We shall use junctures to show that $\OO$ has the structure of a simplicial flag complex. Recall that a simplicial complex is a {\em flag complex} if, whenever a subgraph of the 1-skeleton is isomorphic to the 1-skeleton of a simplex, it is equal to the 1-skeleton of a simplex in the given complex.

\begin{lemma}{orthant complex}
For any orthants $O_1, O_2$, and for any twist normalized representations  $O_1 = \gamma(C^\omega(\bbar \mu, \bbar V, \bbar r))$ and $O_2 = \gamma(C^\omega(\bbar \nu, \bbar W, \bbar s))$, the following hold:

\begin{enumerate}
\item $O_1 = O_2$ if and only if for $\omega$-a.e.\ $n$ the following 
hold: $W^n \isotopic V^n$, and for each $i$ the segments $r^\omega_i, s^\omega_i$ have the same germ in the $\reals$--tree $\AM(\bbar W_i) = \AM(\bbar V_i)$.

\medskip

\item $O_1 \subset O_2$ if and only if for $\omega$-a.e.\ $n$ the 
following hold: each component $W^n_i$ of $W^n$ is isotopic to a component $V^n_i$ of $V^n$ and the segments $r^\omega_i, s^\omega_i$ have the same germ in the $\reals$--tree $\AM(\bbar W_i) = \AM(\bbar V_i)$. In this case we say that $O$ is a \emph{face} of $O'$.

\medskip

\item $O_1 \intersect O_2$ is the maximal common face of $O_1$ and $O_2$.
\end{enumerate}

Moreover, the poset $\OO$ is isomorphic to the poset of simplices of a flag complex $\KK$, 
\index{1aaorthantK@$\KK$, orthant complex (as a flag complex)}
having one simplex of dimension $k-1$ for each $k$-orthant. 
\end{lemma}

\textit{Remarks.} Since equality, subset, and intersection are well-defined set theoretic operations, it follows from items~(1--3) that the face relation and the ``maximal common face'' are well-defined independent of the choice of twist normalized representations, which is not at all clear a priori.

\begin{proof} While items~(1) and~(3) formally follow from item~(2), the proofs of (1--3) will all follow by studying $O_1 \intersect O_2$ using junctures.

Let $C^\omega_1 = C^\omega(\bbar \mu,\bbar V,\bbar r)$ and $C^\omega_2 = C^\omega(\bbar\nu,\bbar W,\bbar s)$. To understand $C^\omega_1\intersect C^\omega_2$, recall from Section \ref{SectionCubesCone} that this intersection is either empty or equal to the common ultralimit of the junctures of the approximating cubes. Since both cubes contain $\seq x$, the empty case cannot occur, and we are left to study the junctures. 

Lemmas \ref{cube junctures} and \ref{cube juncture details} show that the junctures of
the approximating cubes $C^n_1= C(\mu^n,V^n,r^n)$ and
$C^n_2 = C(\nu^n,W^n,s^n)$ are themselves subcubes $C^n_{12}\subset
C^n_1$ and $C^n_{21}\subset C^n_2$, which have the form
$$ C^n_{12} = C(\mu^n,V^n,r'^n),
$$
where $r'^n$ denotes a collection consisting of a subinterval (or point) or each segment of the collection $r^n$, and similarly
$$ C^n_{21} = C(\nu^n,W^n,s'^n),
$$
where $s'^n$ is a collection of subintervals or points of $s^n$. These lemmas also produce an indexing of $V^n$ and $W^n$ and a $k \ge 0$ so that for $\omega$-a.e.\ $n$ we have: $V^n_i \essint W^n_j \ne \emptyset$ if and only if $1 \le i=j \le k$ in which case we set $U^n_i = V^n_i \essint W^n_i$; this occurs only if $r'^n_i$, $s'^n_i$ have the same length (positive or zero), and all other lengths in $r'^n$ and $s'^n$ are zero. Following our usual ultraproduct convention we can say that $\bbar V_i \essint \bbar W_j$ is nonempty if and only if $i=j \in \{1,\ldots,k\}$, in which case it is isotopic to $\bbar U_i$. 

We can and do parametrize each $r'^n_i$ in such a way that $r'^n_i(0)$ is the point nearest $r^n_i(0)$, and so the corner $\kappa(C^n_{12}) = \{\mu^n\}\times \prod r'^n_i(0)$ is the nearest corner to $\kappa(C^n_1)$. Because the limiting cube $C^\omega_{12}$ contains $\seq x$, it must be that the ultralimit $\kappa^\omega(\bbar C_{12}) $ of $\kappa(C^n_{12})$ equals $\seq x$. It follows that for each $i=1,\ldots,k$ we have $r'^\omega_i(0) = r^\omega_i(0)$, or equivalently the subsegment of $r^n_i$ from $r^n_i(0)$ to $r'^n_i(0)$ grows sublinearly. Similar comments applied to $s'^n_i$ yield a corner $\kappa(C^n_{21}) = \{\nu^n\} \times \prod s'^n_i(0)$ such that $\kappa^\omega(\bbar C_{21}) = \seq x$. A note of caution is that this does {\em not} mean that the initial segments of $r^n_i$ and $s^n_i$ overlap --- again, there can be large but sublinearly growing initial segments between $r^n_i(0)$ and $r'^n_i(0)$, and similarly for $s^n_i$ and $s'^n_i$. Note also from Lemma~\ref{cube juncture details} that for each $i=1,\ldots,k$ we have $\Length(r'^n_i) \approx \Length(s'^n_i)$, and so the sequences $\Length(\bbar r'_i), \Length(\bbar  s'_i)$ both grow linearly or both grow sublinearly. 

Hence we conclude that $O_1\intersect O_2$ can be identified with the face of $O_1$ associated to those components $\bbar W_i$ of $\bbar W$ where $\Length(\bbar r'_i)$ grows linearly, by replacing sublinearly growing segments $\bbar r'_i$ with the basepoints of $\bbar r_i$, and replacing linearly growing segments $\bbar r'_i$ with initial segments of $\bbar r_i$ that contain them. The resulting sequence of faces has an ultralimit that coincides with $C^\omega_1\intersect C^\omega_2$ in a neighborhood of $\seq x$, and hence its germ is equal to $O_1\intersect O_2$. 

Furthermore, we claim that for each $i$, $\Length(\bbar r'_i)$ grows linearly if and only if $i \in \{1,\ldots,k\}$, $U^n_i \isotopic W^n_i \isotopic V^n_i$ for $\omega$-a.e.\ $n$, and $r^\omega_i$, $s^\omega_i$ have the same positive length germ in $\AM(\bbar W_i) = \AM(\bbar V_i)$. 

Once this claim is proved, it follows that the intersection $O_1 \intersect O_2$ can be described as the face of $O_1$ associated to those components $\bbar W_i$ such that $\bbar W_i$ isotopic to a component $\bbar V_i$ of $\bbar V$ and $r^\omega_i$, $s^\omega_i$ have the same positive length germ in $\AM(\bbar W_i) = \AM(\bbar V_i)$, and $O_1 \intersect O_2$ is similarly described as a face of~$O_2$. Items~(1), (2) and~(3) are all immediate consequences of this description.

For the ``if'' direction of the claim, suppose that $\Length(\bbar r'_i)$ grows sublinearly, that $i \in \{1,\ldots,k\}$, that $U^n_i \isotopic W^n_i \isotopic V^n_i$ for $\omega$-a.e.\ $n$ --- so $\Length(\bar s'_i)$ also grows sublinearly --- and that $\Length(r^\omega_i)$, $\Length(s^\omega_i)$ are both nonzero, so $\Length(\bbar r_i)$, $\Length(\bbar s_i)$ both grow linearly. By Lemma~\ref{cube juncture details}~(2), for $\omega$-a.e.~$n$ we have $r'^n_i = s'^n_i = r^n_i \intersect s^n_i$ in the tree $\TM(V^n_i) = \TM(W^n_i)$. By truncating sublinearly growing initial subsegments of $r^n_i$ and of $s^n_i$, namely the smallest initial segments containing $r'^n_i$ and $s'^n_i$, respectively, we obtain linearly growing segments $r''^n_i, s''^n_i$ such that $r^\omega_i = r''^\omega_i$ and $s^\omega_i = s''^\omega_i$, and such that $r''^n_i$ and $s''^n_i$ have disjoint interiors, so $r''^\omega_i$, $s''^\omega_i$ have distinct germs in $\AM(\bbar W_i) = \AM(\bbar V_i)$. 

For the ``only if'' direction, suppose that $\Length(\bbar r'_i)$ does grow linearly, implying that $i \in \{1,\ldots,k\}$ and that $\Length(\bbar s'_i)$ also grows linearly. If $U^n_i$ is an annulus then, by Lemma~\ref{cube juncture details}~(1), each of $r'^n_i$ and $s'^n_i$ is a twist segment supported by $U^n_i$, and by Dehn twist normalization it follows that $U^n_i \isotopic W^n_i \isotopic V^n_i$. If $U^n_i$ is not an annulus then we also have $U^n_i \isotopic W^n_i \isotopic V^n_i$, because all these surfaces have~$\xi=1$ and $U^n_i$ is essentially contained in each. Applying Lemma~\ref{cube juncture details}~(2), it follows that $r'^\omega_i$, $s'^\omega_i$ have the same germ in $\AM(\bbar W_i) = \AM(\bbar V_i)$. But $r^\omega_i, r'^\omega_i$ have the same germ, and $s^\omega_i, s'^\omega_i$ have the same germ, so $r^\omega_i$, $s^\omega_i$ have the same germ. This completes the proof of the claim.

We have shown that the intersection of two orthants is an orthant
which is equal to a common face of the two.  Also, the poset of
nontrivial faces of a $k$-orthant, meaning all faces except the
singleton $\{\seq x\}$, is isomorphic to the lattice of nonempty sets
of components of a $k$-component surface, which is isomorphic to the
poset of faces of a $k-1$ simplex.  Having excluded the unique
$0$-orthant $\{\seq x\}$ from $\OO$, this completes the proof that
$\OO$ has the structure of a simplicial complex $\KK$ with one simplex
of dimension $k-1$ for each $k$-orthant.

Now we show that $\KK$ is a flag complex. 

Let $O_1,\ldots,O_k$ be distinct 1-orthants which represent vertices of a complete graph in~$\KK$, and so for each $i,j$ there is a 2-orthant $O_{ij}$ whose faces are $O_i$ and~$O_j$. Choose twist normalized representatives $O_i = \gamma(C^\omega(\bbar \mu_i,\bbar W_i, \bbar r_i))$ where $\bbar W_i$ has a single component. Choose twist normalized representatives $O_{ij} = \gamma(C^\omega(\bbar\nu_{ij},\bbar V_{ij},\bbar s_{ij}))$ where $\bbar V_{ij}$ has two components $\bbar V_i$ and $\bbar V_j$. By item~(2) of the lemma, for $\omega$-a.e.~$n$ and each $i,j$ the surfaces $W^n_i$, $W^n_j$ are isotopic to distinct components of $V^n_{ij}$, and we may choose the notation so that $W^n_i \isotopic V^n_i$, $W^n_j \isotopic V^n_j$. It follows that $W^n_i$ and $W^n_j$ are disjoint and nonisotopic for $n$ in a set $I_{ij}$ of full $\omega$-measure. The intersection  $\intersect I_{ij}$ over all $(i,j)$ still has full $\omega$-meaure, so we conclude that $W^n_1,\ldots,W^n_k$ are pairwise disjoint and non-isotopic for $\omega$-a.e.\ $n$, and we obtain an essential subsurface $\bbar W = \bbar W_1\union \cdots \union \bbar W_k$. Let $\bbar\sigma$ be a marking sequence on $\bbar W^c$ defined as the projection $\pi_{\MM(\bbar W^c)}(\bbar x)$, and let $\bbar r = (\bar r_1,\ldots,\bbar r_k)$. Then we obtain an orthant
$$ O = \gamma({C^\omega}(\bbar \sigma,\bbar W,\bbar r)).
$$

We need to check  that the corner of $O$, namely the limit of $\bbar\kappa =\{\bbar\sigma\}\times \prod \bbar r_i(0)$, equals $\seq x$. But this is a consequence of the quasidistance formula for $d(\kappa^n,x^n)$, in which we separate the terms
$\Tsh A{d_Z(\kappa^n,x^n)}$ according to whether $Z\esssubset (W^n)^c$, $Z\esssubset W^n$, or  $Z\pitchfork \boundary W^n$. The first type of term adds up to an estimate of $d_{\MM(W^n)^c}(\sigma^n,x^n)$, which by definition of $\bbar\sigma$ is bounded. The second type adds up to estimate the finite sum $\sum_id_{\MM(W^n_i)}(r^n_i(0),x^n)$, each of whose terms grows sublinearly since the corner of each 1-orthant $O_i$ is $\seq x$. The third type is estimated, termwise, by $\Tsh A{d_Z(\boundary W^n,x^n)}$, which sum up to estimate $d(x^n,\QQ(\boundary W^n))$, by Lemma \ref{Q product structure}. This again grows sublinearly since $\seq x\in \intersect_i \QQ_\omega(\boundary \bbar W_i) = \QQ_\omega(\boundary\bbar W)$. 
We conclude that $d(\kappa^n,x^n)$ grows
sublinearly, so $\kappa^\omega = x^\omega$. 

This tells us that $O\in \OO$. It is clear by construction and item~(2) of the lemma that 
$O_i$ are the vertices of $O$. This completes the proof. 
\end{proof}

The following observation will be used later when we give a topological characterization of Dehn twist flats.

\begin{corollary}{CorOrthantGivesFlat}
The germ of every Dehn twist $\xi$-orthant is contained in a unique Dehn twist flat germ.
\end{corollary}

\begin{proof}
Consider two Dehn twist flat germs $Q_\omega(\bbar\nu)$, $Q_\omega(\bbar\nu')$ whose intersection contains a common Dehn twist $\xi$-orthant $O$. Let $\bbar W$, $\bbar W'$ be the sequences of annulus neighborhoods of the pants decomposition sequences $\bbar\nu$, $\bbar \nu'$. We obtain two twist normalized expressions $O = \gamma(C^\omega(\bbar\mu, \bbar W, \bbar r)) = \gamma(C^\omega(\bbar\mu', \bbar W', \bbar r'))$. By Lemma~\ref{orthant complex}~(1) the sequences $\bbar W$, $\bbar W'$ are $\omega$-equivalent, which implies that $Q_\omega(\bbar\nu) = Q_\omega(\bbar\nu')$.
\end{proof}

\subsection{Applying local finiteness}

We continue to fix the base point $\seq x \in \AM(S)$.  Consider the
subset of $\GG$ consisting of germs at $\seq x$ of submanifolds of
$\AM$ of dimension $\xi=\xi(S)$.  This subset generates a sublattice
$\FF\subset \GG$ by taking finite unions and intersections.  Since
germs of manifolds are topologically characterizable, it follows that
germs in $\FF$ are topologically characterizable.  Our goal in this
section is to produce finer topological properties in $\FF$, in order
to yield a topological characterization of germs of Dehn twist
flats given in Corollary~\ref{top char of twist germs}.

Let $\hhat\OO\subset \GG$ be the sublattice generated from the orthant
complex $\OO$ at $\seq x$ by taking finite unions and intersections.
By Lemma~\ref{orthant complex}, $\hhat\OO$ is isomorphic to the
lattice of finite subcomplexes of the simplicial complex $\KK$.

The manifold local finiteness theorem, Theorem \ref{local cube
finiteness}, will imply:

\begin{lemma}{F subset hatO}
$\FF \subset \hhat\OO$.
\end{lemma}

\begin{proof} 
Let $M$ be a manifold of dimension $\xi=\xi(S)$ passing through
$\seq x$. Theorem \ref{local cube finiteness} states that there is a
neighborhood $U$ of $\seq x$ such that $M\intersect U$ is contained in a finite union of cubes. After subdivision and possibly replacing $U$ by a smaller open set containing $\seq x$, we obtain a finite collection of cubes $C_1,\ldots,C_I$, each having $\seq x$ as a corner, whose union contains $M \intersect U$. Applying local compactness of $M$, choose an open set $V \subset U$ containing $\seq x$ such that the closure of $M \intersect V$ is a compact subset of $M \intersect U$.

Suppose $M \intersect V$ has nontrivial intersection with the interior of a cube $C_i$ of dimension~$\xi$. For $j \ne i \in \{1,\ldots,I\}$ the intersection $C_i \intersect C_j$ is contained in the boundary of $C_j$, by Lemma \ref{orthant complex}, and so $\interior(C_i)$ is disjoint from the closed set $C_j$. It follows that $M \intersect \interior(C_i)$ is open in~$M$. Combined with invariance of domain it follows that $\interior(C_i) \intersect V \subset M$, and so $C_i \intersect V \subset M$. 

Suppose next that $M \intersect V$ has nontrivial intersection with a cube $C_j$ of dimension $< \xi$. Again by invariance of domain, $M \intersect V$ must meet the interior of some cube $C_i$ of dimension $\xi$ having $C_j$ as a face, and so $C_j \intersect V \subset C_i \intersect V \subset M$. 

We conclude that any germ of a manifold is {\em equal to} a finite union of orthants, and hence $\FF \subset \hhat\OO$.
\end{proof}

To clarify the structure of $\FF$, we introduce some more objects.

Consider a top dimensional orthant $O$ and a twist normalized 
representation $O = \gamma(C^\omega(\bbar \mu, \bbar W, \bbar r))$. 
Let the components of $\bbar W$ be  $\bbar W_1,\ldots,\bbar W_\xi$. 
In each $\bbar W_i$ we have a ray $\seq r_i$ in the associated 
$\reals$--tree $\AM(\bbar W_i)$. Actually we only need to consider the {\em germ} of a ray, but we will still denote it~$\seq r_i$. A component $\bbar W_i$ is called a {\em boundary annulus} if for $\omega$-a.e.\ $n$ the surface $W^n_i$ is an annulus homotopic to the boundary of another component $W^n_j$, necessarily of complexity~1. Let $b(O) = b(\bbar W)$ denote the number of boundary annuli, which is a well-defined function of $O$ by Lemma~\ref{orthant complex}~(1). Note that $b(O)=0$ if and only if all components of $\bbar W$ are annuli, if and only if $O$ is a Dehn twist orthant, if and only if each $\bbar r_i$ is a twist direction --- the first ``if and only if'' is a consequence of top dimensionality, the second is a matter of definition, and the last is a consequence of twist normalization. 

Fix $j \in \{1,\ldots,\xi\}$, and consider the $j^{\text{th}}$ codimension~1 face of $O$, obtained by restricting the ray $\seq r_j$ to its initial point. If $O' = \gamma(C^\omega(\bbar \mu', \bbar W', \bbar r'))$ is another twist normalized top dimensional orthant meeting $O$ along the $j^{\text{th}}$ codimension~1 face, then by applying Lemma~\ref{orthant complex}~(3) it follows that $O'$ is obtained from $O$ in one of the following ways: only the component $\bbar W_j$ is changed; or $\bbar W'$ and $\bbar W$ are equivalent and only the ray germ $\seq r_j$ is changed (all of these changes are up to $\omega$-equivalence). If $\bbar W_j$ is of complexity 1 then
there are infinitely many different choices for $O'$, for example 
there are infinitely many different ray germs $\seq r'_j$ to choose from in the $\reals$--tree $\AM(\bbar W_j)$. If $\bbar W_j$ is a nonboundary annulus then there are also infinitely many different choices for $O'$, for instance there are infinitely many annuli $\bbar W'_j$ that can replace $\bbar W_j$. However, if $\bbar W_j$ is a boundary annulus then the only change we can make is to replace $\seq r_j$ by the unique opposite ray germ $-\seq r_j$ in the line $\AM(\bbar W_j)$. 

We conclude that, along each of the $b=b(O)$ codimension-1 faces of $O$ associated to boundary annuli, there is a unique orthant adjacent to $O$. It follows that any manifold germ $M$ containing $O$ must contain all of these unique neighboring orthants. Furthermore each of these orthants still has the same defining surface $\bbar W$ and the same set of $b$ boundary annuli, and for all the corresponding faces the unique neighboring orthants must be included. We conclude that all $2^b$ orthants obtained in this way must be contained in the germ $M$. We call this set a {\em lune},\index{lune} and refer to the number $b$ as its {\em  rank}.\index{rank of lune}  We note that it is naturally identified with a Euclidean spherical lune $\R^b \times (\R^+)^{\xi-b} \intersect {\mathbb S}^{\xi-1}$, with its subdivision into Euclidean orthants (i.e., spherical simplices). 

\begin{lemma}{lune characterization}
Lunes are precisely the minimal $\xi$-dimensional elements
of the lattice~$\FF$. 
\end{lemma}

\begin{proof}
Consider $L$ a lune of rank $b$. We first show that $L \in \FF$, and that $L$ is a minimal $\xi$-dimensional element of $\FF$. As discussed above, $L$ is a union of $2^b$ orthants, and without loss of generality the associated (germs of) rays for these orthants are of the form
$$ {\seq r}_1^{j_1},\ldots,{\seq r}_b^{j_b},{\seq r}_{b+1}^1,\ldots,{\seq r}_\xi^1,
$$
where  $\bbar W_1,\ldots,\bbar W_b$ are the boundary annuli in a
decomposition $\bbar W$, and for each $i\in\{1,\ldots,b\}$ we have $j_i\in \{1,2\}$ and the ray germs  ${\seq r}_i^1,\,{\seq r}_i^2$ are opposite pairs. 

For each $i\in \{b+1,\ldots,\xi\}$ we define the following objects. First, choose $\seq r^2_i$ to be a ray in $\AM(\bbar W_i)$ with germ distinct from $\seq r^1_i$; if $\xi(\bbar W_i) = 1$ make sure that $\seq r^2_i$ is not a twist direction. Next, let $\bbar V_i$ be $\bbar W_i$ if
$\xi(\bbar W_i) = 1$, and if $\xi(\bbar W_i) = 0$ let $\bbar V_i$ be
the unique sequence (up to the usual ultraproduct equivalence) of $\xi=1$ subsurfaces containing $\bbar W_i$ and disjoint from all the other $\bbar W_k$. We can interpret $\seq r^1_i$, $\seq r^2_i$ as ray germs in $\AM(\bbar V_i)$ via the natural embedding $\AM(\bbar W_i)\to\AM(\bbar V_i)$. Denote $\seq s^1_i = \seq r^1_i$. Finally, choose $\seq s^2_i$ to be a ray germ in $\AM(\bbar V_i)$ which shares its basepoint with $\seq r^1_i$, $\seq r^2_i$ but is distinct from both. Let $\bbar W[i]$ be the subsurface sequence obtained from $\bbar W$ by replacing $\bbar W_i$ by $\bbar V_i$. 

Now for each $i\in \{b+1,\ldots,\xi\}$ and each tuple
$\hat \jmath=(j_1,\ldots,j_\xi)\in\{1,2\}^\xi$, consider the 
orthant $O[i]({\hat\jmath})$ formed from  $\bbar W[i]$ and the ray germs 
$$ \seq r_1^{j_1}, \ldots, \seq r_b^{j_b}, \ldots, \seq r_{i-1}^{j_{i-1}}, \seq s_i^{j_i}, \seq r_{i+1}^{j_{i+1}} \ldots
$$
in other words, we use $\seq r_k^{j_k}$ for all $k$ except $k=i$, where
we use $\seq s_i^{j_i}$. Let 
$$
M[i] = \bigcup_{\hat \jmath} O[i](\hat \jmath).
$$
This is a manifold germ, and our lune $L$ is the intersection 
$$
L = M[b+1]\intersect\cdots\intersect M[\xi].
$$
This shows that $L$ is in $\FF$. Since any $\xi$-dimensional element $C$ of $\FF$ contained in $L$ must contain a top dimensional orthant $O\subset L$, it follows from the paragraph before Lemma~\ref{lune characterization} that $C$ contains $L$, and so $L$ is minimal. Each lune is therefore a minimal $\xi$-dimensional elements of $\FF$.

Now let $C$ be any minimal $\xi$-dimensional element of $\FF$.
Then $C$ must contain a $\xi$-dimensional orthant $O$, by Lemma \ref{F subset hatO}. By the discussion above, $C$ must contain the lune $L$ determined by $O$, and by the minimality of $C$, we have
$C=L$. Each minimal $\xi$-dimensional element of $\FF$ is therefore a lune.
\end{proof}

We will also make use of the following:

\begin{lemma}{LemmaLuneManifold}
For each lune $L$ of rank $b$ there exists a manifold germ containing $L$ which is a union of $2^{\xi - b}$ distinct lunes of rank $b$, no two of which have an orthant in common.
\end{lemma}

\begin{proof} We borrow the notation of Lemma~\ref{lune characterization}, but in this case the construction is somewhat easier. For each tuple $\hat \jmath = (j_{b+1},\ldots,j_\xi) \in \{1,2\}^{\xi - b}$ let $L(\hat\jmath)$ be the lune of rank $b$ which is the union of the orthants associated to ray germs 
$$ \seq r^{j_1}_1, \ldots, \seq r^{j_b}_b, \seq r^{j_{b+1}}_{b+1},\ldots, \seq r^{j_\xi}_\xi
$$
where $(j_1,\ldots,j_b)$ varies freely in $\{1,2\}^b$. The union of these lunes, over all $\hat \jmath \in \{1,2\}^{\xi-b}$, is a manifold germ.

Finally, note that distinct lunes have no orthants in common, by the paragraph before Lemma~\ref{lune characterization}. 
\end{proof}

Since $\hhat\OO$ is isomorphic to the poset of finite subcomplexes of the $(\xi-1)$-dimensional simplicial complex ${\mathcal K}$  of Lemma \ref{orthant complex}, each element $C\in\hhat\OO$ determines a simplicial $(\xi-1)$-chain with $\Z_2$-coefficients in ${\mathcal K}$, namely the formal sum of the simplices corresponding to the top dimensional orthants appearing in $C$.  In what follows we will conflate the chain with $C$  when convenient. Given two chains $\alpha,\beta\in C_{\xi-1}({\mathcal K})$, we say that {\em  $\alpha$ is part of $\beta$} if $\beta=\alpha+\alpha'$ where the chains $\alpha$ and $\alpha'$ have no simplices in common.

We let ${\mathcal L}$ denote the collection of lunes, which by Lemma~\ref{lune characterization} are topologically characterized. Our next goal is a characterization of the rank of lunes as a function on $\LL$:

\begin{lemma}{lune rank}
 The rank is the unique function 
$f\co {\mathcal L}\rightarrow\{0,\ldots,\xi\}$
with the following property:
\begin{itemize}
\item If $C$ is a lune and if $f(C)\leq b \in \{0,\ldots,\xi\}$, then $f(C)=b$ if and only if $C$ is part of a nonzero cycle
$$\ \sum_{i=1}^{2^{\xi-b}}\;C_i,
$$
such that each $C_i$ is a lune satisfying $f(C_i)\leq b$.
\end{itemize}
 \end{lemma}

To prove this, we will need a lemma about
flag complexes:

\begin{lemma}{flag has big cycles}
Every nontrivial reduced
$\Z_2$ $n$-cycle 
 in an $n$-dimensional flag complex has cardinality at
least $2^{n+1}$. 
\end{lemma}

\begin{proof}
The lemma obviously holds for $0$-dimensional flag complexes, since the support of a nontrivial reduced $0$-cycle must contain at least two vertices.

Assume inductively that $n=\dim X>0$, and
 that the lemma holds
for flag complexes of dimension $<n$.
We first observe that the link of any vertex is an 
$(n-1)$-dimensional flag complex,
and hence by the induction assumption, the lemma holds for
links.

Let $M$ be a $\Z_2$ $n$-cycle in $X$. Consider two adjacent
$n$-simplices $\sigma_1,\sigma_2$ meeting at a codimension 1 face $\tau$. 
Let $v_i$ be the vertex of $\sigma_i$ complementary to~$\tau$. The
link of $v_i$ in $M$ is a $\Z_2$ $(n-1)$-cycle, hence by the assumption has
cardinality at least $2^n$. The lemma would follow if we show that the
stars of $v_1$ and $v_2$ do not have common $n$-simplices.

Suppose there is such a simplex. Then $v_1$ and $v_2$ must be joined
by an edge $e$. Now the abstract join $\tau * e$ is an $(n+1)$-simplex all of
whose edges are in the complex $X$. Since $X$ is a flag complex, it
must contain an $(n+1)$-simplex, but this contradicts $\dim X=n$.
\end{proof}

\begin{proof}[Proof of Lemma \ref{lune rank}.] We will refer to the  property stated in Lemma~\ref{lune rank} as {\em Property S}.

We first show that the rank function has Property S. If $C$ is a lune of rank~$b$ then, by Lemma~\ref{LemmaLuneManifold}, $C$ is part of a cycle consisting of $2^{\xi-b}$ lunes of rank $b$. Conversely, suppose $b\in \{0,\ldots,\xi\}$, $C$ is a lune of rank $\leq b$, and $C$ is part of a nonzero cycle $\sum_{i=1}^{2^{\xi-b}}C_i$ where each $C_i$ is a lune of rank $\leq b$, so $C_i$ is composed
of $2^{\rank(C_i)}\leq 2^b$ orthants, and the entire cycle is composed of $\le 2^{\xi}$ orthants. Nonzero cycles require at least $2^\xi$ orthants by Lemma \ref{flag has big cycles}, which implies that
$\rank(C_i)=b$ for all $i$.  Since $C$ is part of $\sum_i C_i$,
the intersection $C\cap C_j$ must contain a top dimensional
orthant for some $j$.  The minimality property for lunes, Lemma~\ref{lune characterization}, implies that $C=C_j$ and hence $\rank(C)=\rank(C_j)=b$. This shows that $\rank$ has Property S.

Now suppose $f:\LL\ra \{0,\ldots,\xi\}$ has Property S,
but is not equal to $\rank$.  Let $b$ be the maximum of the integers
$\bar b\in\{0,\ldots,\xi\}$ such that $f^{-1}(\bar b)\neq\rank^{-1}(\bar b)$.

Suppose $C$ is a lune of rank $b$.  Then $C$ belongs to a 
cycle $\sum_{i=1}^{2^{\xi-b}}C_i$ where $\rank(C_i)=b$. By choice of $b$ we have $f(C_i)\leq b$. Hence by Property S, we
get $f(C)=b$. Thus $\rank^{-1}(b)\subset f^{-1}(b)$.

Now suppose $C\in f^{-1}(b)$.  Then $C$ belongs to a cycle
$\sum_{i=1}^{2^{\xi-b}}C_i$ where $f(C_i)\leq b$.  By the choice
of $b$, we have $\rank(C_i)\leq b$ for all $i$, and by Lemma 
\ref{flag has big cycles} we get $\rank(C_i)=b$ for all $i$. 
We conclude as above, using Lemma~\ref{lune characterization}, that $C=C_j$ for some $j$, and hence $\rank(C)=b$.  Thus $f^{-1}(b)\subset \rank^{-1}(b)$.  This contradicts the choice of~$b$.
\end{proof}

All the pieces are now in place for the main result of this section: 

\begin{corollary}{top char of twist germs}
There is a topological characterization of Dehn twist orthant germs and Dehn twist flat germs in $\AM(S)$.
\end{corollary}

\begin{proof}
Membership in the lattice $\FF$ is topologically characterized, as is the dimension function on $\FF$.  Therefore Lemmas \ref{lune characterization} and \ref{lune rank} give topological characterizations of lunes and lune rank. Dehn twist $\xi$-orthants are the lunes of rank $0$, so these are also topologically characterizable. Dehn twist orthants of arbitrary dimension are topologically characterized as intersections of sets of Dehn twist $\xi$-orthants.

It remains to topologically characterize Dehn twist flat germs --- those configurations of $2^\xi$ orthants associated to a Dehn-twist flat through $\seq x$. This boils down to finding a topological characterization of those antipodal pairs of Dehn twist $1$-orthants (also known as Dehn twist vertices) that occur in Dehn twist flat germs.  

Consider a twist normalized Dehn twist vertex $O = \gamma(C^\omega(\bbar \mu, \bbar W, \bbar r))$, so $W^n$ is connected and is an annulus for $\omega$-a.e.\ $n$. The \emph{antipodal point} of $O$ is defined to be $-O = \gamma(C^\omega(\bbar\mu, \bbar W, - \bbar r))$, where $\bbar r$ and $-\bbar r$ have the same initial point but opposite directions in the line $\TM(\bbar W)$. Note that the antipodal point of $O$ is well-defined, for suppose that $O = \gamma(C^\omega(\bbar \mu', \bbar W', \bbar r'))$ is another twist normalized expression. By Lemma~\ref{orthant complex}~(1), for $\omega$-a.e.\ $n$ we have $W^n \isotopic W'^n$ and $\seq r$, $\seq r'$ have the same direction in the line $\AM(\bbar W) = \AM(\bbar W')$. Moreover, the asymptotic partial markings $\seq \mu, \seq \mu' \in \AM(\bbar W^c)$ are each equal to the projection of~$\seq x$. It follows that the two asymptotic cubes $C^\omega(\bbar \mu, \bbar W, -\bbar r)$, $C^\omega(\bbar \mu', \bbar W', -\bbar r')$ are equal in $\AM(S)$ and so the two expressions for the antipodal point $\gamma(C^\omega(\bbar \mu, \bbar W, -\bbar r))$ and $\gamma(C^\omega(\bbar \mu', \bbar W', -\bbar r'))$ define the same Dehn twist vertex $-O$.

This shows that the relation of ``antipodal point'' gives a decomposition of the set of Dehn twist vertices into pairs. We now show that antipodal pairs are topologically characterizable.

Consider a lune of positive rank $b>0$. The corresponding decomposition $\bbar W$ contains $b$ boundary annuli associated to which are $b$ antipodal pairs of Dehn twist vertices that span a sphere of dimension $b-1$ in the orthant complex, which is subdivided in the standard way into $2^b$ simplices. 

These lune boundary spheres are topologically characterizable: they are precisely the $(b-1)$-dimensional spheres which may be obtained as the intersection of two lunes of rank $b$. Moreover the simplicial decomposition of such a sphere is topologically characterizable, since the simplices are Dehn-twist orthants. Hence the pairs of antipodal vertices of lune boundary spheres are topologically characterizable in terms of this simplicial structure.

Now any Dehn twist vertex of $\OO$ can be placed into such a lune boundary sphere, simply by extending its defining annulus to a decomposition where it is a boundary annulus. Thus antipodal pairs of Dehn twist vertices may be characterized topologically as those vertex pairs which may be embedded as antipodal vertices in a triangulated lune boundary sphere.

To summarize, germs of Dehn twist flats are topologically characterized as those flag subcomplexes $C \subset \OO$ whose vertex set decomposes into a union of $\xi$ antipodal pairs of (twist normalized) Dehn twist vertices 
$$\pm O_i = \gamma(C^\omega(\bbar \mu_i, \bbar W_i, \pm \bbar r_i)), \qquad i=1,\ldots,\xi
$$
such that any choice of one vertex from each pair spans a $\xi$-orthant: $C$ is the union of those $2^\xi$ orthants; by Lemma~\ref{orthant complex}~(2) applied to any such orthant the annuli $\bbar W_1,\ldots,\bbar W_\xi$ form a sequence of pants decompositions $\bbar W$; and the Dehn twist flat germ associated to $\bbar W$ is~$C$.
\end{proof}

\subsection{Characterizing Dehn twist flats}
\label{SectionCharTwistFlats}

We conclude Section~\ref{orthant defs} with the following observation, which gives a local topological characterization of Dehn twist flats.

\begin{lemma}{top char of twist flats}
Suppose $E\subset \AM(S)$ is a connected top-dimensional manifold,
and that for  every
$\seq x\in E$, the germ of $E$ at $\seq x$ is the germ
of a Dehn twist flat.  Then $E$ is contained in a Dehn twist flat.
If in addition $E$ is a closed subset of $\AM(S)$, then 
$E$ is a Dehn twist flat.
\end{lemma}

\begin{proof}  Pick any Dehn twist flat $E'\subset\AM$ such that the interior in $E$ of $E\cap E'$ is a nonempty set $U \subset E$. Suppose $\seq x\in E$ lies in the closure of $U$ in $\AM$.  Since $E'$ is a closed subset of $\AM(S)$, we have $\seq x\in E\cap E'$. By the definition of $U$, the germ of $E \cap E'$ at $\seq x$ has dimension $\xi$, so by Corollary~\ref{CorOrthantGivesFlat} it follows that $E,E'$ have the same germ at $\seq x$, and we conclude that $\seq x\in U$.
Thus $U$ is an open and closed subset of $E$. Since $E$ is connected, we have $E=U\subset E\cap E'\subset E'$.

If $E$ is a closed subset of $\AM(S)$, then $E\cap E'$ is 
open and closed in $E'$, and hence $E\cap E'=E'$.
\end{proof}

\section{Finishing the proofs}
\label{endgame}

\newcommand\Aut{{\rm Aut}}

We are now ready to prove our main theorems on quasi-isometric
classification and rigidity. The proof will follow the general sketch
from the introduction. Let us first state the theorems in their complete form:

\State{Theorem \ref{QI classification}}{Classification of Quasi-Isometries}{
 Suppose that $\xi(S) \ge
2$.

If $S \ne S_{1,2}$ then quasi-isometries of $\MCG(S)$ are uniformly
close to isometries induced by left-multiplication.

That is, given $K,\delta>0$ there exists $D>0$ such that, if $f\from
\MCG(S)\to\MCG(S)$ is a $(K,\delta)$-quasi-isometry then there exists
$g\in\MCG(S)$ such that $$ d(f(x),L_g(x)) \le D \quad \text{for all $x
\in \MCG(S)$}
$$
where $L_g$ is left-multiplication by $g$.

If $S = S_{1,2}$ then the same result holds if we replace $L_g$ by a
quasi-isometry of $\MCG(S_{1,2})$ induced by an element $g \in
\MCG(S_{0,5})$ via the standard index~5 embedding
$\MCG(S_{1,2})/Z(\MCG(S_{1,2})) \hookrightarrow \MCG(S_{0,5})$.
}

Here $Z(G)$ denotes the center of $G$, which is trivial for all 
mapping class groups with $\xi(S)>2$; furthermore, in the 
cases where it is non-trivial, it is finite. 

\State{Theorem \ref{QI rigidity}}{Quasi-Isometric Rigidity}{
If $\Gamma$ is a finitely generated group 
quasi-isometric to $\MCG(S)$, then there exists a finite-index 
subgroup $\Gamma'<\Gamma$ and a homomorphism 
$$\Gamma'\to\MCG(S)/Z(\MCG(S))$$
 with finite kernel and
finite index image. 
}

In the cases of mapping class groups of complexity less than $2$, 
quasi-isometric rigidity is either trivial or well known. On the 
other hand, the analogue of Theorem~\ref{QI classification} fails in 
the case of complexity one, as there are quasi-isometries of the free 
group which are not a bounded distance from isometries.

An immediate consequence of Theorem \ref{QI classification} is the
following characterization of the quasi-isometry group $\QI(\MCG(S))$,
i.e. the group of quasi-isometries of $MCG(S)$ modulo those that are
finite distance from the identity.
\begin{corollary}{QIGroup} If $S$ has complexity at least $2$ then the natural homomorphism
$$\MCG(S)/Z(\MCG(S)) \to \QI(\MCG(S))$$ is an isomorphism except
when $S=S_{1,2}$, in which case it is an isomorphism to a subgroup of 
index~5.
\end{corollary}
The proof of this corollary is embedded in the proof of Theorem
\ref{QI rigidity} in \S\ref{QIR proof} below.

\subsection{Preservation of asymptotic Dehn twist flats}

We begin the proof by showing that the local topological characterization of
Dehn twist flats in the asymptotic cone, given in Lemma~\ref{top char of twist flats}, implies a
global characterization: 

\begin{theorem}{asymptotic flat theorem}
If $\xi(S)\geq 2$, 
any homeomorphism $f\co\AM(S) \to \AM(S)$ 
 permutes the 
Dehn twist flats in $\AM(S)$. 
\end{theorem}

\begin{proof}
By Corollary \ref{top char of twist germs}, any homeomorphism must
preserve the set of Dehn twist flat germs in $\AM(S)$ (with arbitrary
basepoints) since they are topologically characterized. It follows
that, at every point in the image $f(E)$ of a Dehn twist flat $E$, 
its germ is equal to the germ of a Dehn
twist flat. Lemma \ref{top char of twist flats} therefore implies that $f(E)$ is
itself a Dehn twist flat. 
\end{proof}

\subsection{Coarse preservation of Dehn twist flats}
\label{asymptotic to coarse}

We next descend to the group itself, where we  show that a quasi-isometry of
$\MM(S)$ coarsely preserves Dehn twist flats, in a uniform sense.

\begin{theorem}{coarse flat theorem}
If $\xi(S) \geq 2$, then given $K\ge 1$ and $C \ge 0$ there exists $A$ such that, if $f\co\Mod(S)\to\Mod(S)$ is a $(K,C)$-quasi-isometry and $E$ is a Dehn twist flat in $\Mod(S)$ then there exists a unique Dehn twist flat $E'$ such that the Hausdorff distance between $f(E) $ and $E'$ is at most $A$. 
\end{theorem}

\begin{proof} Uniqueness of $E'$ follows from the fact that distinct Dehn twist flats have infinite Hausdorff distance, an immediate consequence of Proposition~\ref{Q product structure} and Lemma~\ref{LemmaFlatSetTh}.

The existence proof is essentially an argument by contradiction, using Theorem \ref{asymptotic flat theorem}.  If there is no uniform control of the Hausdorff distance between Dehn twist flats and their quasi-isometric images, while on the other hand in every limiting situation the Dehn twist flats are preserved in the asymptotic cone, then in a sequence of counterexamples we can carefully select basepoints and scales to get configurations in which the image of a Dehn twist flat is simultaneously very close to two distinct Dehn twist flats.  This contradicts the fact that distinct Dehn twist flats look different at all scales (Lemma~\ref{asymmetric linear spreading}, which is a consequence of Lemma~\ref{juncture}).

We will find the following variation of the Hausdorff metric useful. Given subsets $A,B$ of a
metric space $X$ and a point $p\in X$,
for each $r>0$ we define
$$
D_{r,p}(A,B) = \inf\{ s \ge 0 \suchthat A \intersect \NN_r(p) \subset
\NN_s(B) \quad\text{and}\quad B \intersect \NN_r(p) \subset \NN_s(A) \} 
$$
Notice that, if $\NN_r(p) \subset \NN_{r'}(p')$, then
$D_{r,p} \le D_{r',p'}$.
This is not quite a distance function --- it fails the triangle
inequality --- but it does give a useful criterion for equality of
ultralimits. 

In the following lemma we consider a sequence $(X_i,p_i)$
of  based metric spaces with ultralimit $X_\omega$. No rescaling is
assumed here; in applications below, $X_i$ will be the rescaled
$\Mod(S)$. 

\begin{lemma}{LemmaUltraClose}
Given $(A_i)$, $(B_i)$ two sequences of closed subsets,
$A_\omega = B_\omega$ if and only if 
for each basepoint $(q_{i})$ and 
for each $r \ge 0$ the ultralimit of $D_{r,q_i}(A_i,B_i)$ equals zero. 
\end{lemma}

\begin{proof}
Suppose that for some  
$(q_{i})$ and some 
$r>0$ we have
$D_{r,q_i}(A_i,B_i)\to_\omega \epsilon\in(0,\infty]$. 
Choose $\eta\in(0,\epsilon)$. 
It follows that, for $\omega$-a.e. $i$, one of the following is true:
\begin{description}
\item[(1)] $A_i\intersect \NN_r(q_i) \not\subset \NN_\eta(B_i)$
\item[(2)] $B_i\intersect \NN_r(q_i) \not\subset \NN_\eta(A_i)$
\end{description}
Furthermore, either (1) is true for $\omega$-a.e.\ $i$, or (2) is true for $\omega$-a.e.\ $i$; let us assume the former. From this we conclude that there is a sequence 
$$x_i \in A_i\intersect \NN_r(q_i) \setminus \NN_\eta(B_i)
$$
for which $x_\omega\in A_\omega$ but the distance between 
$x_{\omega}$ and $B_\omega$
is at least $\eta$. Hence $A_\omega\ne B_\omega$.

Suppose next that  $\lim_\omega D_{r,q_i}(A_i,B_i)=0$ for all 
$(q_{i})$ and all $r \ge 0$. 
To prove that $A_\omega \subset B_\omega$, consider
$x_\omega \in A_\omega$ represented by a sequence $(x_i)$ at bounded 
distance from $(q_{i})$, so 
there exists some $r \ge 0$ such that $x_i \in A_i \intersect
\NN_r(q_i)$ $\omega$-almost surely. For any integer $k>0$ it follows that $x_i
\in \NN_{1/k}(B_i)$ $\omega$-almost surely, so we may choose a sequence $y^k_i
\in B_i$ such that $d_i(x_i,y^k_i) < 1/k$ $\omega$-almost surely, and therefore
$y^k_\omega \in B_\omega$ and $d(x_\omega,y^k_\omega) \le 1/k$. The
sequence $y^k_\omega$ therefore converges to $x_\omega$, but this
sequence is in the closed set $B_\omega$, proving that $x_\omega \in
B_\omega$. A symmetric argument proves that $B_\omega \subset
A_\omega$.  
\end{proof}

We turn now to the proof of Theorem~\ref{coarse flat theorem}. Suppose that the theorem is false. Then we may fix $K \ge 1$, $C \ge 0$ so that the following is true: for any $A \ge 0$ there is a $(K,C)$-quasi-isometry $f \from \Mod(S) \to \Mod(S)$, and a Dehn twist flat $\Flat$, such that for any Dehn twist flat $\Flat'$, the Hausdorff distance between $f(\Flat)$ and $\Flat'$ is greater than $A$.

From this symmetric statement we make the further asymmetric
conclusion that for each $s>0$ there is a $(K,C)$-quasi-isometry $f
\from \Mod(S) \to \Mod(S)$ and a Dehn twist flat $F$ such that for all
Dehn twist flats $F'$ we have
$$f(F) \not\subset \NN_s(F')
$$ 
For if not, then there exists $s>0$ such that for all $(K,C)$-quasi-isometries $f$ and all Dehn twist flats $F$ there exists a Dehn twist flat $F'$ such that $f(F) \subset \NN_s(F')$. The closest point projection $\pi$ from the $(K,C)$-quasiflat $f(F)$ to the Dehn twist flat $F'$ moves points a distance at most $s$ and can therefore be regarded as a $(K',C')$-quasi-isometry from $\reals^n$ to $\reals^n$ for constants $K',C'$ that depend only on $K,C,s$. Now any $(K',C')$-quasi-isometric embedding from $\reals^n$ into $\reals^n$ is a $(K',C'')$-quasi-isometry where $C''$ depends only on $K', C'$, and $n$; see, for example, \cite{kapovich:book}. It follows that $\pi$ is uniformly onto, meaning that there exists a constant $B$ depending only on $K',C'$ such that $F'$ is in the $B$ neighborhood of $\pi(f(F))$, and so $F' \subset\NN_{s+B}(f(F))$. This shows that $f(F)$ and $F'$ have Hausdorff distance at most $s+B$, which is a contradiction for $A > s+B$.

Fix a sequence $s_i$ diverging to $+\infinity$, a sequence of
$(K,C)$-quasi-isometries $f_i \from \Mod(S) \to \Mod(S)$, and a sequence
of Dehn twist flats $F_i$, such that for all $i$ and all Dehn twist
flats $F'$ we have
\begin{equation}
\label{EqQFlatNotFlat}
f_i(F_i) \not\subset \NN_{s_i}(F')
\end{equation}
Since there are finitely many $\Mod(S)$-orbits of Dehn twist flats, 
by pre-composing with elements of $\Mod(S)$  and extracting a
subsequence we may assume that the
$F_i$ take a constant value $F$. Fix a base point $p_0 \in F$. By
post-composing with elements of $\Mod(S)$ we may assume that
$f_i(p_0)=p_0$, and in particular $p_0 \in f_i(F)$, for all $i$. We
may therefore pass to the asymptotic cone with base point $p_0$ and
scaling sequence $s_i$ producing a bi-Lipschitz homeomorphism $f_\omega\from \AM(S) \to \AM(S)$ and a Dehn twist flat
$F_\omega=\lim_\omega(F)$, the asymptotic cone of $F$. Applying
Theorem \ref{asymptotic flat 
  theorem}  we obtain a Dehn twist flat
$F'_\omega = \lim_\omega F'_i$ such that $f_\omega(F_\omega) = F'_\omega$. 
By applying Lemma \ref{LemmaUltraClose} it follows that, fixing any $R>0$, 
\begin{equation}
\label{EqEquals_delta_i_s_i}
\frac{1}{s_i}D_{Rs_i,p_0}(f_i(F),F'_i) \to_\omega 0.
\end{equation}
On the other hand, \pref{EqQFlatNotFlat} implies
that there is a point $q_i \in f_i(F) - \NN_{s_i}(F'_i)$, and so for any
$r>0$ the following statement is always true:
\begin{equation}
\label{EqFar}
D_{r,q_i}(f_i(F),F'_i) > s_i
\end{equation}
In order to get a contradiction out of~\pref{EqEquals_delta_i_s_i}
and~\pref{EqFar} we shall reapply Theorem~\ref{asymptotic flat theorem}  to a properly chosen sequence of intermediate basepoints,
near which $f_i(F)$ is still close to $F'_i$, but sufficiently far
that another Dehn twist flat, $F''_i$, is also close to it. The contradiction will then come from the following rigidity property of Dehn twist flats, which expresses that any Dehn twist flat spreads away from any other one in a linear fashion:

\begin{lemma}{asymmetric linear spreading}
There exists $\ep_1\in(0,1)$ such that for any sufficiently large~$r$, any
$x\in\Mod(S)$, and any Dehn twist flats $F_1$, $F_2$, if $F_1 \intersect \NN_{\frac{r}{2}}(x)  \ne \emptyset$ and if
$$
F_1 \intersect \NN_r(x) \subset \NN_{\ep_1 r}(F_2)
$$
then $F_1 = F_2$. 
\end{lemma}

\begin{proof}
Express $F_i$ as $\QQ (P_i)$ for a pants decomposition $P_i$ ($i=1,2$). 
Lemma~\ref{juncture} implies that the junctures of
$\QQ (P_1)$ and $\QQ (P_2)$ are 
$E_1=\QQ (P_1\rfloor P_2)$ and $E_2=\QQ (P_2\rfloor P_1)$. Note that 
$P_1\rfloor P_2$ is a marking with base equal to $P_1$, and a
transveral for each component of $P_1$ that is not a curve of
$P_2$. Hence, assuming $F_1 \ne F_2$, $E_1$ and $E_2$ must be
subflats of strictly smaller dimension. The pair $(F_i,E_i)$ is therefore
uniformly quasi-isometric to the pair $(\reals^\xi,\reals^k)$ for some $k<\xi = \xi(S)$. 

If $F_1$ meets $\NN_{\frac{r}{2}}(x)$ then pick $y\in F_1$ such that 
$d(x,y) \le \frac{r}{2}$. Lemma~\ref{juncture} implies that for any 
$z \in F_1$ the distance $d(z,F_2)$ is bounded below (up to uniform 
coarse-Lipschitz error) by the distance $d(z,E_1)$. Now use the 
elementary fact that for any $y\in\reals^\xi$ and any $R>0$ the 
Euclidean $R$-ball around~$y$ contains a point $z$ that is not 
contained in the Euclidean $\frac{R}{2}$-neighborhood of~$\reals^k$. After adjusting for the multiplicative errors, and setting $r$ large enough to overcome the additive errors, we find that for suitable $\ep'_1, \ep_1>0$ we have
$$F_1\intersect\NN_{\frac{r}{2}}(y)\not\subset \NN_{\ep'_1 r}(E_1)
$$
and so
$$ F_1\intersect\NN_{\frac{r}{2}}(y)\not\subset \NN_{\ep_1 r}(F_2).
$$
Since  $\NN_{\frac{r}{2}}(y)$ is contained
in $\NN_r(x)$, this is what we wanted to prove. 
\end{proof}

Let $\epsilon_2 = \min\{\frac{1}{R}, \frac{\epsilon_1}{2}\}$. Applying~\pref{EqEquals_delta_i_s_i} and~\pref{EqFar}, $\omega$-almost
surely the following two statements are true:

\begin{align*}
D_{Rs_i,p_0}(f_i(F),F'_i) &\le \epsilon_2 R s_i \\
D_{Rs_i,q_i}(f_i(F),F'_i) &> s_i \ge \epsilon_2 R s_i
\end{align*}

For each $i$, consider a sequence $(p_{i,k})$ starting at $p_0$ and ending with $q_i$, with step size $d(p_{i,k},p_{i,k+1})\le 1$. There must be some $j$ such that, labeling $x_i = p_{i,j}$ and $x'_i = p_{i,j+1}$, 
\begin{equation}\label{F' x good}
D_{Rs_i,x_i}(f_i(F),F'_i) \le \ep_2 R s_i
\end{equation}
but such that 
\begin{equation}\label{F' x' bad}
D_{Rs_i,x'_i}(f_i(F),F'_i) \ge \ep_2 R s_i
\end{equation}
Now assuming $\ep_2 R s_i > 1$ (which is true for large enough $i$),
we have $\NN_{Rs_i}(x'_i) \subset \NN_{Rs_i(1+\ep_2)}(x_i) $ and hence
(\ref{F' x' bad}) implies
\begin{equation}\label{F' no good}
D_{Rs_i(1+\ep_2),x_i}(f_i(F),F'_i) \ge \ep_2 R s_i
\end{equation}

Now we apply Theorem~\ref{asymptotic flat theorem} again, this time
using $(x_i)$ as the basepoints, and we conclude via
Lemma~\ref{LemmaUltraClose} that there exists a sequence $(F''_i)$ of
Dehn twist flats such that, for $\omega$-almost every $i$,
\begin{equation}\label{F'' good}
D_{Rs_i(1+\ep_2),x_i}(f_i(F),F''_i) < \ep_2 R s_i
\end{equation}
and in particular $F''_i\ne F'_i$ for $\omega$-a.e.\ $i$,  by (\ref{F' no good}). 
Now  (\ref{F'' good}) implies in particular that 
$$
f_i(F)\intersect
B_{Rs_i(1+\ep_2)}(x_i)\subset \NN_{\ep_2 R s_i}(F''_i).
$$
Moreover by (\ref{F' x good}) we have
$$
F'_i \intersect \NN_{Rs_i}(x_i) \subset \NN_{\ep_2 R s_i}(f_i(F))
$$
and moreover (by triangle inequality)
$$
F'_i \intersect \NN_{Rs_i}(x_i) 
\subset \NN_{\ep_2 R s_i}(f_i(F) \intersect \NN_{Rs_i(1+\ep_2)}(x_i)).
$$
Putting this together we see
\begin{equation}\label{flats too close}
F'_i \intersect \NN_{Rs_i}(x_i) \subset \NN_{2\ep_2 R s_i}(F''_i).
\end{equation}
Now since $x_i\in f_i(F_i)$ we note that $\NN_{\ep_2 R s_i}(x_i)$,
which is contained in $\NN_{\frac{R s_i}{2}}(x_i)$,
intersects $F'_i$ nontrivially. Now (\ref{flats too close}) 
implies, by Lemma \ref{asymmetric linear spreading} 
(noting $2\ep_2 \le \ep_1$), that 
$F'_i=F''_i$, a contradiction. 
\end{proof}

\subsection{Quasi-isometry classification}
\label{SectionQIClassificationProof}

We are now ready to prove Theorem \ref{QI classification}, the
classification of quasi-isometries of $\Mod(S) \approx \MM(S)$.  The
first step is to show that a $(K,C)$--quasi-isometry $f \from \MM(S) \to
\MM(S)$ induces an automorphism of the curve complex $\phi \from
\CC(S) \to \CC(S)$.

First we define $\phi$ on the vertex set $\CC_0(S)$. Given a vertex $c \in \CC_0(S)$, choose a Dehn twist 1-flat $\QQ(\mu)$ representing $c$, as defined in Section~\ref{SectionKFlat}. We claim that there exists a Dehn twist 1-flat $\QQ(\nu)$ at finite Hausdorff distance from $f(\QQ(\mu))$. From this claim, choose a vertex $d \in \CC_0(S)$ represented by $\QQ(\nu)$, define $\phi(c)=d$, and apply Lemma~\ref{LemmaFlatSetTh} \pref{ItemEqCEq} to conclude that $\phi(c)$ is well-defined. To prove the claim, first apply Lemma~\ref{LemmaFlatIntersection} to obtain two maximal Dehn twist flats $\QQ(\mu_0)$, $\QQ(\mu_1)$ whose coarse intersection is represented by $\QQ(\mu)$. Then apply Theorem~\ref{coarse flat theorem} to obtain Dehn twist flats $\QQ(\mu'_0)$, $\QQ(\mu'_1)$ at finite Hausdorff distance from $f(\QQ(\mu))$, $f(\QQ(\mu_1))$ respectively. The coarse intersection of $\QQ(\mu'_0)$ and $\QQ(\mu'_1)$ is represented by $f(\QQ(\mu))$, and also by $\QQ(\mu'_0 \rfloor \mu'_1)$ according to Lemma~\ref{LemmaFlatIntersection}. It follows that $\QQ(\mu'_0 \rfloor \mu'_1)$ is quasi-isometric to $f(\QQ(\mu))$ which is quasi-isometric to $\reals$, and so $\QQ(\mu'_0 \rfloor \mu'_1)$ is a Dehn twist 1-flat.

The proof of the claim in the previous paragraph also shows for any $k=0,\ldots,\xi$ that the image under $f$ of any Dehn twist $k$-flat is coarsely equivalent to a Dehn twist $k$-flat. The dimension $k$ only plays a role in the final sentence of the paragraph, in which $\reals$ is replaced by $\reals^k$.

To prove that $\phi$ is a bijection of $\CC_0(S)$, apply the same process to a coarse inverse $\bar f \from \MM(S) \to \MM(S)$ to obtain a map $\bar\phi \from \CC_0(S) \to \CC_0(S)$. The curve $c'' = \bar\phi(c')$ is represented by a Dehn twist 1-flat $Q(\mu'')$ at finite Hausdorff distance from $\bar f(Q(\mu'))$. The latter is at finite Hausdorff distance from $\bar f(f(Q(\mu))$, which is at finite Hausdorff distance from $Q(\mu)$ because $f$ and $\bar f$ are coarse inverses. The Dehn twist 1-flat $Q(\mu)$ therefore represents both $c$ and $c''$, and it follows that $c=c''$ by applying Lemma~\ref{LemmaFlatSetTh} \pref{ItemEqCEq}. This shows that $\bar\phi \phi$ is the identity, and a similar proof shows that $\phi \bar \phi$ is the identity.

We next show that two vertices $c_0,c_1 \in \CC(S)$ are endpoints of an edge of $\CC(S)$ if and only if $\phi(c_0)$, $\phi(c_1)$ are endpoints of an edge of $\CC(S)$. We need only prove the ``only if'' direction, the converse following by the same argument applied to a coarse inverse for $f$. Assuming $c_0,c_1$ are endpoints of an edge in $\CC(S)$, this edge is represented by some Dehn twist 2-flat $\QQ(\mu)$. The image $f(\QQ(\mu))$ is coarsely equivalent to some Dehn twist 2-flat $\QQ(\mu')$ which represents an edge with endpoints $d_0,d_1 \in \CC(S)$. We must show that $\{\phi(c_0),\phi(c_1)\} = \{d_0,d_1\}$. Let $\QQ(\mu_0)$, $\QQ(\mu_1)$ be Dehn twist 1-flats representing $c_0,c_1$, and let $\QQ(\mu'_0)$, $\QQ(\mu'_1)$ be Dehn twist 1-flats representing $\phi(c_0)$, $\phi(c_1)$. By Lemma~\ref{LemmaFlatSetTh} \pref{ItemContCCont} each of $\QQ(\mu_0)$, $\QQ(\mu_1)$ is coarsely contained in $\QQ(\mu)$, and so each of $\QQ(\mu'_0)$, $\QQ(\mu'_1)$ is coarsely contained in $\QQ(\mu')$, and by Lemma~\ref{LemmaFlatSetTh} \pref{ItemContCCont} again it follows that $\phi(c_0), \phi(c_1) \in \{d_0,d_1\}$. Since $c_0 \ne c_1$ it follows that $\QQ(\mu_0)$, $\QQ(\mu_1)$ are not coarsely equivalent, so $\QQ(\mu'_0)$, $\QQ(\mu'_1)$ are not coarsely equivalent, so $\phi(c_0) \ne \phi(c_1)$, from which it follows that $\{\phi(c_0),\phi(c_1)\} = \{d_0,d_1\}$.

This finishes the description of the automorphism $\phi \from \CC(S) \to \CC(S)$ induced by $f$.

Now suppose that $S \ne S_{1,2}$, returning to the case of $S_{1,2}$ 
in the end.  The theorem of Ivanov, Korkmaz, and Luo shows that $\phi$
is induced by some mapping class $\Phi \in \MCG(S)$ \cite{ivanov:complexes2,korkmaz:complex,luo:complex}. 
We must show
that the quantity
$$d(f,\Phi) = \inf_{x \in \MM(S)} d(f(x),\Phi(x))
$$ 
is bounded by a constant depending only on $K$, $C$, and the topology of $S$.

Consider a Dehn twist flat $Q(\mu)$. Theorem~\ref{coarse flat theorem} tells us that the Hausdorff distance between $f(Q(\mu))$ and some Dehn twist flat $Q(\mu')$ is bounded uniformly in terms of $K$, $C$, and the topology of $S$. We claim that $\mu'=\Phi(\mu)$. To prove this let $\mu = \{c_1,\ldots,c_\xi\}$, and let $c'_i = \Phi(c_i)$, and so $\Phi(\mu) = \{c'_1,\ldots,c'_\xi\}$. Let $Q(\mu_i)$ be Dehn twist 1-flats representing $c_i$, $i=1,\ldots,\xi$, and so $Q(\mu_i)$ is coarsely contained in $Q(\mu)$. The image $f(Q(\mu_i))$ has finite Hausdorff distance from some Dehn twist 1-flat $Q(\mu'_i)$ representing $c'_i$, and so $Q(\mu'_i)$ is coarsely contained in $Q(\mu')$ as well as in $Q(\Phi(\mu))$. It follows that $c'_i \in \mu'$, and so $\Phi(\mu) \subset \mu'$. But $\mu',\Phi(\mu)$ are sets of the same cardinality $\xi$, and so $\Phi(\mu)=\mu'$, proving the claim.

We now find a uniform bound to $d(f(x),\Phi(x))$ for each $x \in \MM(S)$. From Lemma~\ref{LemmaFlatIntersection} we obtain two Dehn twist flats $\QQ(\mu_0)$, $\QQ(\mu_1)$ whose coarse intersection is represented uniformly by the point $x$, with constants depending only on the topology of $S$. Applying $\Phi$, we obtain pants decompositions $\mu'_i = \Phi(\mu_i)$ and the coarse intersection of the corresponding Dehn twist flats $\QQ(\mu'_0)$, $\QQ(\mu'_1)$ is represented uniformly by the point $\Phi(x)$ in terms of the topology of $S$ only. Letting $\nu'_0 = \mu'_0 \rfloor \mu'_1$ and $\nu'_1 = \mu'_1 \rfloor \mu'_0$, the junctures $\QQ(\nu'_0)$, $\QQ(\nu'_1)$ also represent the coarse intersection, and so must be single points $\nu'_0$, $\nu'_1$ respectively. Choose points $y_i \in \QQ(\mu'_i)$ whose distance to $\Phi(x)$ is $\le C$ depending only on the topology of $S$, and so $d(y_0,y_1) \le 2C$. By Lemma~\ref{juncture} \pref{J divergence} we obtain a bound on $d(y_0,\nu'_0)$, and so also a bound on $d(\Phi(x),\nu'_0)$, depending only on the topology on $S$.

On the other hand, applying the map $f$, and using the bound on the Hausdorff distance between $f(\QQ(\mu_i))$ and $\QQ(\mu'_i)$, the point $f(x)$ also uniformly represents the coarse intersection of $\QQ(\mu'_0)$ and $\QQ(\mu'_1)$, the bounds now being in terms of $K$, $C$, and the topology of $S$. The same argument as in the last paragraph, using Lemma~\ref{juncture} \pref{J divergence}, now produces a bound on $d(f(x),\nu'_0)$ depending only on $K$, $C$, and the topology of $S$. 

Combining the last two paragraphs, we obtain the desired bound on $d(f(x),\Phi(x))$, completing the proof when $S \ne S_{1,2}$.

In the case of $S=S_{1,2}$, we still get an automorphism $\phi \from \CC(S_{1,2}) \to \CC(S_{1,2})$ as before, but it may no longer be induced by a mapping class of $S_{1,2}$. This is a finite-index problem caused by the hyperelliptic involution $\tau \from S_{1,2} \to S_{1,2}$ which we now recall. See Luo \cite{luo:complex}.

The map $\tau$ interchanges the two punctures of $S_{1,2}$. The quotient of $S_{1,2}$ by $\tau$, minus the four branch points, is $S_{0,5}$. Since $\tau$ is central in $\Mod(S_{1,2})$, every element descends to $S_{0,5}$ and we get a map $\beta\co\Mod(S_{1,2}) \to \Mod(S_{0,5})$ whose kernel is the center $\langle\tau\rangle \approx \Z/2\Z$. The image has index 5, because an element of $\Mod(S_{0,5})$ lifts if and only if it preserves the puncture which is the image of the two punctures of $S_{1,2}$. Hence $\beta$ is a quasi-isometry, and we let $\beta'\co\Mod(S_{0,5})\to\Mod(S_{1,2})$ be a
quasi-inverse. 

Now any quasi-isometry $f\co \Mod(S_{1,2})\to\Mod(S_{1,2})$ gives rise to a quasi-isometry $f' = \beta\circ f\circ\beta'$, and Theorem
\ref{QI classification} applied to $S_{0,5}$ gives an element $g\in 
\Mod(S_{0,5})$ such that $d(f',L_g)$ is bounded. If $g$ is in the index 5 image of $\Mod(S_{1,2})$ then a preimage $h\in\Mod(S_{1,2})$ works for $h$, i.e., $d(f,L_h)$ is bounded. If not, then at least we can produce the ``almost-geometric'' quasi-isometry $L = \beta' \circ L_g \circ \beta$, and obtain a bound on $d(f,L)$. This completes the proof.

\subsection{Quasi-isometric rigidity}
\label{QIR proof}
We conclude with the proof of Theorem \ref{QI rigidity},
quasi-isometric rigidity of $\Mod(S)$. The argument here is 
well-known, cf.\ \cite{schwartz:rankone}.

Let $G=\Mod(S)$. We may assume $\xi(S)\geq 2$ as the finite and virtually free
cases are already known. 
Left-multiplication gives a homomorphism $\lambda\co G\to \QI(G)$, where 
$\QI(G)$ is the
group of quasi-isometries of $G$ modulo the bounded-displacement subgroup. The
kernel of $\lambda$ is the center $Z=Z(G)$ (in general, $\ker \lambda$ consists of those
elements whose centralizer has finite index in $G$. For $\Mod(S)$ it 
is easy to show that the center are the only such elements).

Now supposing $S\ne S_{1,2}$, 
Theorem \ref{QI classification} implies that $\lambda$ is surjective. Hence we have
$\QI(G) = G/Z$ (this is Corollary \ref{QIGroup}).

Now if $\Gamma$ is quasi-isometric to $G$ then conjugation by the
quasi-isometry $\Phi$ gives an isomorphism $\QI(\Gamma) \homeo \QI(G)$ so we
get a map $\lambda'\co\Gamma \to \QI(G)$. 
Moreover $\ker \lambda'$ is finite: for each $\gamma\in \Gamma$, the quasi-isometry 
$\Phi L_\gamma \Phi^{-1}$ representing $\lambda'(\gamma)$ has uniformly bounded
constants (depending on $\Phi$), and hence by Theorem \ref{QI classification}
is a {\em uniformly} bounded distance from its approximating element of
$G$. Hence if $\gamma \in \ker \lambda'$, the approximating element is in $Z$, and
so $\gamma$ is restricted to a bounded set in $\Gamma$. Thus $\ker \lambda'$ is
finite. 

Finally, the image of $\lambda'$ has finite index in $\QI(G)$: this follows from the
fact that the left-action of $\Gamma$ on itself is transitive, and hence the
conjugated action on $G$ is cobounded. This gives the desired map
$\Gamma\to G/Z$ with finite kernel and finite-index image. 

If $G=\Mod(S_{1,2})$, we observe as in the proof of Theorem \ref{QI
classification} that $G/Z$ injects as a finite-index subgroup of
$G'=\Mod(S_{0,5})$, and hence it inherits the rigidity property from
$G'$ with the additional cost of restricting to a finite-index
subgroup of $\Gamma$.

\printindex

\bibliographystyle{hamsplain}

\bibliography{math}

\end{document}